\documentclass [12pt]{amsart}

\usepackage{amssymb,amsxtra,amsfonts}
\usepackage{epsf}              %
\usepackage{epsfig}             %
\usepackage{graphics}

\usepackage[colorlinks]{hyperref}

\usepackage{verbatim}  
\usepackage{bbding}   

\newenvironment{ppb}[1]
{\ \!\!\!\!\!\!\!\!\!\!\!\!\!\!\!\!\!\!\!\!\!\!\!\!\!\!\!\!\!\!\!\!\!\!\!\!\!\!\!\! {\bf PPB------------------------------------------------------------------------------------------------PPB}\newline \tiny {#1}
\  \newline\normalsize\phantom{f}\!\!\!\!\!\!\!\!\!\!\!\!\!\!\!\!\!\!\!\!\!\!\!\!\!\!\!\!\!\!\!\!\!\!\!\!\!\!\!\! {\bf PPB------------------------------------------------------------------------------------------------PPB}\newline}{}

\def\reE@DeclareMathSymbol#1#2#3#4{%
    \let#1=\undefined
    \DeclareMathSymbol{#1}{#2}{#3}{#4}}
\DeclareSymbolFont{symbolsC}{U}{txsyc}{m}{n}
\SetSymbolFont{symbolsC}{bold}{U}{txsyc}{bx}{n}
\DeclareFontSubstitution{U}{txsyc}{m}{n}
\reE@DeclareMathSymbol{\strictiff}{\mathrel}{symbolsC}{76}

\newcommand\beq{\begin{equation}}
\newcommand\eeq{\end{equation}}
\newcommand\bal{\begin{align*}}
\newcommand\eal{\end{align*}}   
\newcommand\bmx{\left(\begin{matrix}}
\newcommand\emx{\end{matrix}\right)}
\newcommand\bsmx{\left(\begin{smallmatrix}}
\newcommand\esmx{\end{smallmatrix}\right)}

\newcommand{\spq}{/\!\!/}

\providecommand{\spqa}[1]{\underset{#1}{/\!\!/}}
\newcommand{\st}{\ \bigl\vert\ }

\providecommand{\Rad}{\text{\rm Rad}}

\def\part#1{\frac{\partial\phantom{q}}{\partial#1}}

\newcommand {\flb}{\lbrack\!\lbrack}
\newcommand {\frb}{\rbrack\!\rbrack}
\newcommand {\flp}{(\!(}
\newcommand {\frp}{)\!)}

\newcommand{\union}{\cup} 
 
\newcommand{\sdp}{{\ltimes}}
\newcommand{\glu}{\strictiff}
\newcommand{\glue}[1]{\underset{#1}{\strictiff}}
\newcommand{\fus}{\circledast}
\newcommand{\fusion}[1]{\underset{#1}{\circledast}}

\newcommand{\Lie}{{\mathop{\rm Lie}}}

\newcommand{\Vect}{{\mathop{\rm Vect}}}             
\DeclareMathOperator{\Ann}{\mathop{\rm Ann}} 
\DeclareMathOperator{\ISto}{{\IS}to} 

\newcommand{\Gal}{\mathop{\rm Gal}}

\newcommand{\pap}[2]{{\ _{#1}\cA_{#2}}}

\newcommand{\papk}[3]{  \ _{#1}^{\phantom{#3}}\cA_{#2}^{#3}     }
\newcommand{\gah}{\pap{G}{H}}  
\newcommand{\gaho}{\papk{G}{H}{1}}   
\newcommand{\gahr}{\papk{G}{H}{r}}

\newcommand{\Ad}{{\mathop{\rm Ad}}}
\newcommand{\ad}{{\mathop{\rm ad}}}

\DeclareMathOperator{\pr}{pr}

\newcommand{\Prod}{\prod}

\newcommand{\tr}{{\mathop{\rm Tr}}}

\DeclareMathOperator{\Hom}{Hom}         
\DeclareMathOperator{\Aut}{\mathop{\rm Aut}}

\newcommand{\SL}{{\mathop{\rm SL}}}

\newcommand{\GL}{{\mathop{\rm GL}}}

\renewcommand{\Im}{\mathop{\rm Im}}

\renewcommand{\Re}{\mathop{\rm Re}}

\renewcommand{\ker}{\mathop{\rm Ker}}

\newcommand{\Stab}{{\mathop{\rm Stab}}}

\newcommand{\reg}{{\mathop{\footnotesize\rm reg}}}

\newcommand{\ba}{{\bf a}}

\newcommand{\bd}{{\bf d}}

\newcommand{\bD}{{\bf D}}

\newcommand{\bH}{{\bf H}}
\newcommand{\bG}{{\bf G}}

\newcommand{\bK}{{\text{\bf K}}}

\newcommand{\bM}{{\bf M}}

\newcommand{\bs}{{\bf S}}
\newcommand{\bS}{{\bf S}}

\newcommand{\bT}{{\bf T}}

\newcommand{\IA}{\mathbb{A}}
\newcommand{\IB}{\mathbb{B}}
\newcommand{\IC}{\mathbb{C}}
\newcommand{\ID}{\mathbb{D}}

\newcommand{\IH}{\mathbb{H}}

\newcommand{\IL}{\mathbb{L}}
\newcommand{\IM}{\mathbb{M}}

\newcommand{\IP}{\mathbb{P}}                                     
                           
\newcommand{\IR}{\mathbb{R}}                           
\newcommand{\IS}{\mathbb{S}}

\newcommand{\IT}{\mathbb{T}}

\newcommand{\IZ}{\mathbb{Z}}

\newcommand{\bP}{{\bf P}}

\newcommand{\A}{\mathcal{A}}
\newcommand{\cA}{\mathcal{A}}
\newcommand{\cB}{\mathcal{B}}
\newcommand{\cC}{\mathcal{C}}

\newcommand{\cD}{\mathcal{D}}
\newcommand{\cE}{\mathcal{E}}

\newcommand{\cH}{\mathcal{H}}
\newcommand{\ch}{\eta}     

\newcommand{\cK}{\mathcal{K}}

\newcommand{\cN}{\mathcal{N}}
\newcommand{\cO}{\mathcal{O}}

\newcommand{\cP}{\mathcal{P}}

\newcommand{\R}{\mathcal{R}}
\newcommand{\cR}{\mathcal{R}}

\newcommand{\cU}{\mathcal{U}}

\newcommand{\g}{       \mathfrak{g}     }

\newcommand{\lt}{\mathfrak{t}}
\newcommand{\lh}{\mathfrak{h}}

\newcommand{\h}{\mathfrak{h}}

\newcommand{\lu}{\mathfrak{u}}

\newcommand{\wt}{\widetilde}

\newcommand{\wh}{\widehat}

\newcommand{\al}{\alpha}

\newcommand{\be}{\beta}
\newcommand{\ga}{\gamma}
\newcommand{\de}{\delta}
\newcommand{\De}{\Delta}

\newcommand{\Ga}{\Gamma}

\newcommand{\la}{\lambda}
\newcommand{\La}{\Lambda}

\newcommand{\si}{\sigma}

\newcommand{\Si}{\Sigma}
\renewcommand{\th}{\theta}

\renewcommand{\bar}{\overline}

\makeatletter
 \newlength{\typesize}
\makeatother

\newlength{\vvoff}
\newlength{\hhoff}

\def\mapright#1{\smash{
        \mathop{\longrightarrow}\limits^{#1}}}

\def\mapdown#1{\Big\downarrow
        \rlap{$\vcenter{\hbox{$\scriptstyle#1$}}$}}

\def\underset#1#2{\ \smash{\mathop{ #2 }\limits_{#1}}\ }

\newcommand{\pf}{\begin{bpf}}

\newcommand{\pfms}{\begin{bpfms}}
\newcommand{\epf}{\end{bpf}\hfill$\square$\\}           
\newcommand{\epfms}{\end{bpfms}\hfill$\square$\\}       

\newcommand{\idea}{\begin{bidea}}

\newcommand{\eidea}{\end{bidea}\hfill$\square$\\}           

\newcommand{\sk}{\begin{bsk}}    

\newcommand{\esk}{\end{bsk}\hfill$\square$\\}           
\newcommand{\sketch}{\begin{bsketch}}

\newcommand{\esketch}{\end{bsketch}\hfill$\square$\\}

\theoremstyle{plain}  \newtheorem{hypo}{Hypothesis}[section]
\newtheorem{thm}[hypo]{Theorem}
\newtheorem{prop}[hypo]{Proposition}

\newtheorem{cor}[hypo]{Corollary}
\newtheorem{lem}[hypo]{Lemma}

\newtheorem {defn}[hypo]{Definition}

\theoremstyle{definition} \newtheorem{rmk}[hypo]{Remark}

\newtheorem{eg}[hypo]{Example}

\begin{document}

\title[Geometry and  braiding of Stokes data]{Geometry and  braiding of Stokes data; \\
{Fission and wild character varieties}}
\author{P. P. Boalch}

\dedicatory{To Robbie}

\begin{abstract}
A family of new algebraic Poisson varieties will be constructed, generalising the complex character varieties of Riemann surfaces.
Then the well-known (Poisson) mapping class group actions on the character varieties will be generalised.
\end{abstract}

\maketitle

\setcounter{tocdepth}{1}
\tableofcontents

\renewenvironment{ppb}[1]{}{}

\newpage
\section{Introduction}

Given a Riemann surface  $\wh\Si$ (with boundary) many people have studied moduli spaces
\beq\label{eq: pi1reps}
\Hom(\pi_1(\wh\Si), G)/G\eeq
 of representations of the fundamental group of $\wh\Si$ in a Lie group $G$, the {\em character varieties} (cf. \cite{sikora-charvars}).
If $G$ is a complex reductive group
with a chosen 
symmetric nondegenerate invariant
bilinear form
on its Lie algebra, then 
\eqref {eq: pi1reps} has an algebraic Poisson structure and the symplectic leaves are given by fixing the conjugacy classes around each component of the boundary. This provides a large class of holomorphic symplectic manifolds, which often have complete hyperk\"ahler metrics (and are then diffeomorphic to certain moduli spaces of meromorphic Higgs bundles 
\cite{Hit-sde, Sim-hboncc, Nak}).

Perhaps the best explanation as to {\em why} such spaces of fundamental group representations with fixed conjugacy classes have holomorphic symplectic structures is 
because they arise as symplectic quotients of the infinite dimensional affine space of all $C^\infty$ connections on a fixed $G$-bundle on $\wh\Si$
(this is the extension to surfaces with boundary, involving loop groups,
of the complexification of the viewpoint of Atiyah--Bott \cite{AB83}, described e.g. in \cite{Aud95long}).
Goldman \cite{Gol84} explained how this may also be understood in terms of the cup product in group cohomology.

The   quasi-Hamiltonian approach \cite{AMM} yields an alternative, algebraic construction of such symplectic manifolds as  finite-dimensional ``multiplicative'' symplectic quotients of a smooth affine variety, as follows.
Suppose $\wh\Si$ has $m\ge 1$ boundary circles $\partial_i$ and choose a basepoint $b_i\in \partial_i$ in each component.
Let $\Pi$ denote the fundamental groupoid of $\wh\Si$ with basepoints 
$\{b_1,\ldots,b_m\}$.
Then the space
$$\Hom(\Pi,G)$$ of homomorphisms from the groupoid $\Pi$ to the group $G$
is a smooth affine variety which naturally has the structure of quasi-Hamiltonian $G^m$-space, and so in particular has an action of $G^m$ and a group valued moment map
$$\mu:\Hom(\Pi,G)\to G^m.$$
The quotient $\Hom(\Pi,G)/G^m$ then inherits a Poisson structure and is isomorphic to \eqref{eq: pi1reps}.
Alternatively, if $\cC=(\cC_1,\ldots,\cC_m)\subset G^m$ is a conjugacy class (i.e. the choice of a conjugacy class $\cC_i\subset G$ for each boundary component), then the quasi-Hamiltonian reduction (the  multiplicative symplectic quotient) 
$$\mu^{-1}(\cC)/G^m$$
inherits a holomorphic  symplectic structure (where it is a manifold) and is isomorphic to a symplectic leaf of \eqref{eq: pi1reps}.

The first aim of this article is to give a similar algebraic approach
to more general spaces of monodromy-type data  
classifying irregular meromorphic connections on bundles on Riemann surfaces,
thereby  constructing many new algebraic symplectic manifolds
generalising the character varieties.
Secondly we will consider varying the initial data, leading to the braiding of the title,
generalising the much-studied mapping class group actions on
the character varieties.
 
To see the natural generalisation to the irregular case, recall first the relation between the above spaces and regular singular connections.
For this one starts with a smooth compact complex algebraic curve $\Si$
with $m$ distinct marked points $\ba=(a_1,\ldots,a_m)$ and
defines $\Si^\circ=\Si\setminus \{a_i\}$ to be the corresponding punctured curve.
Then, if $G=\GL_n(\IC)$,
Deligne's Riemann--Hilbert correspondence \cite{Del70} implies that
the  $G$ orbits in $\Hom(\pi_1(\Si^\circ), G)$
correspond bijectively to
isomorphism classes of  connections on rank $n$ 
algebraic vector bundles on $\Si^\circ$ with regular singularities
at each point $a_i$, and a similar statement holds for other groups $G$.
The condition of regular singularities means that the bundles have extensions across the punctures for which the connections have only simple poles (and local horizontal sections have at most polynomial growth as they approach the singularities).
Thus one sees a large class of generalisations may be obtained 
by relaxing this regularity assumption 
(and still lead to hyperk\"ahler manifolds \cite{wnabh} 
which again are often complete).
The irregular Riemann--Hilbert correspondence 
(on curves with $G=\GL_n(\IC)$) 
was worked out several decades ago (see \cite{malg-book})
but is not as well-known as its regular singular cousin: 
in brief the fundamental group representation is enriched by 
adding ``Stokes data'' at each singularity and there are various ways of thinking about this extra data: 
for instance as elements of a certain nonabelian 
cohomology space (Malgrange--Sibuya, cf. \cite{BV89}) or sheaf-theoretically (Deligne \cite{DMR-ci}), yielding an equivalence of categories.
The approach used here is closer to that of Martinet--Ramis \cite{MR91} describing the Stokes data as elements of certain unipotent subgroups of $G$
(in turn using Ecalle's work on multisummation); this viewpoint has the benefit of being as explicit as possible and amounts to having preferred cocycles representing each of the Malgrange--Sibuya cohomology classes (\cite{L-R94}).
The exact groupoid approach we use looks to be new however.
It is also useful when we vary the initial data.

\subsection{Statement of main results}

Fix a connected complex reductive group $G$ and a maximal torus $T\subset G$ with Lie algebras $\lt \subset \g$.
It is convenient to define  an ``irregular curve'' to be a smooth curve $\Si$ with marked points $\ba$ as above together with the extra data of an ``irregular type'' $Q_i$ at each marked point: if $z$ is a local coordinate on $\Si$ vanishing at $a_i$
then  
$$Q_i = \frac {A_{r_i}}{z^{r_i}} + \cdots \frac {A_1}{z}$$
for some  elements $A_i\in \lt$.
Given an irregular curve we will consider monodromy/Stokes data of connections on $G$-bundles which are locally isomorphic to 
$$dQ_i + \text{ less singular terms}$$
at each $a_i$, so that fundamental solutions involve essentially singular 
terms of the form $\exp(Q_i)$ near $a_i$.
(It is known that any meromorphic connection takes the above form after passing to a finite cover.)%

Then consider the real two-manifold with boundary
$$\wh \Si\to \Si$$  
obtained by taking the real oriented blow-up of $\Si$ at each marked point, thus replacing each point $a_i$ with a circle $\partial_i$.
The basic facts (see Definitions \ref{def: sing dirn}, \ref{def: stokes gp}) 
then are that $Q_i$ determines:

1) a connected complex reductive group $H_i\subset G$,
the stabiliser of $Q_i$,
 
2) a finite set $\IA_i\subset \partial_i$ of singular directions at $a_i$, 

\noindent
and for each singular direction $d\in \IA_i$ 

3) a unipotent group $\ISto_d(Q_i)\subset G$, normalised by $H_i$.

Then we puncture $\wh \Si$ once in its interior along each singular direction (sufficiently near the corresponding boundary component) and let 
$\wt \Si\subset \wh \Si$ denote the resulting punctured surface.
Choose a basepoint $b_i\in \partial_i$ in each boundary component and now let 
$\Pi$ denote the fundamental groupoid of $\wt \Si$ with basepoints 
$\{b_1,\ldots,b_m\}$. 
(If each irregular type 
$Q_i$ is zero there are no singular directions and we are in the original regular singular situation.)
Then consider the subset of {\em Stokes representations}
$$\Hom_\IS(\Pi,G)\subset \Hom(\Pi,G)$$
consisting of homomorphisms $\rho$ from $\Pi$ to $G$ which satisfy the following two conditions:
1) $\rho$ takes the loop $\partial_i$ based at $b_i$ into the group $H_i$, and 
2) for each singular direction $d\in \IA_i$, $\rho$ takes the small loop 
based at $b_i$ which goes around $\partial_i$ until direction $d$ before encircling the puncture in the direction $d$ and then retracing its path to $b_i$, to the group $\ISto_d$.
The main result is then:

\begin{thm}
The space $\Hom_\IS(\Pi,G)$ of Stokes representations is a smooth affine variety and is a quasi-Hamiltonian $\bH$-space, 
where $\bH=H_1\times \cdots \times H_m\subset G^m$.
\end{thm}

This implies that the quotient $\Hom_\IS(\Pi,G)/\bH$, which classifies
meromorphic connections with the given irregular types, inherits a Poisson structure, and its symplectic leaves are obtained by fixing a conjugacy class 
$\cC_i\subset H_i$ for each $i=1,\ldots,m$.
We will also characterise the stable points of $\Hom_\IS(\Pi,G)$ in the sense of geometric invariant theory (for the action of $\bH$),
show there are lots of examples when the quotients are well-behaved
and describe the irregular analogue of the Deligne--Simpson problem.

Such Poisson structures may be obtained analytically from an irregular analogue of the Atiyah--Bott approach (as in \cite{smid, wnabh}) and the quasi-Hamiltonian approach was worked out previously in \cite{saqh} 
in the case when the most singular coefficient
of each irregular type was regular (off of all the root hyperplanes).
The spaces of Stokes data are much more complicated when this assumption is removed and the main work of the present article is to develop an inductive approach to build (quasi-Hamiltonian) spaces of Stokes data out of simpler pieces.
This fits in well with the quasi-Hamiltonian philosophy of building moduli spaces of flat connections from simple pieces, and with the idea of 
the factorisation theorem of Ramis \cite{ramis-factn};
in effect we construct some new building blocks (higher fission spaces, \S\ref{sn: fission spaces})
and show that the spaces of Stokes representations may be built out of these using the quasi-Hamiltonian fusion and reduction operations.
That one can do such an induction at the quasi-Hamiltonian level is perhaps the main discovery of this article\footnote{\cite{fission} discussed the possibility of ``fusion on the other side of the analytic halo''---in effect here we do fusion {\em within} the analytic halo.}.
Surprisingly it follows immediately that our 
building blocks may be used to construct many 
other holomorphic symplectic manifolds, beyond the quotients of the spaces of Stokes representations.
For example one may glue various surfaces $\wt \Si$ along their boundaries provided the groups $H_i$ match up.
Also one may obtain Van den Bergh's quasi-Hamiltonian spaces \cite{vdb-doublepoisson, vdb-ncqh} from the higher fission spaces and thus all of the so-called multiplicative quiver varieties.

\subsection{Varying the initial data}

In brief the above result implies that the choice of an irregular curve canonically determines a Poisson variety 
$\Hom_\IS(\Pi,G)/\bH$.
For the second main result
we will
define in \S\ref{sn: adm def} the notion of an 
``admissible family'' of irregular curves over a base $\IB$,  
generalising the notion of 
deforming a smooth curve with marked points 
such that the curve remains smooth and none of the points coalesce.
Then we will show (Theorem \ref{thm: loc syst of vars}) 
that the corresponding family of 
Poisson varieties assemble into a ``local system of Poisson varieties'' over $\IB$, i.e. into a nonlinear fibration with a flat (Ehresmann) connection which  integrates along any path in the base to yield algebraic Poisson isomorphisms between fibres.
This leads to an algebraic Poisson action of the fundamental 
group of $\IB$ on any fibre $\Hom_\IS(\Pi,G)/\bH$.
Such deformations have been considered briefly by Witten \cite{Witten-wild} \S6 in relation to geometric 
Langlands and $\cN=4$ super Yang-Mills theory.
Mathematically these flat nonlinear connections should be viewed as the irregular analogue of the Gauss--Manin connection on nonabelian cohomology
 (\cite{smid} \S7).

Basic examples of admissible deformations in the irregular case were considered by Jimbo et al \cite{JMU81};
they looked at the case $G=\GL_n(\IC)$ when the most singular coefficient at each pole had distinct eigenvalues. This was extended to other reductive groups in \cite{bafi} (keeping the most singular coefficient at each pole off all the root hyperplanes) and it was shown that in the simplest case the resulting Poisson action of the $G$-braid group coincides with the quasi-classical limit of the quantum-Weyl group action of Lusztig, Soibelman, Kirillov--Reshetikhin\footnote{This quasi-classical limit, an action of the $G$-braid group on the dual Poisson Lie group, was explicitly computed by 
De\! Concini--Kac--Procesi \cite{DKP}.
A key point of the geometrical approach is the identification \cite{smapg, bafi} of a simple space of Stokes data with the Poisson Lie group dual to $G$.}.
One impetus 
(cf. \cite{smapg, smid}) for this line of thinking came from trying to understand  the Poisson braid group actions in Dubrovin's work on semisimple Frobenius manifolds, related to the Markoff polynomial in the 3-dimensional case (\cite{Dub95long} p.243)---see also \cite{cec-vafa-nequals2classn, Hit95long, Ugag, bondal-gpoid.braids, xu-poisson.inv}.
More generally,
in the integrable systems literature interest in symplectic structures on spaces of 
Stokes data goes back at least to Flaschka--Newell 
\cite{FN82}; the more recent work of  Woodhouse \cite{Woodh-sat} and  Krichever \cite{Krich-imds}
 also computed such 
symplectic forms explicitly (in the $\GL_n(\IC)$ case with distinct leading eigenvalues) and our more general formulae were found similarly. The quasi-Hamiltonian approach here yields an algebraic proof that such two-forms are indeed  symplectic, and in the general linear case a quite simple proof that  for generic parameters the wild character varieties are smooth symplectic algebraic varieties (see Corollary \ref{cor: good quotients}). This last result alone probably justifies the quasi-Hamiltonian approach and was one of our main aims.

Note that if for example one is interested in complete hyperk\"ahler manifolds it makes little sense to restrict attention to the regular singular case: there are many examples of isomorphisms between the moduli spaces 
that arise in the irregular case and the regular singular case (some examples appear in \cite{szabo-nahm} and in such cases the hyperk\"ahler metrics match up), but it is not true that all irregular cases are isomorphic to a regular singular case---on the contrary it seems that if one counts in each dimension the number of deformation classes of complete hyperk\"ahler manifolds arising from Hitchin's self-duality equations, then most classes only have irregular representatives.

The layout of this article is as follows. 
Sections 2--6 are basically pure quasi-Hamiltonian geometry, first giving the background definitions, then directly establishing the new spaces we will need, then establishing many properties of them. 
Section 7 defines the spaces of Stokes data corresponding to 
connections on a disc and shows that, with suitable framings, they are quasi-Hamiltonian.
This is used in Section 8 to prove that the space of Stokes/monodromy data attached to a (global) irregular curve is indeed an algebraic Poisson variety.
Section 9 discusses various aspects of stability (one related to differential Galois theory), defines the irregular  Deligne--Simpson problem
and gives many examples when the quotients are well-behaved.
Finally Section 10 considers admissible families of irregular curves 
and shows that the corresponding family of Poisson varieties fit together into a Poisson local system.
It also mentions the link to Baker functions and integrable hierarchies.
(Some aspects of the irregular Riemann-Hilbert correspondence are discussed in Appendix A, to help motivate the basic definitions.)
Note that we have focused on the new features that occur in the present context and  some of the results from \cite{smid, bafi, saqh} whose generalisation is routine have been omitted.

\noindent{\bf Acknowledgments.}
This work was partially supported by ANR grants 
08-BLAN-0317-01/02 (SEDIGA), 09-JCJC-0102-01 (RepRed).
The author is grateful to the referee for several very helpful suggestions.

\section{Quasi-Hamiltonian geometry}

Some familiarity with quasi-Hamiltonian geometry of Alekseev--Malkin--Meinrenken \cite{AMM} will be assumed. This section will recall (the holomorphic analogue of) the basic results.
In essence this theory is a multiplicative version of the usual Hamiltonian theory, with moment maps taking values in Lie groups rather than the dual of the Lie algebra.
The axioms for the analogue of the symplectic form and its interaction with the group action and the moment map are  more complicated.
The upshot is a direct and explicit algebraic approach to constructing certain quite exotic symplectic manifolds, 
previously constructed via infinite dimensional techniques.
The motivation in \cite{AMM} was to give a finite dimensional algebraic construction of the symplectic structure on moduli spaces of flat connections on bundles over Riemann surfaces with fixed local monodromy conjugacy classes. 
We have found this theory is also useful to construct new moduli spaces.

\subsection{Notation}
Let $G$ be a connected complex reductive 
group\footnote{The main results are new even for  $G=\GL_n(\IC)$ so the reader could restrict to that case, but one gets a richer class of braid group actions in general (see the examples in \S\ref{ssn: example base} and \cite{bafi} which motivated us to define $G$-valued Stokes multipliers).}
 with Lie algebra $\g$.
(The group $G(\IC)$ of $\IC$ points, will often also be denoted by $G$.) 
Suppose
we have chosen a symmetric nondegenerate invariant
bilinear form $(\ ,\ ):\g\otimes\g\to \IC$ (this choice will be tacitly assumed throughout). 
The Maurer--Cartan forms on $G$ are denoted $\theta,\overline\theta\in\Omega^1(G,\g)$ 
respectively (so in any representation 
$\theta=g^{-1}dg, \overline\theta=(dg)g^{-1}$).
Generally if $\A,\cB,\cC\in\Omega^1(M,\g)$ are $\g$-valued
holomorphic one-forms on a complex manifold $M$ then
$(\A,\cB)\in\Omega^2(M)$ and 
$[\A,\cB]\in\Omega^2(M,\g)$ are defined
by wedging the form parts and pairing/bracketing the Lie algebra parts.
Define $\cA^2:= \frac{1}{2}[\A,\A]\in \Omega^2(M,\g)$ (which works out
correctly in
any representation of $G$  using matrix multiplication).
Then one has 
$d\theta=-\theta^2, d\overline\theta = \overline\theta^2$.
Define $(\A\cB\cC) = (\A, [\cB,\cC])/2\in\Omega^3(M)$ (which is
invariant under all permutations of $\A, \cB, \cC$). The canonical 
bi-invariant three-form on $G$ is then
$\frac{1}{6}(\theta^3)$.
The adjoint action of $G$ on $\g$ will be 
denoted $gXg^{-1}:=\Ad_gX$ for any
$X\in\g, g\in G$.
If $G$ acts on  $M$, the fundamental vector field $v_X$ of $X\in\g$ is minus
the tangent to the flow 
$(v_X)_m = -\frac{d}{dt} (e^{Xt}\cdot m)\bigl\vert_{t=0}$ so that the map
$\g\to\Vect_M; X\to v_X$ is a Lie algebra homomorphism. 
(This sign convention differs from \cite{AMM} 
leading to sign changes in the quasi-Hamiltonian axioms and
the fusion and equivalence theorems.)

Recall that a complex manifold $M$ is a 
{\em complex quasi-Hamiltonian $G$-space}
if there is an action of $G$ on $M$, 
a $G$-equivariant map $\mu:M\to G$ (where $G$ acts on itself by
conjugation) and a $G$-invariant holomorphic two-form
$\omega\in \Omega^2(M)$ such that:

\noindent(QH1). 
The exterior derivative of $\omega$ is the pullback along the moment map of the canonical three-form on $G$:
$d\omega = \mu^*(\th^3)/6.$

\noindent(QH2).
For all $X\in \g$,
$\omega(v_X,\cdot\,) = \frac{1}{2}\mu^*(\theta+\overline\theta, X)
\in \Omega^1(M).$

\noindent(QH3).
At each point $m\in M$: 
$\ker \omega_m \cap \ker  d\mu = \{0\} \subset T_mM$. 

It is possible to show that if (QH1) and (QH2) hold, then (QH3) is equivalent to the condition 
$\ker(\omega_m)=\left\{(v_X)_m\ \bigl\vert\ X\in\g \text{ satisfies } \Ad_g X = -X \text{ where } g:=\mu(m)\right\}$ 
(cf. \cite{ABM-purespinors} Remark 5.3).

\begin{rmk}\label{rmk: ab implies symp}
Observe that if $G$ is abelian 
(and in particular if $G=\{1\}$ is trivial) then these axioms imply that the two-form $\omega$ is a complex symplectic form.
The reduction procedure (see below) yields many symplectic manifolds in this way.
\end{rmk}

\begin{eg}[Conjugacy classes \cite{AMM}]
Let $\cC\subset G$ be a conjugacy class, with the conjugation action
of $G$ and moment map $\mu$ given by the inclusion map.
Then $\cC$ is a quasi-Hamiltonian $G$-space with two-form
$\omega$ determined by 
$$\omega_g(v_X,v_Y) = 
\frac{1}{2}\bigl((X,gYg^{-1}) - (Y,gXg^{-1})\bigr) $$
for any $X,Y\in \g, g\in \cC$.
\end{eg}

Other basic examples of quasi-Hamiltonian spaces appear as moduli spaces of holomorphic connections on Riemann surfaces with boundary, with a framing at one point on each boundary component.
(If one instead chooses a framing on all of the boundary then one obtains an infinite dimension symplectic manifold with a Hamiltonian loop group action changing the
 framing---one way to understand the quasi-Hamiltonian axioms is via the way in which such loop group spaces yield finite dimensional quasi-Hamiltonian spaces in \cite{AMM} \S8, by forgetting the framing at all but one point on each component the boundary.)
For example the annulus has two boundary components so corresponds naturally to a quasi-Hamiltonian $G\times G$-space,
which may be written explicitly as follows.

\begin{eg}[The double \cite{AMM}]
The space $\bD= G\times G$ 
is a quasi-Hamiltonian $G\times G$-space 
with $(g,k)\in G\times G$ acting as 
$(g,k)(C,h)=(kCg^{-1}, khk^{-1})$, with moment map 
$$\mu(C,h)= (C^{-1}hC,h^{-1})\in G\times G$$
and with two-form $\omega$ such that
\beq  \label{eq: double omeg def}
2 \omega = 
(\bar\ga,\Ad_h\bar\ga) 
+ (\bar\ga,\bar\ch + \ch)
\eeq
where $\bar\ga = C^*(\bar\th), \ch = h^*(\th), \bar\ch = h^*(\bar\th).$
\end{eg}

The notation $C^*(\bar\th)$ here means that we view $C$ as a map from $\bD$ to $G$ and pull back the right-invariant Maurer-Cartan form $\bar\th$ to obtain a $\g$-valued one-form on $\bD$. 
Similarly the one-holed torus leads to the following space (the name will be explained below).

\begin{eg}[Internally fused double \cite{AMM}]
The space $\ID= G\times G$ is a quasi-Hamiltonian $G$-space 
with $G$ acting by diagonal conjugation 
($g(a,b)=(gag^{-1}, gbg^{-1})$), moment map given by the group commutator
$$\mu(a,b)= aba^{-1}b^{-1}$$
and two-form
$$\omega_\ID=
-\frac{1}{2}(a^*\theta,b^*\overline\theta)
-\frac{1}{2}(a^*\overline\theta,b^*\theta)
-\frac{1}{2}((ab)^*\theta,(a^{-1}b^{-1})^*\overline\theta).$$
\end{eg}

These are both special cases of the following result. Let $\Si$ be a compact connected Riemann surface with boundary. Let $g$ be the genus of $\Si$ and let  
$m$ be the number of boundary components (we assume $m\ge 1$).
Choose a basepoint $b_i$ on the $i$th boundary component for each $i=1,\ldots, m$.
Let $$\Pi=\Pi_1(\Si,\{b_1,\ldots,b_m\})$$
be the fundamental groupoid of $\Si$ with basepoints $\{b_i\}$, i.e. the groupoid of homotopy classes of paths in $\Si$ whose endpoints are in the set of chosen basepoints.

\begin{thm}\label{thm: amm gpoid spaces}
The space $\Hom(\Pi,G)$ of homomorphisms from the groupoid $\Pi$ to the group $G$ is a smooth quasi-Hamiltonian $G^m$-space.
\end{thm}
\pf
This is just a slightly more intrinsic restatement of \cite{AMM} Theorem 9.1. Choosing suitable paths generating $\Pi$ identifies
$\Hom(\Pi,G)$  with $G^{2(g+m-1)}$ and in turn with the reduction of  the fusion product (see below) $\ID^{\fus g}\fus \bD^{\fus m}$ by the diagonal action of $G$ at the identity value of the moment map. Then one can check the result is independent of the chosen generating paths.
\epf

If $\Si$ is disconnected (and each component has at least one boundary component) the same result holds, taking the product of the quasi-Hamiltonian spaces attached to each connected component.

\subsection{Operations}

The fusion product, which puts a ring structure on the  
category of quasi-Hamiltonian $G$-spaces, is defined as follows. 

\begin{thm}[\cite{AMM}]\label{thm: fusion}
Let $M$ be a quasi-Hamiltonian $G\times G\times H$-space, 
with moment map $\mu=(\mu_1,\mu_2,\mu_3)$. 
Let $G\times H$ act by 
the diagonal embedding $(g,h)\to (g,g,h)$. 
Then $M$ with two-form 
\begin{equation} \label{eqn: fusion 2form}
\wt{\omega}= \omega - \frac{1}{2}(\mu_1^* \theta, \mu_2^* \overline \theta)
\end{equation}
and moment map
$$\wt{\mu} = (\mu_1\cdot \mu_2,\mu_3):M\to G\times H$$
is a quasi-Hamiltonian $G\times H$-space. 

\end{thm}

We will refer to the extra term subtracted off in \eqref{eqn: fusion 2form} as
the ``fusion term''. 
If $M_i$ is a quasi-Hamiltonian $G\times H_i$ space for $i=1,2$ their 
fusion product $$M_1\fus M_2$$ 
is defined to be the
quasi-Hamiltonian $G\times H_1\times H_2$-space 
obtained from the
quasi-Hamiltonian $G\times G \times H_1 \times H_2$-space 
$M_1\times M_2$ by
fusing the two factors of $G$.

This is set up so that it corresponds to gluing two boundary components  into two of the holes of a three-holed sphere. More precisely, suppose 
$\Si_1$ is a (possibly disconnected) surface with boundary and we choose two distinct boundary components and construct a new surface $\Si_2$ by gluing the two chosen boundary 
components of $\Si_1$ into two of the holes of a three-holed sphere.
Let $M_i$ be the quasi-Hamiltonian space attached to $\Si_i$ via Theorem \ref{thm: amm gpoid spaces} (repeating for each connected component if necessary) for 
$i=1,2$. Then $M_2$ is isomorphic to the space obtained by fusing the two $G$ factors of $M_1$ corresponding to the two chosen boundary components. (For example, if $\Si_1$ is the annulus, this explains the name ``internally fused double''.)

Now let us recall the quasi-Hamiltonian reduction theorem:

\begin{thm}[\cite{AMM}]
Let $M$ be a quasi-Hamiltonian $G\times H$-space with moment map
$(\mu,\mu_H):M\to G\times H$ and suppose
that the quotient by $G$ of the inverse image $\mu^{-1}(1)$  
of the identity under the first moment map is a manifold.
Then the restriction of the two-form $\omega$ to $\mu^{-1}(1)$ 
descends to the {\em reduced space}
\begin{equation}\label{eq: qh quot}
M\spq G := \mu^{-1}(1)/G 
\end{equation}
and makes it into a quasi-Hamiltonian $H$-space. 
In particular, if $H$ is abelian (or in particular trivial) 
then $M\spq G$ is a complex symplectic manifold. 
\end{thm}

Next we recall that the quotient of a quasi-Hamiltonian $G$-space by $G$ is Poisson. This result will be used in the following form:

\begin{prop}\label{prop: poisson quotients}
Suppose $M$ is a smooth affine variety with the structure of quasi-Hamiltonian $G$-space.
Then the (geometric invariant theory) quotient $M/G$ is a Poisson variety.
\end{prop}
\pf
It is well-known (\cite{AMM} 4.6) that the ring of $G$-invariant functions on $M$ is a Poisson algebra (see also \cite{AKM} \S6, \cite{ABM-purespinors} \S5.4). 
But by definition the geometric invariant theory quotient is the affine variety associated to the 
ring of $G$-invariant functions on $M$, so is Poisson.
\epf

Note that the points of the  geometric invariant theory quotient correspond bijectively to the {\em closed} $G$-orbits in $M$, and so in general it is different to the set-theoretic quotient. Alternatively one may view the 
points of the  geometric invariant theory quotient as parameterising the quotient of $M$ by a stronger equivalence relation than orbit equivalence (S-equivalence): two points of $M$ are S-equivalent if their orbit closures intersect.
(One may also consider other geometric invariant theory quotients, by using a nontrivial linearisation, but we will ignore these here for brevity.)
Note that unless otherwise stated $M/G$ will denote the geometric invariant theory quotient, and $M\spq G$ will denote the complex quasi-Hamiltonian quotient \eqref{eq: qh quot} (i.e. the geometric invariant theory quotient of the  subvariety $\mu^{-1}(1)\subset M$).

\subsection{Gluing $M\glue{} N$}
It is convenient to formalise the (well-known) notion of gluing 
quasi-Hamiltonian spaces, as follows. 
Given a quasi-Hamiltonian $G\times G\times H$-space $M$, we may fuse the two $G$ factors to obtain a quasi-Hamiltonian $G\times H$ space. Then, if the quotient is well-defined, we may
reduce by the $G$ factor (at the identity of $G$) to obtain a quasi-Hamiltonian $H$-space, the {\em gluing} of the two $G$-factors.
Thus for example if $M_i$ is a  quasi-Hamiltonian $G\times H_i$-space for $i=1,2$ then $M_1$ and $M_2$ may be glued to obtain a quasi-Hamiltonian $H_1\times H_2$ space (if it is a manifold) by gluing %
 their product:
$$M_1\glue{G} M_2 := (M_1 \fusion{G} M_2)\spq G.$$

If the factors being glued are clear from the context this will be abbreviated to $M_1\glu M_2$.
In most of the cases we will consider, the $G$-action will be free with a global slice so there is no problem performing the gluing.
Note that whereas fusion is only commutative up to isomorphism, the gluing operation is actually commutative (when it is defined).

\subsection{Van den Bergh's quasi-Hamiltonian spaces $\cB(V,W)$}

Choose two finite dimensional complex vector spaces $V,W$ and define
\beq\label{eq: vdb space}
\cB(V,W)=\{(a,b)\in \Hom(W,V)\oplus\Hom(V,W)\st \det(1+ab)\ne 0\}.
\eeq
The automorphism groups $\GL(V)$ and  $\GL(W)$ of $V$ and $W$ naturally induce an action of $\GL(V)\times\GL(W)$ on $\cB(V,W)$.
\begin{thm}(\cite{vdb-doublepoisson, vdb-ncqh, yamakawa-mpa})
$\cB(V,W)$ is a quasi-Hamiltonian $\GL(V)\times\GL(W)$-space.
The moment map is
\beq\label{eq: vdb mmap}
(a,b)\mapsto ((1+ab)^{-1}, 1+ba)\in \GL(V)\times\GL(W)
\eeq
and the two-form is
\beq\label{eq: vdb 2form}
\omega = \frac{1}{2}\left(
\tr_V(1+ab)^{-1}da\wedge db  -
\tr_W(1+ba)^{-1}db\wedge da \right).
\eeq
\end{thm}
An alternative proof will follow from 
Theorems \ref{thm: qh 1 level} and \ref{thm: isom to vdb} below.
Henceforth the notation ``$(a,b)\in \cB(V,W)$" will always mean
$a\in \Hom(W,V)$ and $b\in\Hom(V,W)$ with $1+ab$ invertible,
analogously to the  convention that $(p,q)\in T^*\Hom(V,W)$ means  $q\in\Hom(V,W)$ and
 $p\in \Hom(W,V) = \Hom(V,W)^*$.
Note that such spaces (without the quasi-Hamiltonian structure) are familiar from the explicit description of the local classification of regular holonomic $\cD$-modules on a curve 
(see e.g. \cite{bdm-disk, maisonobe-couple, malg-book}). %

\section{Higher fission spaces} \label{sn: fission spaces}

In this section we will describe some new algebraic quasi-Hamiltonian spaces.
Later
it will be explained how these spaces arise from considering the explicit local classification (in terms of Stokes data) of some simple irregular connections on curves, and how  %
they may be glued together to construct some much more complicated spaces of Stokes data classifying irregular meromorphic connections on curves.

Let $G$ be a connected complex reductive group and choose
a parabolic subgroup $P_+\subset G$ and a Levi subgroup $H\subset P_+$.
Let $P_-\subset G$ be the (unique) opposite parabolic with Levi subgroup $H\subset P_-$ so that $P_-P_+$ is dense in $G$ (see e.g. \cite{Bor91} p.199).
Let $U_\pm\subset P_\pm$ be the corresponding unipotent radicals.
For example
if $G$ is a general linear group so that $G=\GL(V)$ for a finite dimensional complex vector space $V$,
then choosing a parabolic subgroup of $G$ together with a Levi subgroup is equivalent to choosing an ordered grading of $V$, i.e. 
to choosing a direct sum decomposition 
$$V=\bigoplus_1^k V_i$$
for some integer $k\ge 1$. By convention we then take the subgroup 
$P_+ \subset G$ to be the subgroup stabilising the flag
$$F_1\subset F_2\subset \cdots \subset F_k=V$$
where  $F_i= V_1\oplus\cdots \oplus V_i$, and $H=\Prod \GL(V_i)\subset G$.
Thus in an adapted basis $U_+$ is the subgroup of block upper triangular matrices with $1$'s on the diagonal and $H$ is the block diagonal subgroup.

Now fix an integer $r\ge 1$ and define a space
$$\gahr := G \times (U_+\times U_-)^{r}\times H.$$
In the case $r=1$ the superscript will be omitted so that $\gah=\gaho$.

A point of $\gahr$ is given by specifying $C\in G, h\in H$ and 
$\bS\in (U_+\times U_-)^{r}$ with $\bS = (S_1,\ldots,S_{2r})$ where
$S_{\text{even}}\in U_-$ and $S_{\text{odd}}\in U_+.$
There is an action of $G\times H$ on $\gahr$ given by 
$$(g,k)(C,\bS,h)  = (kCg^{-1}, k\bS k^{-1}, khk^{-1})$$
where $(g,k)\in G\times H$ and $k\bS k^{-1} = (kS_1k^{-1},\ldots,kS_{2r}k^{-1})$.

Note that the only choices involved here are $G, P_+$ and the Levi subgroup $H$ (i.e. a lifting to $P_+$ of the Levi factor $P_+/U_+\cong H$). 
Sometimes it will be convenient to denote the same space also by 
$\!\!\!\papk{H}{G}{r}$.
In the general linear case, 
an ordered graded vector space $V$ 
thus determines a space $\gahr$ for each integer $r\ge 1$,
and we will denote them by 
$$\cA^r(V_1,\ldots,V_k) = \gahr,$$
or by $\cA^r(V)$ if the choice of ordered grading of $V$ is clear (and as before $\cA(V) = \cA^1(V)$).
The main result we will prove in this section is the following.
\begin{thm}\label{thm: qh 1 level}
Suppose $G$ is a complex reductive group 
and $H\subset P_+\subset G$ are chosen as above. 
Then $\gahr$ is a quasi-Hamiltonian $G\times H$-space, with moment map
$$\mu(C,\bs,h) = (C^{-1} h S_{2r}\cdots S_2 S_1 C,\,\, h^{-1}) \in G\times H.$$
\end{thm}

In general these spaces $\gahr$ will be referred to as 
``{\em higher fission spaces}'' 
(enabling one to break the group from $G$  to $H$). 

\noindent
{\bf Examples.}\ \ 

1) If $P_+=G=H$ then $\gahr$ is the double $G\times G$ of \cite{AMM}.

2) When $P_+$ is a Borel subgroup, so that $H$ is a maximal torus, we may pass to a covering to obtain the spaces $\wt \cC = G\times (U_+\times U_-)^r\times \lh$ of \cite{saqh}, where $\lh = \Lie(H)$.

3) When $r=1$ (and $P_+$ arbitrary) the spaces $\gahr$ specialise to the spaces $\gah$ of \cite{fission}.
They always have dimension $2\dim(G)$ and (up to passing to a covering) interpolate between the double $G\times G$ (which appears when $H=G$) and $G\times G^*$ where $G^*$ is the dual Poisson Lie group of $G$ (which appears when $H$ is a maximal torus).

4) Next suppose $r=2$ and we perform the reduction by $G$ of $\gahr$ at the value $1$ of the moment map.
The resulting space is 
\beq\label{eq: r=2 reduced space}
\left\{(S_1,S_2,S_3,S_4,h)\in U_+\times U_-\times  U_+\times U_- \times H \st  hS_4S_3S_2S_1 =1 \right\}.\eeq
This inherits the structure of  quasi-Hamiltonian $H$-space
with moment map $h^{-1}$. (It is symplectic if $H$ is a maximal torus as in \cite{saqh}).
By forgetting $h,S_3,S_4$ this space  \eqref{eq: r=2 reduced space} embeds as the (open) subset
of $ U_+\times U_-$ consisting of pairs $(S_1,S_2)$ 
such that $S_2S_1\in P_+P_-\subset G$.
In \S\ref{sn: derive vdb} Theorem \ref{thm: isom to vdb}
we will establish the following statement:

{\em If $G$ is a general linear group and $P_+$ is a maximal proper parabolic subgroup then the spaces \eqref{eq: r=2 reduced space} coincide with the Van den Bergh spaces \eqref{eq: vdb space}.}

One of the main features that appears in the general parabolic case is the possibility to glue such fission spaces end to end. 
If one does this for a decreasing sequence of nested Levi subgroups (and decreasing exponents $r$) then, as will be shown in \S\ref{sn: local Stokes}, 
all the complicated spaces which arise in this way actually appear as spaces of (framed) Stokes data for
meromorphic connections  on curves.

\pfms (of Theorem \ref{thm: qh 1 level}).
Define maps $C_i: \gahr \to G$ by
$$C_i = S_i\cdots S_2S_1C$$
so that $C=C_0$.
Define $b = hS_{2r}\cdots S_2S_1: \gahr \to G$  so that
the $G$ component of the moment map $\mu$  is $C^{-1} b C$.
This enables us to define the following $\g$-valued one-forms on $\gahr$:
$$\ga_i = C_i^*(\th),\qquad\bar\ga_i = C_i^*(\bar\th),\qquad
\eta = h^*(\th_H),\qquad
\bar\be = b^*(\bar\th)$$
where $\th, \bar \th$ are the Maurer--Cartan forms on $G$ (and 
$\th_H, \bar \th_H$ are the Maurer--Cartan forms on $H$).
We may then define a two-form $\omega$ on $\gahr$ by the formula 
\beq  \label{eq: omeg def}
2 \omega = 
(\bar\ga,\Ad_b\bar\ga) 
+ (\bar\ga,\bar\be)
+ (\bar\ga_m, \ch) -
\sum_{i=1}^m(\ga_{i},\ga_{i-1})
\eeq
where $m=2r$ and $\bar\ga = \bar\ga_0$, and the brackets $(\ ,\ )$ 
denote the bilinear form on $\g$.
We will show $\gahr$ is a quasi-Hamiltonian $G\times H$-space with this two-form.
First, since we will use it often, note that the invariance of the bilinear form implies
the pairing $(\ ,\ ):\g\otimes\g\to \IC$ 
restricts to zero on 
$\lu_\pm\otimes (\lh\oplus \lu_\pm)$ and is nondegenerate on 
$\lh\otimes \lh$, where $\lh=\Lie(H), \lu_\pm=\Lie(U_\pm)$.
\ppb{
\end{lem}
\pf
Choose a Cartan subalgebra $\lt\subset\lh\subset \g$.
Let $\cR\subset \lt^*$ be the roots of $\g$,
let $\g=\lt\oplus\bigoplus_{\al\in \cR}\g_\al$ be the root space decomposition, and write $\g_0=\lt$.
Then observe that if $X\in \g_\al,Y\in \g_\be$ and $(X,Y)\ne 0$ for some $\al,\be\in \cR\union\{0\}$, then for any $C\in \lt$:
$\al(C)(X,Y) = ([C,X],Y) = -(X,[C,Y]) = -\be(C)(X,Y)$
so that $\al+\be =0$.
\epf
}

\noindent {\bf Proof of (QH1).} 
To simplify the notation write $M = \gahr$ and $m=2r$.
\begin{lem}\label{lem: qh1 lem}
Suppose $A,B,S:M\to G$ are maps such that $A=SB$ and $S$ takes values in a fixed unipotent (isotropic) subgroup of $G$.
Let $\al=A^*(\th), \be = B^*(\th)$.
Then $$3 d(\al,\be) = (\be^3) - (\al^3)\in \Omega^3(M,\g).$$
\end{lem}
\pf
Since $S=AB^{-1}$ the element $\si:=S^*(\th)$ is conjugate to $\al-\be$ and so
$$0 = (\si^3) = (\al^3)- (\be^3) +3(\al\be^2) - 3(\al^2\be).$$
Then observe the Maurer--Cartan equations imply 
$d(\al,\be) = -(\al^2\be) + (\al\be^2)$.
\epf
\begin{cor}\label{cor: dsum simplifies}
$3\sum_{i=1}^m d(\ga_{i},\ga_{i-1}) = (\ga^3) - (\ga_m^3)$ where $\ga=\ga_0$.
\end{cor}
\pf 
Lemma \ref{lem: qh1 lem} implies $3d(\ga_{i},\ga_{i-1})= (\ga_{i-1}^3)-(\ga_i^3)$ and so the sum collapses.%
\epf

For the proof of QH1, write $\mu = (\mu_G,\mu_H)$ where $\mu_G = C^{-1}bC, \mu_H = h^{-1}$.
Then $\mu_G^*(\th)$ is conjugate to $\be +\bar\ga -b^{-1}\bar\ga b$ and a direct calculation shows
$$\mu_G^*(\th^3) = (\be^3)  + 3(\bar\ga\be^2) + 3(\bar\ga^2\be) 
+3d(\bar\ga,\bar\be+\Ad_b\bar\ga). $$
Thus from this and Corollary \ref{cor: dsum simplifies}, showing $\mu^*(\th^3) = 6d\omega$ reduces to verifying that
$$
(\be^3)  + 3(\bar\ga\be^2) + 3(\bar\ga^2\be) -(\ch^3) = 
(\ga_m^3) - (\ga^3) + 3(\bar\ga_m\ch^2) + 3(\bar\ga_m^2\ch).$$
But, swapping the sides of $(\ch^3)$ and $(\ga^3)$, this amounts to showing $((\be +\bar\ga)^3) = ((\ch + \bar\ga_m)^3)$, and this is a simple consequence of the fact that (by definition)  $bC=hC_m$.

\noindent {\bf Proof of (QH2).}
First, considering just the $G$ action, since $\mu_G= C^{-1}bC$ we have 
\beq\label{eq: gt+bt}
\mu_G^*(\th+\bar\th) = 
C^{-1}(\be + \bar\be)C
+
C^{-1}(b\bar\ga b^{-1} - b^{-1}\bar\ga b)C.\eeq
If $X\in \g$ and we use primes to denote derivatives along the corresponding fundamental vector field $v_X$ then
$\ga':=\langle v_X,\ga\rangle = X, \bar\ga' = CXC^{-1}, \be'=\eta'=0$.
This enables us to compute 
$$2\omega(v_X,\cdot) = 
(X,\Ad_{(C^{-1}b)} \bar\ga -\Ad^{-1}_{(bC)}\bar\ga 
+ \Ad^{-1}_{C}\bar\be +\Ad^{-1}_{C_m}\ch +\ga_m -\ga).$$
Comparing with \eqref{eq: gt+bt} we see this agrees with  $(X,\mu_G^*(\th+\bar\th))$ provided the relation
$\Ad^{-1}_{C_m}\ch +\ga_m -\ga = C^{-1}\be C$ holds---but this follows easily from the fact that  $hC_m=bC$.
Secondly, for the $H$ action, since $\mu_H= h^{-1}$ we have 
\beq\label{eq: ht+bt}
\mu_H^*(\th+\bar\th) =  -\ch -\bar\ch.
\eeq
Now if  $X\in \h$  then
$\bar\ga_i'=- X, \ga_i' = -C_i^{-1}XC_i, 
\bar\be'=bXb^{-1}-X, \eta'=X-h^{-1}Xh$.
This enables us to compute that 
$2\omega(v_X,\cdot) = (X,\al)$ where
$$\al:=
-\Ad_b\bar\ga-\bar\be+\bar\ga-\ch-\bar\ga_m+\Ad_h\bar\ga_m+
\sum_1^m \Ad_{C_i}\ga_{i-1} -\Ad_{C_{i-1}}\ga_i.$$
Now $\Ad_{C_i}\ga_{i-1} = \bar\ga_i-\bar\si_i$ and  
$\Ad_{C_{i-1}}\ga_i = \bar\ga_{i-1}+\si_i$ 
where $\si_i = S_i^*(\th)$, so (since $\si_i$ pairs to zero with $X\in\lh$) the summation in the expression for $\al$ may be simplified to $\sum_1^m \bar\ga_{i} -\bar\ga_{i-1} = \bar\ga_m-\bar\ga$, which cancel with other terms in $\al$, yielding $2\omega(v_X,\cdot) = -(X, \Ad_b\bar\ga+\bar\be+\ch-\Ad_h\bar\ga_m ).$ In turn this equals $-(X, \ch+\bar\ch)$ (as expected from \eqref{eq: ht+bt}) since $h=bCC_m^{-1}$. This completes the proof of (QH2).

\noindent {\bf Proof of (QH3).}
Fix a point $p\in M$ and a vector $v\in \ker(\omega)\cap \ker(d\mu) \subset T_pM$.  Thus our aim is to show that $v$ is zero.
If we let primes denote components along $v$
then $\ch':= \langle \ch , v\rangle =0$ since $v\in \ker d\mu_H$.
Similarly since  $v\in \ker d\mu_G$ we have
\beq\label{eq: be expn1}
\be' = b^{-1}\bar\ga' b -\bar\ga'.
\eeq
In order to use the condition that $v\in \ker \omega$ we note that $\omega$ may be expanded as follows.
Write $[ij] = S_iS_{i-1}\cdots S_j$ for any $i \ge j$.

\begin{lem}
$$2\omega = 
(\bar\ga,\Ad_b\bar\ga) + (\bar\ga,\bar\ch)
+(\Ad_{[m1]}\bar\ga,\ch)
$$
$$
+\sum_{i=1}^m (\Ad_{[mi]}\si_i,\ch)
+(\bar\ga,\Ad_{h[mi]}\si_i)
- (\bar\si_i,\Ad_{[i1]}\bar\ga)
-\sum_{j <i} (\bar\si_i,\Ad_{[ij]}\si_j)
$$
\end{lem}
\pf
Since $C_i=S_iC_{i-1}$, 
$\ga_i = \ga_{i-1}+C^{-1}_{i-1}\si_i C_{i-1}$ and so 
$(\ga_i,\ga_{i-1}) = (\si_i,\bar \ga_{i-1}),$ and inductively
$\bar\ga_i = \Ad_{[i1]}\bar\ga + \sum_{j=1}^i\Ad_{[ij]}\si_j.$
Thus the term $-\sum_1^m(\ga_i,\ga_{i-1})$ yields the last two terms of the displayed expression,
and $(\bar\ga_m,\ch)$ yields the third and fourth terms.
Finally since $b=hS_m\cdots S_1$ it follows that
$\bar\be = \bar\ch + \sum_1^m \Ad_{h[mi]}\si_i$ and so 
$(\bar\ga,\bar\be)$ yields the second and fifth terms.
(For later use note $\be = \Ad^{-1}_{[m1]}\ch + \sum \Ad^{-1}_{[i1]}\bar\si_i$.)
\epf

Thus if $u\in T_pM$ is arbitrary, it follows directly that
$$2\omega(v,u) = (\Ga,\dot{\bar\ga}) + (\cH, \dot\ch) + 
\sum (\De_i, \dot{\si_i})$$
where the dots denote $u$ components and
$$\Ga = b^{-1}\bar\ga' b - b\bar\ga' b^{-1} 
- \sum \Ad_{h[mi]}\si_i' -\sum \Ad^{-1}_{[i1]}\bar\si_i',$$
$$\cH = \Ad^{-1}_h\bar\ga' + \Ad_{[m1]}\bar\ga' + \sum \Ad_{[mi]}\si_i',$$
$$\De_i = \Ad^{-1}_{h[mi]}\bar\ga'+\Ad_{[i-1,1]}\bar\ga' 
-\sum_{j>i}\Ad^{-1}_{[ji]}\bar\si'_j + \sum_{j<i} \Ad_{[i-1,j]}\si_j'.$$
Now from \eqref{eq: be expn1} and its conjugate by $b$ (and the expansion of $\be$) it follows that  $\Ga=0$.
Similarly it follows that $\cH = 2h^{-1}\bar\ga' h$ and so the vanishing of $\omega(v,u)$ for all $u$ implies $\pi_\lh(\bar\ga')$ is zero, where $\pi_\lh:\g\to \lh$ is the projection.
The only other information we have is that $\De_i$ is orthogonal to 
$\lu_i$ for each $i$, i.e. $\De_i\in \lh\oplus\lu_i$ (where $\lu_i\subset \g$ is the Lie algebra containing $\dot\si_i$). This implies $v=0$, as follows.
First %
observe the expression for $\De_i$ implies
\beq\label{eq: oldobs}
\Ad_{S_i}\De_i - \De_{i+1} = -\bar\si'_i - \si_{i+1}'
\eeq
for $i=1,\ldots,m-1$.
Applying $\pi_\lh$ to this we see that 
$\pi_\lh(\De_i)=\pi_\lh(\Ad_{S_i}\De_i) = \pi_\lh(\De_{i+1})$ and we define $\kappa=\pi_\lh(\De_i)\in \lh$ to be this common value.
Now define
$$T_i =  \Ad_{S_i}\De_i  +  \si_{i+1}' =  \De_{i+1} -\bar\si'_i.$$
Then due to the orthogonality conditions on the $\De_i$ it follows that
$T_i = \si_{i+1}' -\bar\si'_i +\kappa$.
Thus $\De_i +\si_i' = \Ad^{-1}_{S_i}(\kappa)$.
Taking $i=1$ and expanding $\De_i$, this says 
$$b^{-1}\bar\ga' b + \bar\ga' - \sum_{j>1}\Ad^{-1}_{[j1]}\bar\si'_j + \si_1' = \Ad^{-1}_{S_1}(\kappa).$$
In turn since $\be' = \sum\Ad^{-1}_{[j1]}\bar\si'_j=b^{-1}\bar\ga' b - \bar\ga'$ this reduces to 
\beq\label{eq: si1 eqn}
2(\bar\ga' + \si_1') = \Ad^{-1}_{S_1}(\kappa).
\eeq
Then taking the $\lh$ component implies $\kappa=0$, and thus 
$\De_i = -\si_i'$.
Then from \eqref{eq: oldobs} $-\bar\si_i' +  \si_{i+1}' = 
-\bar\si_i' -  \si_{i+1}'$ and so $\si_{i+1}'=0$ (for $i=1,\ldots,m-1$).
In particular $\De_m=0$, and the expression for $\De_m$ now simplifies to 
$$\De_m = \Ad^{-1}_{hS_m}\bar\ga'+\Ad_{[m-1,1]}\bar\ga' +
\Ad_{[m-1,1]}\si_1'=0.$$
Thus $b^{-1}\bar\ga' b + \bar\ga' + \si_1' =0$.
Using \eqref{eq: si1 eqn} this implies $\bar\ga'=0$ and in turn $\si_1'=0$. Thus all the components of $v$ are zero and (QH3) is established.\epfms

\subsection{Pictures}

We will draw pictures of two ways one might want to think about the higher fission spaces. 
The first is helpful to keep track of the monodromy relations, and the second expresses the possible ways to glue surfaces together.
For general linear groups an alternative viewpoint is possible in terms of graphs (see \cite{cqv}).

First recall the standard way to obtain the double $G\times G$ (from \cite{AMM}): one takes an annulus $\Si$ (in the complex plane) with
a marked point on each boundary circle. 
Then the moduli space of flat $C^\infty
$ connections on $G$ bundles over $\Si$ together with a 
framing at each marked point, is naturally a quasi-Hamiltonian $G\times G$ space. 
Choosing two paths generating the fundamental groupoid of $\Si$ (based at the two marked points) enables one to identify the moduli space with $G\times G$ (taking the monodromy of the connections along the chosen generating paths). 
\begin{figure}[ht]
	\centering
	\input{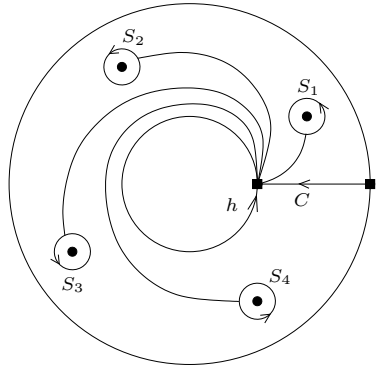}
	\caption{One way to picture the fission spaces}\label{fig: Pi gens}
\end{figure}

One may {\em picture} the fission spaces $\gahr$ in a similar way by puncturing the annulus at $2r$ equally spaced points in its interior, as in Figure \ref{fig: Pi gens}.
(The usual picture for the double is obtained if $r=0$.)
Note that this is not quite ``what is actually happening'' in the derivation of these spaces in terms of the Stokes phenomenon,
but nonetheless is sometimes useful.  
In some sense, in the Stokes phenomenon, the inner circle should be shrunk to zero and the punctures pushed into origin, forming part of the ``analytic halo'' of \cite{MR91}.
One of the surprises of the present article is that one does in fact, nonetheless, get genuine quasi-Hamiltonian spaces in this way.
(So the quasi-Hamiltonian framework does indeed go beyond the context of flat connections on surfaces with boundary.)
In later sections we will in fact take this viewpoint quite seriously and use it (in \S\ref{sn: irreg curves}) to define the basic topological objects: {\em Stokes representations}  
and {\em Stokes $G$-local systems}.%

The second picture explains the name fission, and 
occurs if the group $H$ may be written as a product
$H=H_1\times H_2$ of two groups (the generalisation to arbitrarily many factors is immediate). For example if $G$ is a general linear group this is always the case, unless $H=G$. On the Lie algebra level, the Dynkin diagram of $H$ arises by deleting some nodes in the Dynkin diagram of $G$ (i.e. breaking that of $G$ into pieces).
Then it is more  accurate to replace the annulus by the product of a circle and a $Y$ shaped piece, since one may glue on both a quasi-Hamiltonian 
$H_1$-space and a quasi-Hamiltonian $H_2$-space.
Indeed if we draw a dashed line to represent the analytic halo (through the punctures added above), then after crossing this halo the pieces of surface may drift apart yielding the following picture. 

\begin{figure}[ht]
	\centering
	\input{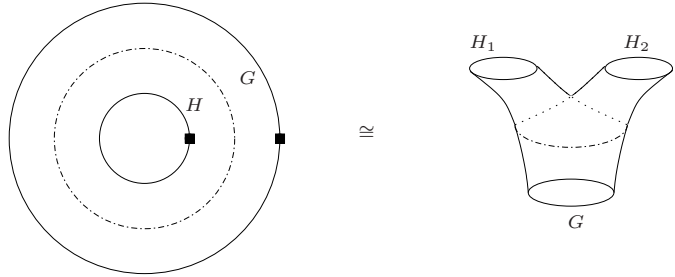}
	\caption{Another way to picture the fission spaces}
\end{figure}

\subsection{Fission varieties}

Thus we now have a large supply of quasi-Hamiltonian spaces and may glue them together and perform the reduction operation to obtain many symplectic manifolds.
\begin{defn}
A ``\,{\em fission variety}\!'' is a symplectic or quasi-Hamiltonian variety obtained 
via the operations of fusion and reduction on spaces of the form

a) conjugacy classes $\cC\subset G$ in arbitrary complex reductive groups $G$,

b) fission spaces $\gahr$, and

c) tame fission spaces $\IM\cong (G\times P)/U$ of \cite{logahoric} Theorem 9.
\end{defn}

\noindent
Many classes of algebraic symplectic manifolds arise as examples:

1) The result of \S\ref{sn: derive vdb} below implies that all of the multiplicative quiver varieties of Crawley--Boevey and Shaw 
\cite{CB-Shaw} (attached to arbitrary graphs) are examples of fission varieties.

2) Since the double $G\times G$ is a special case of a fission space, the complexification of the results of \cite{AMM} show that
moduli spaces of $G$ valued representations of the fundamental group of Riemann surfaces with boundary, are fission varieties.

3) More generally one may consider the parabolic and parahoric extensions of the spaces in 2), by considering tame meromorphic connections on parahoric torsors over smooth compact algebraic curves.
The tame fission spaces enable us to capture the corresponding Betti spaces as fission spaces too (cf. \cite{logahoric});
for $G$ a general linear group such spaces appear in Simpson's Riemann--Hilbert correspondence \cite{Sim-hboncc} revisited from a quasi-Hamiltonian viewpoint by Yamakawa \cite{yamakawa-mpa}.

4) Still more generally one may consider spaces of monodromy and Stokes data classifying arbitrary (unramified) meromorphic connections on $G$-bundles on smooth algebraic curves (and, in turn, on parahoric torsors). In later sections of this paper we will show that these are also examples of fission varieties. This example was our main motivation, and provides the symplectic Betti description of all the hyperk\"ahler manifolds of \cite{wnabh}, together with the analogous spaces for arbitrary reductive groups $G$.

Many other, more exotic, examples are possible however
(an explicit example is described in \cite{fission}).
To see how special the above spaces 1-4) are, note that 
the fission varieties in 1) only involve products of general linear groups (they provide a link to graphs, and then to Kac--Moody root systems).
On the other hand the fission varieties in 2), 3) and 4) privilege one fixed reductive group $G$.
\section{Derivation of $\cB(V,W)$} \label{sn: derive vdb}

Now we will specialise to $r=2$ and show that the Van den Bergh spaces $\cB(V,W)$  arise from the fission spaces.
Recall that if $r=2$ and we perform the reduction by $G$ of $\gahr$ at the value $1$ of the moment map then the resulting space is 
\beq\label{eq: r=2 reduced space v2}
\left\{(S_1,S_2,S_3,S_4,h)\in U_+\times U_-\times  U_+\times U_- \times H \st  hS_4S_3S_2S_1 =1 \right\}.\eeq
This inherits the structure of  quasi-Hamiltonian $H$-space
with moment map $h^{-1}$. (It is symplectic if $H$ is a maximal torus as in \cite{saqh}).
By forgetting $h,S_3,S_4$ this space  \eqref{eq: r=2 reduced space v2} embeds as the (open) subset
of $ U_+\times U_-$ consisting of pairs $(S_1,S_2)$ 
such that $S_2S_1\in P_+P_-\subset G$.
The quasi-Hamiltonian two-form on \eqref{eq: r=2 reduced space v2} (obtained by restricting that on $_G\cA_H^2$ to the subset where $C=b=1$) is
\beq\label{eq: restricted 2form}
\omega = \frac{1}{2}\bigl( (\ga_1,\ga_2) + (\ga_2,\ga_3)\bigr)
\eeq
with $\ga_i= C_i^*(\th)$ and now $C_i=S_i\cdots S_1$ (since $C=1$).
Now we will further specialise to the case of a maximal proper parabolic subgroup  of a general linear group.
Thus suppose we are given two complex vector spaces $V,W$ and take $G=\GL(V\oplus W)$ and choose $P_+$ to be the block upper triangular subgroup (and $H\cong  \GL(V)\times \GL(W)$ to be the block diagonal subgroup).
Then we may write
$$
S_1 = \bmx 1& a \\ 0  &1\emx, \  
S_2 = \bmx 1& 0 \\ b  &1\emx, \ 
S_3 = \bmx 1& c \\ 0  &1\emx, \ 
S_4 = \bmx 1& 0 \\ f  &1\emx, \ 
h = \bmx x & 0 \\ 0  & y\emx$$
where $a, c \in \Hom(W,V), b,f\in \Hom(V,W), x\in \GL(V), y\in\GL(W)$.
One may then check that the relation 
$hS_4S_3 =   S_1^{-1}S_2^{-1} $ from \eqref{eq: r=2 reduced space v2}
is equivalent to the equations
\beq\label{eq: reduced eqns}
x=1+ab, \quad y= (1+ba)^{-1},\quad c=-x^{-1}a,\quad f= -(b+ bab).
\eeq
Thus in this case the space \eqref{eq: r=2 reduced space v2} is isomorphic to the space
$$\{(a,b)\in \Hom(W,V)\oplus\Hom(V,W)\st \det(1+ab)\ne 0\}=\cB(V,W).$$
Further note that the moment map $h^{-1}$ on \eqref{eq: r=2 reduced space v2} has components $(x^{-1}, y^{-1})$ coinciding with the moment map \eqref{eq: vdb mmap} on $\cB(V,W)$. 
Finally we verify the two-forms agree:

\begin{lem}
In the coordinates $(a,b)$
the two-form \eqref{eq: restricted 2form}
equals
$$\frac{1}{2}\bigl(\tr_V (1+ab)^{-1} da\wedge db - \tr_W (1+ba)^{-1}db\wedge da\bigr).$$
\end{lem}
\pf
This is a direct calculation, which may be done completely algebraically: we start with a sum of  words in the noncommuting symbols $da,db,a,b, y$ (noting that $x^{-1} = 1-ayb$) 
and then use 
$bay = 1-y = yba$ to simplify the result (also using the cyclicity of the trace). We find\footnote{Using a symbolic manipulation program.}
$(\ga_1,\ga_2) = \tr (da\wedge db),$ 
$(\ga_2,\ga_3) = -\tr(y db\wedge da) - \tr ( ayb da\wedge db)$
and the result follows since $x^{-1} = 1-ayb$.
\epf

This agrees with the quasi-Hamiltonian two-form \eqref{eq: vdb 2form} on $\cB(V,W)$.
In other words we have established the following:

\begin{thm}\label{thm: isom to vdb}
If $G=\GL(V\oplus W)$ and $H=\GL(V)\times \GL(W)$ then the map 
$$\cA^2(V,W)\spq G\to \cB(V,W);\qquad (h,S_1,S_2,S_3,S_4) \mapsto (a,b)$$  is an isomorphism of  quasi-Hamiltonian $H$-spaces, where $a,b$ are the nontrivial matrix entries of $S_1,S_2$ respectively.
\end{thm}

The spaces $\cB(V,W)$ are the basic building blocks for the multiplicative quiver varieties of Crawley-Boevey--Shaw 
\cite{CB-Shaw}.

\begin{cor}
Any multiplicative quiver variety (in the sense of Crawley-Boevey--Shaw) is a fission variety. 
\end{cor}
\pf
The multiplicative quiver varieties are obtained 
by fusing together many copies of spaces of the form $\cB(V,W)$
and then reducing at certain (central) conjugacy classes
to obtain a symplectic manifold.
So the result is immediate from Theorem \ref{thm: isom to vdb}.
\epf

Theorem \ref{thm: isom to vdb} suggests how to find other building blocks, such as $\cB(U,V,W):=\cA^2(U,V,W)\spq G$, to construct more general multiplicative quiver varieties. This will be taken up in \cite{cqv}.
(It also suggests how one might define analogues of multiplicative quiver varieties for groups besides general linear groups, by considering the reductions $\cB^r:=\gahr\spq G$ in general, e.g. if $H=H_1\times H_2$.)

\begin{rmk}
Note that Van den Bergh \cite{vdb-doublepoisson, vdb-ncqh}
constructs some ``noncommutative quasi-Hamiltonian spaces'' such that after choosing a linear representation one obtains a genuine
quasi-Hamiltonian structure. 
Presumably there are noncommutative quasi-Hamiltonian spaces which yield the higher fission spaces (in the general linear case) upon choosing a representation (and maybe one can prove this by rewriting the algebraic  proofs of QH1-3). However it is unclear if one could obtain the  higher fission spaces for other reductive groups in this way, and this seems important for some applications. 
\end{rmk}

\section{Basic properties}\label{sn: basic props}

In this section  we will establish some basic properties of the spaces 
$\gahr$.

\subsection{Isomonodromy isomorphisms I}

Suppose $P_+\subset G$ is a parabolic with Levi subgroup 
$H\subset P_+$ and opposite parabolic $P_-$.
Write $\cA = \gahr$ for the fission space associated to $(P_+,H)$ and $\cA'$ for the corresponding space associated to $(P_-,H)$ i.e. with the roles of 
$P_+$ and $P_-$ swapped.

An isomorphism $\Theta$ between the spaces $\cA$ and $\cA'$ may be defined as follows.
Suppose 
$(C,h,\bS) \in \cA$, then define  
$(D,g,\bT)=\Theta(C,h,\bS)\in \cA'$ by the formulae
$$T_i = S_{i+1}\qquad\text{for $i=1,\ldots, m-1$}$$
$$T_m = h^{-1} S_1 h,\quad
D= S_1 C,\qquad
g=h$$
where $m=2r$.
Then it is clear that 
$D^{-1}hT_m\cdots T_1D = C^{-1}hS_m\cdots S_1C$
so the moment maps match up.

\begin{figure}[ht]
	\centering
	\input{isomisom9.1.pstex_t}
	\caption{}\label{fig: isomisom}
\end{figure}

\begin{prop} \label{prop: imd im}
The map $\Theta$ is an isomorphism of quasi-Hamiltonian $G\times H$ spaces.
\end{prop}

\pf
Let $D_i = T_i\cdots T_1 D$ so that 
$D_i = C_{i+1}$ if $i<m$ and 
$D_m = h^{-1}S_1hC_m$.
Let $p=hT_m\cdots T_1 = S_1 b S_1^{-1}$.
Thus the aim is to verify that the expression
\eqref{eq: omeg def}
for $2\omega$ on $\cA$ equals 
\beq\label{eq: new 2om}
(\bar\de,\Ad_p\bar\de) 
+ (\bar\de,\bar\cP)
+ (\bar\de_m, \ch) -
\sum_{i=1}^m(\de_{i},\de_{i-1})
\eeq
when the variables are related in this way, where $\de_i = D_i^*(\bar\th)$ and $\bar\cP = p^*(\bar\th)$ etc.
First $\bar\de = \bar\si_1 + S_1\bar\ga S_1^{-1}$, 
where $\bar\si_1 = S_1^*(\bar\th)$, 
so that 
\beq\label{eq: imdim1}
(\bar\de,\Ad_p\bar\de) =  
(\si_1,b\si_1b^{-1}) + 
(\bar\ga,b\si_1b^{-1})+
(\si_1,b\bar\ga b^{-1}) +
(\bar\ga,b\bar \ga b^{-1}).
\eeq
Next since 
$\bar\cP = \bar \si_1 + S_1\bar\be S_1^{-1} - \Ad_{(S_1b)}\si_1$
it follows that
\beq\label{eq: imdim2}
(\bar\de,\bar\cP) =  
(\si_1,\bar\be)
-(\si_1,b\si_1b^{-1})
+(\bar\ga,\si_1)
-(\bar\ga,b\si_1b^{-1})+
(\bar\ga,\bar\be).
\eeq
Upon summing these two expression, four terms cancel, the first two terms of $2\omega$ are obtained and there are three remaining terms.
Now 
$$\de_m = 
\ga_m 
+ C_m^{-1}\eta C_m 
+ \Ad^{-1}_{(hC_m)}\si_1
- \Ad^{-1}_{(S_1hC_m)}\bar\eta
$$
so the
 next terms of \eqref{eq: new 2om} are 
\beq\label{eq: imdim3}
(\bar\de_m,\eta) = 
 (\Ad_{(S_1h)}\bar\ga_m,\bar\eta)
\eeq
since the other terms are zero, and
\beq\label{eq: imdim4}
-\sum_1^{m-1} (\de_i,\de_{i-1}) = -\sum_2^m (\ga_i,\ga_{i-1})
\eeq
and finally
\beq\label{eq: imdim5}
- (\de_m,\de_{m-1}) = -(\de_m,\ga_m) = 
-(\eta,\bar\ga_m)
-(h^{-1}\si_1 h,\bar\ga_m) +
(\Ad^{-1}_{(S_1h)}\bar\eta,\bar\ga_m).
\eeq
Now \eqref{eq: imdim3} cancels the last term of   
\eqref{eq: imdim5}, and clearly all the remaining terms of $2\omega$ are obtained except $-(\ga_1,\ga_0)$. But this equals 
$-(\si_1,\bar\ga)$ which is in \eqref{eq: imdim2}.
Thus to finish we need to check the left over terms disappear. These are $(\si_1,b\bar\ga b^{-1}+\bar\be - h\bar\ga_m h^{-1})$
and they vanish since $hC_m = bC$ (and $\bar\eta$ pairs to zero with $\si_1$).
\epf

In fact the same proof %
also yields a more general statement (one may first break up the unipotent groups into a direct spanning decomposition)---see Section \ref{ssn: rimd}.
Note that $\Theta^{-m}$ (the $m$th power of the inverse of $\Theta$) coincides with the ``twist automorphism'' of \cite{AMM} Theorem 4.5, related  to the Dehn twist of two-manifolds.

\begin{eg} \label{eg: cbs twist}
Suppose $r=2$ and we take $G=\GL(V\oplus W),H=\GL(V)\times\GL(W)$ 
with $P_+$ the block upper triangular subgroup of $G$, as in Section \ref{sn: derive vdb}.
Then the isomonodromy isomorphism of Proposition \ref{prop: imd im}
is
\begin{align}
\cA\quad & \to \quad\cA'; \label{eq: example imd im}\\
(C,h,S_1,S_2,S_3,S_4)&\mapsto (S_1C, h ,S_2,S_3,S_4, h^{-1} S_1 h).\notag
\end{align}
Reducing by $G$ thus yields a commuting square of isomorphisms of quasi-Hamiltonian $H$-spaces:
$$
\begin{array}{ccc}
  \cA\spq G  & \mapright{\cong} & \cA'\spq G \\
\mapdown{\cong} &&  \mapdown{\cong} \\
  \cB(V,W)  & \mapright{\cong } & \cB(W,V)\\
&&\\
(a,b) & \mapsto & (b, -(1+ab)^{-1}a) \\
\end{array}
$$
where the top isomorphism is induced from 
\eqref{eq: example imd im}, the left-hand isomorphism is that of 
Theorem \ref{thm: isom to vdb} (taking the nontrivial matrix entries $(a,b)$ of $S_1,S_2$), and the right-hand isomorphism is the analogue for $\cA'$, i.e. taking the nontrivial matrix entries $(b,c)$ of $S_2,S_3$.
Thus the isomorphism along the bottom is $(a,b)\mapsto (b,c)$ where
$c=-(1+ab)^{-1}a$ as in \eqref{eq: reduced eqns}.
This is the map used in \cite{CB-Shaw} to reverse the orientation of the edges of multiplicative quiver varieties.
\end{eg}

\subsection{Conjugacy isomorphisms}

Suppose that $P\subset G$ is a parabolic subgroup 
with Levi subgroup $H\subset P$.
Suppose $\phi$ is an automorphism of $G$ preserving $H$ and the inner product on $\g$.
Then $Q=\phi(P)$ is again a parabolic of $G$ with Levi subgroup $H$.
A basic example is if $P$ is a Borel and $H$ is a maximal torus and
$\phi$ is the inner automorphism obtained by conjugating by an element of the normaliser of $H$ (representing an element of the Weyl group $N(H)/H$).

Write $\cA = \gahr$ for the space associated to $(P,H)$ and $\cA'$ for the corresponding space attached to $(Q,H)$.
Then the map $\cA\to \cA'$ defined by $\phi$ is clearly an isomorphism of spaces and relates the quasi-Hamiltonian two-forms (but it is not an isomorphism of quasi-Hamiltonian spaces since it does not relate the moment maps if $\phi$ acts nontrivially on $G\times H$). Nonetheless upon reduction the fact the two-forms are related implies the resulting symplectic manifolds will be isomorphic.

For example if  $H$ is a maximal torus then $\cA/G$ and $\cA'/G$ will be isomorphic Poisson manifolds, with an isomorphism given by $\phi$.
(Up to a covering, if $r=1$ these spaces are two realisations of the dual Poisson Lie group $G^*$, and examples of such isomorphisms appear in \cite{bafi} Lemma 3.5.)

\subsection{Inversion anti-isomorphisms}

We will say a map $\phi:M\to N$ between two quasi-Hamiltonian $G$-spaces is an {\em anti-isomorphism} if it is an isomorphism of spaces
and negates the two-form and inverts the moment map, i.e.
$$\phi^*(\omega_N) = -\omega_M,\qquad
 \mu_N\circ \phi = \iota\circ \mu_M:M\to G
$$
where $\iota:G\to G$ is the inverse map.
A simple example is the `flip' anti-isomorphism 
\beq\label{eq: flip anti-isom}
\cB(V,W)\to\cB(W,V); \ (a,b)\mapsto(b,a).
\eeq

Now suppose $P_+\subset G$ is a parabolic subgroup with Levi subgroup $H\subset P$ and opposite parabolic $P_-$.
Write $\cA = \gahr$ for the space associated to $(P_+,H)$ and $\cA'$ for the corresponding space associated to $(P_-,H)$ i.e. with the roles of $P_+$ and $P_-$ swapped.

\begin{prop} \label{prop: inversion anti-isom}
The map $\cA\to \cA'; \ (C,h,\bS)\mapsto (C,h^{-1},\bT)$
where 
$$T_i = hS^{-1}_{m+1-i}h^{-1}$$
is an anti-isomorphism of quasi-Hamiltonian $G\times H$-spaces.
\end{prop}

Geometrically this corresponds to reversing the orientation.

\pf
Suppose we change coordinates on $\cA$ as follows:
$$d_i = h^{-1}S^{-1}_{m+1-i}h,\qquad \ e_i = S_i$$
for $i=1,\ldots,r$ (with $m=2r$ as usual)
and set
$$D_i = d_{i}\cdots d_1 C,\qquad
E_i = e_{i}\cdots e_1 C.$$
In particular now $C=E_0=D_0$ 
and we define $E=E_r, D=D_r$
so that the $G$ component of the moment map is
$\mu_G = D^{-1}h E.$
Then one may check the expression \eqref{eq: omeg def}
for $2\omega$ on $\cA$
equals
$$
(\bar\cD,\Ad_b\bar\cE) 
+ (\bar\cD,\bar\eta)
+ (\bar\cE,\eta) 
+
\sum_{i=1}^r(\cD_{i},\cD_{i-1})-(\cE_{i},\cE_{i-1})
$$
where $\cD_i = D^*(\th), \cE_i=E^*(\th)$ etc.
Now the map to $\cA'$ corresponds to swapping   $D_i$ 
and $E_i$ for each $i$ (i.e. to swapping $d_i$ and $e_i$) and inverting $h$.
Then it is clear that  the expression for $\omega$ is negated.
\epf

Note that the maps
\beq\label{eq: inner twist}
\cA\to \cA; \ (C,h,\bS)\mapsto (hC,h,h\bS h^{-1})
\eeq
and 
\beq\label{eq: outer twist}
\cA\to \cA;  \ (C,h,\bS)\mapsto (b^{-1}C,h,\bS)
\eeq
are both quasi-Hamiltonian automorphisms of $\cA=\gahr$
where $b=h S_m\cdots S_1$.
(If $p= (C,h,\bS)\in \cA$ then the first map is 
$p\mapsto \mu_H(p)\cdot p$ and the second is $p\mapsto \mu_G(p)\cdot p$, so this follows from \cite{AMM} Remark 4.2.)
Thus one can conjugate the above inversion anti-isomorphism in various ways to get equivalent versions.

\ppb{

\begin{prop}
The map $\cA\to \cA'; \ (C,h,\bS)\mapsto (D,g,\bT)$
where 
$$D=hS_m\cdots S_1C,\  g=h^{-1}, \ T_i = hS^{-1}_{m+1-i}h^{-1}$$
is an anti-isomorphism of quasi-Hamiltonian $G\times H$-spaces.
\end{prop}

Geometrically this corresponds to reversing the orientation.

\pf

??? leave as exercise

???pass through almost symmetric form of $\omega$ in first derivation, then its quite easy

Suppose we change coordinates in $\cA$ as follows:

$$v_i = S_{r+1-i}\qquad \ u_i = h^{-1}S^{-1}_{r+i} h$$
for $i=1,\ldots,r$
and set
$$D_i = u_{i+1}\cdots u_r D$$
$$E_i = v_{i+1}\cdots u_r C$$
where $D = bC$ with $b=hS_m\cdots S_1$ as usual.
In particular now $C=E_r, D=D_r$ and we define $E=C=E_r$.
Then one may check the expression \eqref{eq: omeg def}
for $2\omega$ on $\cA$
equals
$$
(\bar\cE,\Ad_b\bar\cE) 
+ (\bar\cE,\bar\be)
+ (\eta, \bar\cE_0) 
+
\sum_{i=1}^r(\cD_{i-1},\cD_{i})+(\cE_{i},\cE_{i-1})
$$
where $\cD_i = D^*(\th), \cE_i=E^*(\th)$ etc.
Now the map to $\cA'$ corresponds to swapping the roles of $D_i$ 
and $E_i$ for each $i$ (and also $u_i$ and $v_i$) and inverting $h$.
Then it is clear that the sum in the expression for $\omega$ is negated and it remains to see what happens to the three initial terms, which is easy using the relations $D=bE, D_0=hE_0$. 
\epf

???do we want to do this check of the expressions

??? not D,E not as in saqh, but might be nice to reorder to be similar

}

\section{Further properties}

This section will establish some further properties, which will be useful later on.
First it is convenient to consider more general products of unipotent groups.
Fix a connected complex reductive group $G$ and a subgroup $H\subset G$, such that $H$ is a Levi subgroup of some parabolic subgroup of $G$.
Then define an ``ordered set of unipotent subgroups of $G$ subordinate to $H$'' to be  
a collection $\cU = (U_1,\ldots U_m)$ of unipotent subgroups $U_i\subset G$ such that 

1) for each $i$ there is a parabolic subgroup $P_i\subset G$ such that $H\subset P_i$ is a Levi subgroup, and $U_i\subset \Rad_u(P_i)$ is a subgroup of the unipotent radical of $P_i$, and

2) each $U_i$ is normalised by $H$, i.e. conjugation by any $h\in H$ preserves each  $U_i$. 

Eventually more specific collections of subgroups will be considered, but this definition is convenient to establish some inductive statements needed later.
Given such an ordered set of unipotent subgroups $\cU$, consider the space
$$\cA = G\times H\times \cU$$
with two-form $\omega$ specified (as usual) by 
\beq  \label{eq: omeg def 2}
2 \omega = 
(\bar\ga,\Ad_b\bar\ga) 
+ (\bar\ga,\bar\be)
+ (\bar\ga_m, \ch) -
\sum_{i=1}^m(\ga_{i},\ga_{i-1})
\eeq
i.e. as in \eqref{eq: omeg def}, but with $S_i\in U_i$ now.
(The other notations are the same, so for example
$b=hS_m\cdots S_1, C_i = S_i\cdots S_1C $.)
Note that there is still an action of $G\times H$  on $\cA$  (via the same formulae as before) which preserves $\omega$. Moreover the map
$\mu = (C^{-1}bC, h^{-1})$ from 
$\cA$ to $G\times H$ is well-defined and  equivariant.
Of course at this level of generality $\omega$ will not necessarily be a quasi-Hamiltonian two-form on $\cA$.

\subsection{Direct spanning equivalence}\label{ssn: direct spanning equiv}

Suppose now that  $V_1,\ldots,V_k$ are consecutive subgroups in $\cU$ for some integer $k$ (i.e. for some $i$, that  
$V_1=U_i, V_2=U_{i+1}, \ldots, V_k=U_{i+k-1}$).
Then define 
$$W=V_k\cdots V_1 = \{ S_k\cdots S_1 \in G \st S_i\in V_i\}\subset G.$$
Suppose further that $W$ is ``directly spanned'' by the $V_i$, i.e. 
$W$ is a unipotent subgroup of $G$ and the product map 
$V_k\times\cdots\times V_1\to W$ is an isomorphism of varieties (but not necessarily of groups) cf. Borel \cite{Bor91} \S14.3.
Then we can consider a new  ordered set of unipotent subgroups $\cU'$
by deleting $U_i,\ldots,U_{i+k-1}$ from $\cU$ and inserting $W$ in their place.
Correspondingly there is a space $\cA' = G\times H\times \cU'$ which again has a two-form, denoted $\omega'$.
The product map $V_k\times\cdots \times V_1\to W$ extends to give an isomorphism 
$\cA\to \cA'$ (which is the identity on $G,H$ and the other unipotent groups).

\begin{lem} \label{lem: ds lem}
The isomorphism $\cA\to \cA'$ given by the product map relates the two-forms $\omega$ and $\omega'$.
\end{lem}
\pf
Since everything else is unchanged, this follows from the general fact that 
 $$(\ga_k,\ga_0) = (\ga_k,\ga_{k-1})+\cdots + (\ga_1,\ga_0)$$
if $\ga_i=C_i^*(\th)$ where $\th$ is the left-invariant Maurer-Cartan form on $G$ and  $C_i:V_1\times\cdots\times V_k\times G\to G$ 
is the map taking $(S_1,\ldots,S_k,C)$ to $C_0 = C, C_i=S_i\cdots S_1C$, respectively.
By induction it is sufficient to show 
$(\ga_i,\ga_0) = (\ga_i,\ga_{i-1})+  (\ga_{i-1},\ga_0)$.
Now if $S_f=S_{i-1}\cdots S_1$ then
since $C_{i-1}=S_f C$ it follows that 
$\ga_{i-1} = C^{-1}\si_f C + \ga_0$ where $\si_f=S_f^*(\th)$, 
so that 
$$(\ga_i,\ga_{i-1}) = (\ga_i,C^{-1}\si_f C) + (\ga_i,\ga_0),
\qquad\text{and}\qquad
(\ga_{i-1},\ga_0) = (\si_f,\bar\ga_0).$$
Thus we should show that $(\si_f,\bar\ga_0-C\ga_iC^{-1})=0$.
However if we define $S_g = S_i\cdots S_2S_1$ then $S_g= C_iC^{-1}$ so 
$\si_g=C\ga_iC^{-1}-\bar\ga_0$, and the result follows from the fact that $(\si_f,\si_g)=0$ (since they take values in the isotropic subspace $\Lie(W)$ of $\g$).
\epf

More generally define two sets $\cU$ and $\cU'$  of ordered unipotent
subgroups to be ``direct spanning equivalent'' if they are related by a sequence of such isomorphisms (or their inverses).
Lemma \ref{lem: ds lem} implies any isomorphism $\cA\to \cA'$ coming from a direct spanning equivalence, relates their two-forms.
Thus the following is now immediate.
\begin{cor} \label{cor: dseq}
If $(\cA,\omega)$ and $(\cA',\omega')$ are direct spanning equivalent, and
$(\cA,\omega)$ is a quasi-Hamiltonian $G\times H$-space, then so is $(\cA',\omega')$.
\end{cor}
(It is clear that the actions and moment maps match up.)
For example in the case where $H$ is a maximal torus of $G$ it follows from \cite{Bor91} \S14.4 that each $U_i$ in any ordered set $\cU$ of unipotent groups is directly spanned, in any order, by all the (one dimensional) root groups it contains. Clearly one could subsequently re-assemble the resulting  (ordered set of) root groups into another ordered set of larger unipotent groups in many different ways.

The notion of direct spanning equivalence may be viewed as an abstraction of the relation between ``Stokes matrices'' and ``Stokes factors'' appearing in the $\GL_n$ case in Balser--Jurkat--Lutz \cite{BJL79} \S4.
The extension to other groups is in \cite{bafi} and the observation that this situation is well understood in terms of direct spanning subgroups
is in \cite{smid} (p.156) and \cite{bafi}.
It will become clear in \S\ref{ssn: nesting} (and \S\ref{ssn: qh on stokes})
that there are many quasi-Hamiltonian spaces of the form 
$\cA=G\times H\times \cU$ that are not direct spanning equivalent to one of the fission spaces $\papk{G}{H}{r}$.

\subsection{Isomonodromy isomorphisms II}\label{ssn: rimd}

Suppose $\cA = G\times H\times \cU$ as above with $\cU=(U_1,\ldots,U_m)$ an arbitrary ordered set of unipotent subgroups, and we define
$\cA' = G\times H\times \cU'$ where 
$$\cU' = (U_2\ldots,U_m,U_1)$$
i.e. $U'_i = U_{i+1}$ for all $i\ne m$, and $U_m'=U_1$.
Then we may define an isomonodromy isomorphism $\Theta:\cA\to \cA'$ as in Section \ref{sn: basic props}, i.e.
suppose 
$(C,h,\bS) \in \cA$ and define  
$(D,g,\bT)=\Theta(C,h,\bS)\in \cA'$ by the formulae
$$T_i = S_{i+1}\qquad\text{for $i=1,\ldots, m-1$}$$
$$T_m = h^{-1} S_1 h,\quad
D= S_1 C,\qquad
g=h.$$
Then it is clear that 
$D^{-1}hT_m\cdots T_1D = C^{-1}hS_m\cdots S_1C$
so the moment maps match up.
Examples of such (refined) isomorphisms were considered in \cite{bafi} Proposition 3.7.
The proof of Proposition \ref{prop: imd im} now goes through verbatim to establish:

\begin{prop} \label{prop: imd im2} 
The map $\Theta$ relates the two-forms on $\cA$ and $\cA'$.
\end{prop}
Thus if both of the spaces $\cA$ and $\cA'$ are quasi-Hamiltonian then $\Theta$ is a quasi-Hamiltonian isomorphism.

For example this can be used to show that up to isomorphism the fission spaces  do not depend on the choice of parabolic subgroup $P_+$,
given a fixed Levi subgroup $H\subset G$ 
(this follows from Theorem \ref{thm: local qh imd conn} below).

\subsection{Nesting} \label{ssn: nesting}
Now consider the situation where there is an intermediate
 reductive group $K$ between $G$ and $H$
$$H\subset K\subset G,$$
such that we can consider ordered sets of unipotent groups 
$$\cU,\  \cU', \ \cU'',\qquad \text{subordinate to}\qquad 
H\subset G,\  H\subset K, \  K\subset G$$ 
respectively 
(so that $K$ is again the Levi subgroup of a parabolic of $G$).
As above we thus get spaces 
$$
\cA = G\times H\times \cU,\quad 
\cA' = K\times H\times\cU',\quad
\cA'' = G\times K \times \cU''$$ 
with two forms $\omega,\omega',\omega''$ respectively.

Without loss of generality we may assume (by introducing some trivial unipotent groups $\{1\}$ via direct spanning equivalence) that  each of 
$\cU, \cU'$ and $\cU''$ contains the same number $m$ of groups.
Suppose now that  $U_i$ is directly spanned by $U_i',U_i''$ for each $i$, where $U_i$ is the $i$th group in $\cU$ etc.
(In the present situation, since $K$ normalises $U_i''$, this direct spanning condition just means $U_i$ is the semidirect product of $U_i'$ and $U_i''$.)
Then we can identify $\cA$ with the ``gluing'' of $\cA'$ and $\cA''$ along $K$, as follows. 
We have not assumed the spaces are quasi-Hamiltonian, but the definition of gluing is as one would guess:
Let $(D,h,\{A_i\})$ denote a point of $\cA'$ (with $A_i\in U'_i$)
and let $(C,k,\{B_i\})$ denote a point of $\cA''$ (with $B_i\in U_i''$)
then the gluing is defined to be
$$\cA' \glue{K} \cA'' =
 \{(D,h,\{A_i\},C,k,\{B_i\})\in \cA'\times \cA''\st k=h A_m\cdots A_1 \}/K. $$
Since the action of $K$ is free, we can remove it by setting $D=1$, and thus identify this gluing with
$$
 \{(C,h,\{A_i\},\{B_i\})\}= G\times H\times \cU'\times \cU''.$$
Since each  $U_i$ is directly spanned by $U_i',U_i''$ 
the product $\cU'\times \cU''$ is isomorphic to $\cU$, and so the gluing is isomorphic to $\cA$, but we will choose the isomorphism so that the two-forms match up, as follows.
Define $D_0 =D= 1$ and 
$$D_i = A_i\cdots A_1\in K,$$
then define $S_i\in U_i$ to be 
$$S_i = D_iB_iD_{i-1}^{-1}.$$
Said differently $S_i = A_i\wt B_i$
where  
$\wt B_i := D_{i-1} B_i D^{-1}_{i-1}$.
Clearly  $\wt B_i\in U_i''$ since  $K$  normalises $U''_i$, and so in turn $S_i\in U_i$ by our direct spanning assumption.
Thus taking the point $(C,h,\{S_i\})\in \cA$ defines an isomorphism between the gluing and $\cA$.

\begin{thm} \label{thm: nesting}
This isomorphism $$\cA' \glue{K} \cA'' \quad \cong \quad \cA$$
relates the two-form $\omega$ on $\cA$ to that induced on the left-hand side from  $\omega'$ and $\omega''$.
In particular if $\cA'$ and $\cA''$ are quasi-Hamiltonian, then so is $\cA$.

\end{thm}
\pf
The ``induced two-form'' on the left-hand side is just the restriction of the two-form  on the product $\cA'\times \cA''$ to the subvariety where $D=1$ and $k=hA_m\cdots A_1$ (this subvariety is then identified with $\cA$).  
Thus we must show there is the following equality of two forms on this subvariety:
$$
({\bar\de},\Ad_{k}\bar \de) 
+ ({\bar \de},\bar\kappa)
+ (\bar\de_m, \ch) -
\sum_{i=1}^m(\de_{i},\de_{i-1})$$
$$+({\bar\cE},\Ad_p\bar\cE) 
+ ({\bar \cE},\bar\cP)
+ (\bar\cE_m,\kappa) -
\sum_{i=1}^m(\cE_{i},\cE_{i-1})  $$
$$=
({\bar\ga},\Ad_b\bar\ga) 
+ ({\bar \ga},\bar\be)
+ (\bar\ga_m, \ch) -
\sum_{i=1}^m(\ga_{i},\ga_{i-1}) $$
where $p=kB_m\cdots B_1$, $k = hA_m\cdots A_1$,
$\kappa= k^*(\theta),$
$\ga_i=C_i^*(\th), \de_i=D_i^*(\th), \cE_i = E_i^*(\th), 
E_i =  B_i\cdots B_1C.$
Firstly $D=1$ so $\bar \de=0$, simplifying the first line.
Secondly $p=kB_m\cdots B_1 
= hA_m\cdots A_1B_m\cdots B_1=\cdots = hS_m\cdots S_1=b$
so some terms cancel and consequently we must show 
$$\sum_{i=1}^m(\de_{i},\de_{i-1})+(\cE_{i},\cE_{i-1})  =
(\bar\cE_m, \kappa)+ (\bar\de_m-\bar\ga_m, \ch) +
\sum_{i=1}^m(\ga_{i},\ga_{i-1}). $$
Now since $C_i = D_iE_i$ we find 
$\ga_i = E_i^{-1}\de_i E_i + \cE_i$
and so $$(\ga_i,\ga_{i-1}) = 
(\cE_i,\cE_{i-1}) + (\de_i,\Ad_{B_i} \de_{i-1})
+(\Ad^{-1}_{B_i}(\bar\cE_i),\de_{i-1}) + (\de_i,\Ad_{B_i}\bar\cE_{i-1})
$$
using the fact that $B_i = E_iE_{i-1}^{-1}$.
Now observe that 
$(\de_i,\Ad_{B_i} \de_{i-1})=(\de_i,\de_{i-1})$; 
indeed if we %
consider 
$K\sdp U_i'' \hookrightarrow G$ then $\Ad_{B_i} \de_{i-1}$ takes values in  $\Lie(K)\oplus \Lie(U''_i)$ and has $\Lie(K)$ component  $\de_{i-1}$---but this decomposition is orthogonal so the claim follows. 
On the other hand using $C_m=D_mE_m$ to expand $\bar\ga_m$
and $k=hD_m$ to expand $\kappa$ we see 
$(\bar\cE_m, \kappa)+ (\bar\de_m-\bar\ga_m, \ch)
= (\bar\cE_m, \de_m)$
and so we are reduced to showing:
$$\sum_{i=1}^m (\Ad^{-1}_{B_i}(\bar\cE_i),\de_{i-1}) + (\de_i,\Ad_{B_i}\bar\cE_{i-1}) =-(\bar\cE_m, \de_m).$$ 
But this statement holds for all $m$, and can be proved by induction on $m$ as follows.  For $m=0$ it is trivial as both sides are zero. Otherwise assume it holds for $m=n-1$, so the statement for $m=n$ reduces to showing 
\beq\label{eq: last step induction}
-(\bar\cE_n, \de_n) = (\Ad^{-1}_{B_n}(\bar\cE_n),\de_{n-1}) + (\de_n,\Ad_{B_n}\bar\cE_{n-1})-(\bar\cE_{n-1}, \de_{n-1}).
\eeq
Now $E_n=B_nE_{n-1}$ so 
$\bar \cE_n = \bar\cB_n+\Ad_{B_n}\bar\cE_{n-1}$, where 
$\bar\cB_n = B_n^*(\bar\th)$.
Now $(\de_n,\bar\cB_n)=0$ (since 
$\Lie(K)\perp\Lie(U_n'')$)  so the left-hand side of 
\eqref{eq: last step induction} equals 
$-(\Ad_{B_n}\bar\cE_{n-1},\de_n)$, which cancels with the second term on the right-hand side.
Substituting $\bar \cE_n$, the remaining terms on the right-hand side simplify to
$
(\cB_n,\de_{n-1})$, which is zero since 
$\Lie(K)\perp\Lie(U_n'')$.
\epf

In general  the gluing of two such quasi-Hamiltonian spaces end to end will be referred to as ``nesting''.
This will be used as the key inductive step in Section \ref{sn: local Stokes}
to establish the quasi-Hamiltonian structure on spaces of unramified Stokes data.  
It may also be used to establish various isomorphisms between  fission spaces, as follows.

\ppb{
???draw annulus picture:
 gluing 2 different levels: 
Consider two concentic annuli and Stokes data as above on both annuli (without loss of generality---by allowing redundant Stokes factors---we may assume the two annuli have the same anti-Stokes directions).
Write $A_i=S_i^{(1)}$ for the Stokes factors on the inner annuli and let $S_i^{(2)}$ be the Stokes factor on the outer annuli.
}

Suppose $H\subset K\subset G$ as above   and we choose parabolic subgroups
$P_+ \subset G, P_+'\subset K, P_+''\subset G$ with Levi subgroups $H,H,K$ respectively. Then for any integer $r\ge 1$ the fission spaces 
$$\cA = \papk{H}{G}{r},\quad  \cA'=\papk{H}{K}{r},\quad \cA''=\papk{K}{G}{r} $$
are well defined (using these choices of parabolic subgroups).
Now suppose further that the parabolics have been chosen so that $P_+$ is the semidirect product $P'_+\sdp U_+'' $ of $P_+'$ and the unipotent radical $U_+''$ of $P_+''$ (i.e. $P_+\cap K = P_+'$). 
In this situation we thus have

\begin{cor}\label{cor: nesting fission spaces}
The nesting of two fission spaces of the same level $r$ is again a fission space: 
$$\papk{H}{K}{r}\glue{}\!\!\papk{K}{G}{r}\ \cong \papk{H}{G}{r}$$
as quasi-Hamiltonian $G\times H$-spaces.
\end{cor}
\pf
This follows from Theorem \ref{thm: nesting} since $U_\pm = U'_\pm\sdp U''_\pm$.
\epf

\section{Stokes data for complex reductive groups}\label{sn: local Stokes}

Fix a connected complex reductive group $G$ with Lie algebra $\g$.
In brief, the aim of this section is to define the space of Stokes data 
$\ISto(Q)$ attached to an {\em irregular type} $Q$;
such $Q$ determines a finite set of singular directions $\IA\subset S^1$ and a (complicated) unipotent group $\ISto_d\subset G$ 
for each singular direction $d\in \IA$, and $\ISto(Q)$ is the product of these unipotent groups. Then we will define a slightly larger space $\cA(Q)$ and show it is quasi-Hamiltonian.

Choose a maximal torus $T\subset G$ and let $\lt\subset \g$ denote the corresponding Lie algebras.
Let $\Delta$ be a complex disc and let $a\in \Delta$ be a marked point.
Let $\wh \cO$ denote the formal completion at $a$ of the ring of holomorphic functions on $\Delta$ and let $\wh \cK$ denote its field of fractions.

\begin{defn} \label{def: irtype}
An {\em (unramified)  irregular type} at $a$ is an element
$$Q \in \lt(\wh \cK) /\lt(\wh \cO).$$ 
\end{defn}

One may think of an irregular type as a $\lt$-valued meromorphic function germ, well defined modulo holomorphic terms. 
Explicitly, if we choose a  local coordinate $z$ on $\Delta$ vanishing at $a$,
then 
$\wh\cO= \IC\flb z \frb, \wh \cK = \IC\flp z \frp$, and so 
then an irregular type $Q$ may be  written in the form
$$Q = \frac{A_r}{z^{k_r}} + \cdots + \frac{A_{1}}{z^{k_1}}$$
for integers $0< k_1 < \cdots < k_r$
and elements $A_i\in \lt\subset \g$ for $i=1,\ldots, r$.
(The more abstract definition is  coordinate independent and so will be useful later.)

Let $H\subset G$ be the stabiliser of $Q$ under the adjoint action i.e. $H=\{g\in G \st \Ad_g(A_i)=A_i$ for all $i\ge 1\}$, so $H$ is again a connected complex reductive group with maximal torus $T$.
We will (abusively) write $H=C_G(Q)$ and call it the 
centraliser of $Q$.
Let $\cR\subset \lt^*$ be the set of roots of $\g$ relative to $\lt$ and recall the root space decomposition 
$$\g = \lt \oplus\bigoplus_{\al\in \cR} \g_\al$$
where $\g_\al = \{ X\in \g \st [Y,X] = \al(Y)X\text{ for all $Y\in \lt$}\}$ is the (one dimensional) root space of $\al\in \cR$.
Thus for each root $\al \in \cR$, we may define
$$q_\al = \al \circ Q$$
which is just a meromorphic function (modulo holomorphic terms).
Define the degree $\deg(q_\al)$ of $q_\al$ to be its its pole order at $a$;
using the coordinate $z$ we may identify $q_\al$ with an element of 
$z^{-1}\IC[z^{-1}]$ and so $\deg(q_\al)$ is the degree of the polynomial 
$q_\al(1/z)$ (it is an integer $\ge 0$, 
equal to zero if $q_\al$ does not have a pole at $a$).

Now we start to describe the space of Stokes data attached to $Q$.
This is abstracted from the $\GL_n(\IC)$ case studied in
\cite{MR91, L-R94}, as was previously done in \cite{bafi} in the case when $A_r$ is regular semisimple (and again the use of the notion of direct spanning subgroups of unipotent groups simplifies things).  
Let $\wt \Delta\to\Delta$ denote the real oriented blow-up at $a$, replacing $a$ with the circle $S^1$ of real oriented tangent directions at $a$.
Given a root $\al\in \cR$ consider the  function 
$\exp(q_\al(z))$
as $z$ approaches zero along rays in various directions $d\in S^1$.

\begin{defn} \label{def: sing dirn}
A direction $d\in S^1$ will be said to be a {\em singular direction supported by $\al$,} (or an {\em anti-Stokes direction})  if
$\exp(q_\al(z))$ 
has maximal decay as $z\to 0$ in the direction $d$.
\end{defn}

(Thus if $c_\al/z^k$ is the most singular term of $q_\al$, these are the directions along which the function $c_\al/z^k$ is real and negative.)
Let $\IA\subset S^1$ be the finite set of singular directions (for all roots $\al$).
If $d\in \IA$ let $$\R(d)\subset \cR$$ denote the (nonempty) subset of roots supporting $d$.
Further, given an integer $k$, let 
$$\cR(d,k)\subset \cR(d)$$
denote the subset of roots $\al\in \cR(d)$ such that $\deg(q_\al) = k$.
For example (similarly to \cite{bafi}) if $r=k_1=1$ and $\Delta$ is the unit disc in $\IC$, and $Q=-A/z$ (so that $dQ=Adz/z^2$) then $\IA$ is the set of rays from $0$ to the nonzero points  in the set
$\langle \cR , A \rangle \subset \IC$,
where the angled brackets %
denote the natural pairing between $\lt^*$ and $\lt$.
Further $\cR(d)$ is then the set of roots landing 
on the ray $d$ (and in $\IC^*$)
under the map $\langle \,\cdot\,, A \rangle:\cR\to \IC$. 

\begin{lem}\label{lem: stokes gps}
Each  of the sets $\cR(d), \cR(d,k)$ is a closed subset of some system  of positive roots in $\cR$ (i.e. they are ``special'' in the sense of Borel \cite{Bor91} \S IV.14.5).
\end{lem}

\pf 
See appendix \ref{sn: stokes gps}.
\epf

Thus if 
$U_\al=\exp(\g_\al)\subset G$ 
is the root group corresponding to $\al \in \cR$, 
it follows from  \cite{Bor91} \S IV.14.5 
that $\{ U_\al \st \al \in \cR(d)\}$ ``directly spans'' in any order a unipotent subgroup of $G$. 
This means that if we choose any total ordering of $\cR(d)$ and consider the product map $\phi : \Prod_{\al\in \cR(d)} U_\al\to G$ (with the product taken in the chosen order) then $\phi$ is an algebraic isomorphism (of spaces, not groups) onto its image and this image is a well defined subgroup of $G$ independent of the chosen order of the factors.
\begin{defn} \label{def: stokes gp}
The {\em Stokes group} $\ISto_d$ 
associated to the singular direction $d\in \IA$ is 
the unipotent subgroup of $G$ corresponding to $\cR(d)\subset \cR$:
$$\ISto_d=  \phi\left(\Prod_{\al\in \cR(d)} U_\al\right) \subset G.$$
It has  dimension $\dim_\IC\ISto_d= \#\cR(d)$ and has 
Lie algebra 
$\bigoplus_{\al\in \cR(d)}\g_\al\subset \g$.
\end{defn}

Similarly, for any integer $k$,  define the level $k$ Stokes group, $\ISto_d(k)$, to be the image of 
$\Prod_{\al\in \cR(d,k)} U_\al $ in $G$ and deduce 
(from the direct spanning property, again using \cite{Bor91} \S IV.14.5) 
that the product map gives an isomorphism of spaces:
$$\ISto_d(k_1)\times\cdots  \times \ISto_d(k_r)\cong \ISto_d.$$
(As a group  $\ISto_d$ is the semidirect product of the groups on the left, with $\ISto_d(k_i)$ acting by conjugation on $\ISto_d(k_j)$ for $i<j$.)
Again since the root groups are one-dimensional 
 $\dim\ISto_d(k)= \#\cR(d,k)$.
We will say that the level $k$ is {\em effective} along direction $d$ if
$\cR(d,k)$ is nonempty. 
(We leave it as an exercise to check that 
 the above definitions do not change if we pass to a different 
Cartan subalgebra $\lt$ containing each coefficient $A_i$ of $Q$.)

\begin{defn}
The {\em space of Stokes data}, $\ISto(Q)$, associated to $Q$ is the product of all the Stokes groups:
\beq \label{eq: stokes data sto Q}
\ISto(Q) = \Prod_{d\in \IA}\ISto_d.
\eeq
\end{defn}

(Here we do not take the image of the product in $G$.)
The appearance of such spaces in the local classification of meromorphic connections will be described in Appendix \ref{apx: stokes from connections}.

\subsection{Quasi-Hamiltonian structure on Stokes data}\label{ssn: qh on stokes}

Now choose a singular direction $d_1\in \IA$ and label the other singular directions $d_2,\cdots,d_s$ so that $d_{i+1}$ is next after $d_i$ when turning in a positive sense, so that 
$\ISto(Q) =  \Prod_1^s \ISto_i$
where $\ISto_{i}= \ISto_{d_i}$ and we will
denote elements of $\ISto_{i}$ by $S_i$.
(Beware some references work with singularities at $\infty$, so the notion of positive sense is reversed.)
The aim of the rest of this section is to establish:

\begin{thm} \label{thm: main qh thm}
The space 
$$\cA(Q) = G\times H\times  \ISto(Q)$$
is a quasi-Hamiltonian $G\times H$-space with moment map 
$\mu:\cA(Q)\to G\times H$ given by
$$\mu(C,h,S_1,\ldots,S_s) = (C^{-1}hS_s\cdots S_2S_1C, h^{-1})$$
and two-form $\omega$, where 
\beq  \label{eq: omeg def2}
2 \omega = 
(\bar\ga,\Ad_b\bar\ga) 
+ (\bar\ga,\bar\be)
+ (\bar\ga_s, \ch) -
\sum_{i=1}^s(\ga_{i},\ga_{i-1})
\eeq
with
$\ga_i = C_i^*(\th),\bar\ga_i = C_i^*(\bar\th),
\eta = h^*(\th_H),\bar\be = b^*(\bar\th)$
where $\th, \bar \th$ are the Maurer--Cartan forms on $G$ (and 
$\th_H, \bar \th_H$ are the Maurer--Cartan forms on $H$), and 
$C_i=S_i\cdots S_1 C, b=hS_s\cdots S_2S_1$.
\end{thm}

The key strategy is that such spaces may be obtained by gluing simpler spaces of the form $\papk{G}{H}{r}$ for a nested sequence of reductive groups (and decreasing integers $r$).
If we write $\mu=(\mu_G,\mu_H)$, the element $h=\mu_H^{-1}\in H$ is often referred to as the formal monodromy (clearly in general it is not conjugate in $G$ to the local monodromy $\mu_G$).

\subsection{Stokes data by level}\label{ssn: stokes by level}

To reorganise the Stokes data according to levels it is visually helpful 
to interpret it in terms of a $G$-local system on a punctured disc.
In particular this helps keep track of the various monodromy relations.
(We will call this the ``punctured disc model'' of an irregular connection.)
Identify $\Delta$ with the unit disc (centred at $z=0$) and puncture it at the origin and at $r$ equally spaced points along each anti-Stokes direction $d\in \IA$, yielding a punctured disc $\Delta'$.
Thus $\Delta'$ is the union of a small punctured disk $\Delta_0$ and concentric annuli $\Ann_1,\Ann_2,\ldots,\Ann_r$ (with increasing diameters) so that each annulus has exactly one puncture in each direction $d\in \IA$.
(We will refer to $\Ann_i$ as the ``level-$i$ annulus''.)
Choose a direction $p\in S^1$ somewhere  between $d_s$ and $d_1$ (i.e. in a small negative sense from $d_1$).
Choose a basepoint $*\in \Delta_0$ in the direction $p$.

Fix a point $(C,h,S_1,\ldots,S_s)\in \cA(Q)$
and recall that each Stokes multiplier $S_i\in \ISto_i$ may be uniquely written 
as 
$$S_i = S_i^1 S_i^2 \cdots S_i^r$$
with $S_i^j\in \ISto_{d_i}(k_j)$.
Thus we may define a homomorphism $\rho :\pi_1(\Delta',*)\to G$ 
from the fundamental group $\pi_1(\Delta',*)$ of $\Delta'$ based at $*$ by the formulae
$$\rho(\ga_i^j) = S_i^j,\qquad
\rho(\ga_0) = h$$ 
where $\ga_0$ is the loop in $\Delta_0$ going once in a positive sense around $0$, and 
the loop $\ga_i^j$ follows a positive arc around $0$ in $\Delta_0$ until just before $d_i$ then follows a ray straight out until $\Ann_j$, before making a small positive loop around the $i$th puncture on $\Ann_j$ and then retracing its steps back to $*$.
The map $\rho$ is well-defined since these loops freely generate $\pi_1(\Delta',*)$.

Now let $\ga_r$ be the loop which goes out in direction $p$ to the outer boundary of $\Delta$, goes around the boundary once in a positive sense then returns to $*$ in the direction $p$.
One may then check that we have
\begin{lem}
\beq\label{eq: orig mon prod}
\rho(\ga_r) = hS_s\cdots S_2S_1.
\eeq
\end{lem}
\pf 
For each $i$ the product $\wh \ga_i:=\ga_i^1\cdots\ga_i^r$ goes around to just before $d_i$ out to $\partial \Delta$ around a positive arc just crossing $d_i$, back in to $\Delta_0$ then back along a negative arc to $*$. Thus $\rho(\wh \ga_i) = S_i$. Next we observe that $\ga_r = \ga_0 \wh \ga_s\cdots \wh \ga_1$, and so the result follows.\epf

Next note that there is a nested  chain of connected complex reductive subgroups 
\beq\label{eq: chain of levis}
H=H_1\subset H_2\subset \cdots H_r \subset G
\eeq
defined by $H_i = \Stab(A_r,\ldots,A_i)$,
each a Levi subgroup of a parabolic of $G$.
Said differently 
for each index $i=r, r-1,\ldots,2,1$  there is a vector space decomposition
$$\h_{i+1} = \h_{i} \oplus  \h_i'$$ %
where $\h_i=\Lie(H_i)$ with $\h_{r+1} = \g$ by convention and 
$\h_i' = \Im(\ad(A_i)\bigl\vert_{\h_{i+1}})
= [A_i,\h_{i+1}]\subset\h_{i+1}$.
Thus 
\beq\label{eq: g decomp}
\g = \h \oplus \h_1' \oplus \h_2' \oplus\cdots \oplus \h_r'
\eeq
and each root space $\g_\al$ occurs in precisely one such component, and
if $\g_\al\subset \h'_i$ then $\deg(q_\al)=k_i$.
(In brief $\h_r'$ is the sum of the root spaces on which $A_r$ acts nontrivially, then $\h_{r-1}'$ is the sum of the remaining root spaces on which $A_{r-1}$ acts nontrivially, etc.)

\begin{lem}
If $\ga\in \pi_1(\Delta',\ast)$ is a loop which does not stray into $\Ann_j$
then $\rho(\ga) \in H_j$.
\end{lem}
\pf
We have $h\in H_1$ so it is enough to check $S_i^j \in H_{j+1}$ for all $i,j$. Now $S_i^j\in \ISto_i(k_j)$ whose Lie algebra is spanned by root spaces $\g_\al$ 
of certain roots $\al$ with $\deg(q_\al)= k_j$. Thus 
$\g_\al\subset \h'_j\subset \h_{j+1}$ and so $\ISto_i(k_j)\subset H_{j+1}$ as required.
\epf

Now we will pass to new generating loops of the fundamental group (still based at $*$).
Define loops $\be_i^j$ going out in the direction $p$ until the inner boundary of $\Ann_j$, then around a positive arc until $d_i$, then around a small positive loop around the $i$th puncture of $\Ann_j$, then retracing the same arc back around to $p$ then back in to $*$.
Then we define $$B_i^j = \rho(\be_i^j)\in G$$
to be the corresponding ``twisted'' Stokes multiplier.

\begin{lem} The twisted Stokes multiplier $B_i^j$ is in $\ISto_i(k_j)$.
\end{lem}
\pf
In the fundamental group we have $\be_i^j = x^{-1}\ga_j^i x$ where $x$ is a loop not entering $\Ann_j$. (Explicitly 
$$x= \wh \ga^j_{i-1} \cdots \wh \ga^j_{2} \wh \ga^j_{1}$$
where      
$\wh\ga_k^j = \ga_k^1\ga_k^2\cdots\ga_k^{j-1}$.) Thus $\rho(x)\in H_j$. 
Now we observe that $H_j$ normalises  $\ISto_i(k_j)$ in $G$, i.e. 
$g\ISto_i(k_j)g^{-1} = \ISto_i(k_j)$ for any $g\in H_j$, which immediately implies the desired result, as $S_i^j\in \ISto_i(k_j)$.
To see this note that $H_j$ is generated by $T$ and the root groups $U_\al$ for roots $\al$ with $\g_\al\subset \h_j$. We must check each $U_\al$ normalises $\ISto_i(k_j)$ (since it is clear for $T$).
For this, by \cite{Bor91} Proposition 14.5 (3), it is enough to check that 
any root of the form $\ga = n\al + m \be$ where $n,m\in\IZ_{>0}$ and 
$\be\in\cR(d_i,k_j)$ is actually in  $\cR(d_i,k_j)$.
This is clear however since $q_\al$ will have lower degree than $q_\be$,
and so the leading term of $q_\ga$ is $m$ times that of $q_\be$, which implies immediately that $\ga \in \cR(d_i,k_j)$.
\epf

Thus specifying all the Stokes multipliers 
$S_i^j\in\ISto_i(k_j)$ is equivalent to specifying all the twisted Stokes multipliers $B_i^j\in\ISto_i(k_j)$.

\begin{lem}
In terms of the twisted Stokes multipliers the monodromy around the outer boundary of $\Delta$ is
\beq \label{eq: mon prod}
\rho(\ga_r) = h (B^1_s\cdots B^1_1) (B^2_s\cdots B^2_1) \cdots
(B^r_s\cdots B^r_1).
\eeq
\end{lem}
\pf
Let $\ga_i$ denote the loop which goes out along $p$ to the outer boundary of $\Ann_i$, around this boundary circle in a positive sense, then back to $*$ along $p$.
The result will follow immediately from the inductive step:
$$\rho(\ga_i) = \rho(\ga_{i-1})B^i_s\cdots B_1^i.$$
In turn this is easily established by drawing a picture of $\Ann_i$.
\epf

Now if we define 
$$h_i = \rho(\ga_{i-1}) = h (B^1_s\cdots B^1_1) (B^2_s\cdots B^2_1) \cdots
(B^{i-1}_s\cdots B^{i-1}_1)$$
for each $i$, then $h_i\in H_i$ (since the corresponding loop does not enter $\Ann_i$).
Now define
$$
\cA(i):= 
H_{i+1}\times H_i\times  \Prod_{j=1}^s \ISto_j(k_i),$$
and assume for the moment the following lemma.

\begin{lem}\label{lem: qh str on one level}
$\cA(i)$ is a quasi-Hamiltonian $H_{i+1}\times H_i$-space
with moment map
$$\mu_i(C_i,h_i,B^i_1,\ldots, B^i_s) = 
( C_i^{-1}h_iB^i_s\cdots B^i_1C_i,\  h_i^{-1})$$
where $C_i\in H_{i+1}$ etc. 
\end{lem}
\noindent
Then  Theorem \ref{thm: main qh thm} follows from the following proposition.
\begin{prop}
There is an explicit isomorphism
$$
\cA(Q)\quad\cong \quad \cA(1)\glue{H_2}\cdots \glue{H_r} \cA(r)
$$
and so $\cA(Q)$ is a quasi-Hamiltonian $G\times H$-space.
\end{prop}
\pf
The action of  $(k,g)\in H_i\times H_{i+1}$ on $\cA(i)$ is given by
$$(k,g) (C_i,h_i,B^i_1,\ldots, B^i_s) = 
(kC_ig^{-1},kh_ik^{-1},kB^i_1k^{-1},\ldots, kB^i_sk^{-1}).$$
Since the action of $H_i$ on $\cA(i-1)$ is free, 
we can identify the gluing  $\cA(i)\glue{H_i} \cA(i-1)$
with the subvariety of the product $\cA(i)\times \cA(i-1)$ 
where $C_{i-1}=1$ and 
$$h_i = h_{i-1}B_s^{i-1}\cdots B_1^{i-1}.$$
Thus we may remove both factors of $H_i$ in the product to see 
$\cA(i)\glue{H_i} \cA(i-1)$ is isomorphic to 
$$H_{i-1}\times H_{i+1}\times \Prod_{j=1}^s \ISto_j(k_{i-1})\times\ISto_j(k_i)$$
which is thus a quasi-Hamiltonian $H_{i+1}\times H_{i-1}$-space with moment map 
$$ 
\left(C_i^{-1}h_{i-1}(B^{i-1}_s\cdots B^{i-1}_1)(B^i_s\cdots B^i_1)C_i,
\  h^{-1}_{i-1}\right).$$
Repeating this gluing process yields the result.
Note that the moment map on $\cA(Q)$ is
$(C^{-1}\rho(\ga_r)C , h^{-1})$ where $C=C_r\in G, h=h_1\in H$ and 
$\rho(\ga_r)$ is as in \eqref{eq: mon prod}, in terms of the twisted Stokes multipliers, or as in \eqref{eq: orig mon prod} in terms of the original Stokes multipliers.
\epf

The fact that the quasi-Hamiltonian two-form is given by the formula in the statement of Theorem \ref{thm: main qh thm} now follows from repeated use of the nesting result (Theorem \ref{thm: nesting}).
Finally Lemma \ref{lem: qh str on one level} will be established.

\pfms (of  Lemma \ref{lem: qh str on one level}).
This is true since the spaces $\cA(i)$ are isomorphic to higher fission spaces;
there is an isomorphism
$$\cA(i) \cong \papk{H_i}{H_{i+1}}{k_i}$$
as follows.
First observe that $\cA(i)$ is isomorphic to the space $\cA(Q_i)$ attached to $Q_i=A_i/z^{k_i}$ where $A_i$ is viewed as an element of the Cartan subalgebra $\lt$  of $\lh_{i+1}$.
To see this note that if $\deg(q_\al)= k_i$ then 
$q_\al$ equals $\al\circ Q_i$ plus less singular terms, and so the nontrivial Stokes groups occurring in the definition of $\cA(i)$ are precisely those occurring in $\cA(Q_i)$.

Thus we can reduce to the case of one level with $Q=Q_i, k=k_i, A=A_i, H=C_G(A)$ etc, and the aim is to show $\cA(Q) = \papk{G}{H}{k}$.
This is just (the parabolic extension of)
\cite{bafi} Lemma 2.4, Lemma A.3, and the proof is the same.
Here is the idea, for completeness.

Recall the singular directions may be described as follows:
each $d\in \IA$ is supported by some root $\al\in \cR$, and 
$d\in \IA$ is supported by $\al$ if and only if
$q_\al=\al(A)/z^{k}$
is real and $<0$ (for $z$ in the direction $d$).
Thus $\IA$ is invariant by rotation by $\pi/k$, 
and so $l:=s/2k=\#\IA/2k$ is integral.    
Define a ``half-period'' to be an ordered $l$-tuple of consecutive singular 
directions in $\IA$. 
\begin{lem}
If $\bd \subset \IA$ is a half-period 
then the subgroups $\{\ISto_d \st d\in \bd\}$ directly span  the unipotent radical of a parabolic subgroup of $G$ with Levi subgroup $H$ (and
rotating $\bd$ by $\pi/k$ yields the unipotent radical of the opposite parabolic with Levi subgroup $H$).
\end{lem}
\pf
The standard fact we wish to use is that any element
$\la\in \lt_\IR$ 
determines a parabolic subgroup $P_\la\subset G$ by defining
$$P_\la = \{g\in G\st z^\la g z^{-\la}\text{ has a limit as $z\to 0$ along any ray}\}.$$
Equivalently $P_\la$ is generated by its Levi subgroup and its unipotent radical
$U_\la$, and in turn the
 Levi subgroup of $P_\la$ is generated by $T$ and the root groups 
$U_\al$ for all the roots with $\al(\la)=0$, and 
$U_\la$ is generated by the 
root groups $U_\al$ for all the roots $\al$ with $\al(\la)>0$.
Now  let 
$\th(\bd)$ be the ray bisecting the sector spanned by the half-period $\bd$ 
and take 
$$\la = -\Re (Q(z))\in \lt_\IR$$
for any $z\ne 0$ on the ray $\theta(\bd)$.
Then one notices that, for each root $\al$, the ``sine-wave function'' 
$f_\al(\phi) = -\Re(q_\al(z))\bigl\vert_{z=\exp(i\phi)}:S^1\to \IR$ 
is either identically zero, or has period $2\pi/k$ and is maximal 
on singular directions supported by $\al$.
Thus $f_\al(\th(\bd)) =\al(\la)$ is strictly positive if and only if 
there is a singular direction supported by $\al$ within  $\pi/2k$ of $\theta(\bd)$.
Thus the roots supporting the directions $d\in \bd$ are precisely those whose root groups generate $U_\la$. 
Further the Levi subgroup of $P_\la$ is just the centraliser of $A$.
\epf

Thus taking $\bd = (d_1,\ldots,d_l)$ yields a parabolic $P_+\subset G$ with Levi subgroup $H$. Denote its unipotent radical $U_+$, and let $U_-$ be the unipotent radical of the opposite parabolic, associated to $(d_{l+1},\ldots,d_{2l})$.
Thus, in this one level case, we have a direct spanning equivalence
$$ \cA(Q) = G\times H\times \Prod_1^s \ISto_j \to 
\papk{G}{H}{k} = G\times H\times (U_+\times U_-)^k$$
induced by the product isomorphisms
$$
\ISto_{(n+1)l}\times\cdots \times\ISto_{nl+2} \times\ISto_{nl+1}
\to U_+$$
for $n\ge 0$ even (and to $U_-$ for $n$ odd). 
Thus by Corollary \ref{cor: dseq} $\cA(Q)$ is quasi-Hamiltonian, as desired.
\epfms

\section{Irregular curves and associated Betti spaces}\label{sn: irreg curves}

In this section we will define the notion of an irregular curve and show how to associate a quasi-Hamiltonian space to an irregular curve with some tangential basepoints, and in turn how to canonically associate a Poisson variety to an irregular curve. 

\subsection{Irregular curves}

Fix a connected complex reductive group $G$ and a maximal torus $T\subset G$.

\begin{defn}\label{defn: irreg curve}
An ``irregular curve'' (or ``wild Riemann surface'') is a smooth compact Riemann surface  $\Si$ (possibly with boundary) together with a finite number of distinct marked points $a_1,a_2,\ldots$ in the interior of $\Si$, 
and an irregular type $Q_i$ (in the sense of Definition \ref{def: irtype}) at each marked point.
\end{defn}

For example if the boundary is empty and each irregular type is zero then an irregular curve is essentially the same thing as a smooth complex algebraic curve with some ordered marked points.
(In general we will say an irregular curve is {\em algebraic} if its boundary is empty---it may still have marked points.)
If $m$ denotes the number of marked points plus the number of boundary components, we will always assume $m>0$.

Given an irregular curve $\Si$ let $\wh \Si\to \Si$ denote the real 
two-manifold with boundary obtained by taking the real 
oriented blow up of $\Si$ at each marked point, i.e. replacing each marked point $a_i$ with the circle of oriented real tangent directions at $a_i$. 
Label the boundary circles of $\wh \Si$ as $\partial_1,\ldots, \partial_m$.
Thus
$Q_i$ determines a subgroup $H_i=C_G(Q_i)\subset G$, 
singular directions
$\IA_i\subset \partial_i$ and Stokes groups 
$\ISto_d\subset G$ for each $d\in \IA_i$, as in Definitions \ref{def: sing dirn} and \ref{def: stokes gp}
(where we set $Q_i=0$ if  $\partial_i$ was already a boundary component of $\Si$).

Now puncture $\wh \Si$ once in its interior near each singular direction
$d\in \IA_i, i=1,\ldots, m$, 
and 
 draw small cilia (eyelashes)
on the surface $\wh \Si$ between each puncture and the corresponding 
singular direction $d\in \IA_i\subset \wh \Si$ 
(such that none of the cilia cross).
The cilia are just to help keep track of the punctures. 
Let $\wt \Si\subset \wh \Si$ denote the corresponding punctured surface
(see Figure \ref{fig: irreg curve}).

\begin{figure}[h]
	\centering
	\input{irreg.curve7.1.pstex_t}
	\caption{The surface $\wt\Si$.}\label{fig: irreg curve}
\end{figure}

Now choose a marked point $b_i\in \partial_i$ in each boundary component of $\wt \Si$ and
define $\Pi$ to be the fundamental groupoid of 
$\wt \Si$ based at $\{b_1,\ldots,b_m\}$:
$$\Pi = \Pi_1(\wt \Si, \{b_1,\ldots,b_m\}),$$
consisting of homotopy classes of paths $\ga$ in $\wt \Si$
such that the endpoints of $\ga$ are in the set 
$\{b_1,\ldots,b_m\}$.
(More precisely the set of objects of $\Pi$  is the finite set $\{b_1,\ldots,b_m\}\subset \wt\Si$ and a morphism in $\Pi$ from 
$b_i$ to $b_j$ is a homotopy class of paths in $\wt\Si$ starting at $b_i$ and ending at $b_j$.)

Thus we may consider the space $\Hom(\Pi,G)$ of morphisms from the groupoid $\Pi$ to the group $G$.
(Recall a group is a groupoid with only one object and by `morphism' we mean a functor between the two corresponding categories.)
Explicitly an element $\rho\in \Hom(\Pi,G)$ consists of a choice of an element $\rho(\ga)\in G$ for each path $\ga$ in $\Pi$, such that for composable paths $\ga_1,\ga_2$ one has
$$\rho(\ga_1\circ \ga_2) = \rho(\ga_1)\rho(\ga_2).$$
We will also refer to $\rho$ as a ``$G$-valued representation of $\Pi$''.
Now consider the subspace
$$\Hom_\IS(\Pi,G)\subset \Hom(\Pi,G)$$
of ``Stokes representations'' $\rho$ obeying the following conditions (for any $i=1,\ldots,m$):

SR1) if $d\in \IA_i$ and $\wh \ga_d$ is any loop based at $b_i$ which goes around (in any direction) $\partial_i$ to the direction $d$, and then  loops once around the puncture on the cilium emanating from the direction $d$ (without crossing any other cilia), before retracing its path back to $b_i$, then $\rho(\wh\ga_d)\in \ISto_d$, and

\label{page: SR1,2 v2}

SR2) if $\ga_i$ is the simple closed loop based at $b_i$ going once in a positive sense around $\partial_i$
then $\rho(\ga_i) \in H_i$.

There is an action of the group $\bH:=H_1\times \cdots \times H_m\subset G^m$ on the space of Stokes representations, as follows: 
$(k_1,\ldots,k_m)\in \bH$ sends $\rho$ to the representation $\rho'$ such that 
$$\rho'(\ga) = k_j \rho(\ga) k_i^{-1}$$
for any path $\ga\in \Pi$ from $b_i$ to $b_j$.

\begin{thm}\label{thm: qh str on streps}
The space $\Hom_\IS(\Pi,G)$ of Stokes representations of $\Pi$ in $G$ is a smooth complex affine variety and is (canonically) a quasi-Hamiltonian $\bH$-space, where 
$\bH=H_1\times \cdots \times H_m\subset G^m$.
\end{thm}
\pf
First we will establish this in the case when $\Si$ is a disk with one marked point in its interior.
In this case, if we choose some paths generating $\Pi$, as for example in Figure \ref{fig: Pi gens}, with $\partial_1$ the inner boundary and 
$\partial_2$ the outer boundary, and number  the singular directions $d_1,d_2,\ldots, d_s$  correspondingly, then we obtain an isomorphism
 $$\Hom_\IS(\Pi,G)\cong\cA(Q_1)=G\times H_1\times \Pi_1^s \ISto_{d_i}.$$
This gives $\Hom_\IS(\Pi,G)$ a quasi-Hamiltonian structure 
(by Theorem \ref{thm: main qh thm}) and we should check it is independent of the choice of generating paths.
Draw a concentric circle $\IH$ through all the punctures. 
We only consider sets of generating paths of the following form 
1) a path $\cP$ from $b_2$ to $b_1$, only crossing $\IH$ once, between two punctures, which will then be labelled $d_s,d_1$,
2) a simple loop $\ga_0$ based at $b_1$ around the inner boundary,
3) loops $\wh \ga_d$ based at $b_1$ around the puncture in the direction $d$, so that none of the paths cross (so the paths look as in Figure \ref{fig: Pi gens}).
Upto homotopy the only choice here is the path $\cP$, and the choice for this path is the number of times it winds around before and after crossing $\IH$, and the choice of where it crosses $\IH$. The quasi-Hamiltonian structure is independent of these choices due to the automorphisms \eqref{eq: inner twist} and \eqref{eq: outer twist} (enabling one to undo the winding of $\cP$ on the inside and outside of $\IH$), and the isomonodromy isomorphisms (Propositions \ref{prop: imd im},\ref{prop: imd im2})
(enabling one to change the choice of where $\IH$ is crossed).

Now consider the general case. Removing  disks around each marked point reduces to the cases of a disk (already dealt with above) and the case with only trivial irregular types considered in \cite{AMM}.
This yields an intrinsic quasi-Hamiltonian structure.
To see it is a smooth affine variety, note that
upon choosing appropriate generating paths we may identify $\Hom_\IS(\Pi,G)$ with the reduction by $G$ (at the identity value of the moment map) of the fusion product
\beq\label{eq: big fusion prod}
\ID^{\fus g}\fusion{G}\cA(Q_1)
\fusion{G}\cdots\fusion{G} \cA(Q_m),
\eeq 
where $\ID$ is the internally fused double and $g$ is the genus of $\Si$.
Let $\mu_G$ denote the $G$ component of the moment map, 
from  \eqref{eq: big fusion prod} to $G$.
Since $m\ge 1$ the action of $G$ on \eqref{eq: big fusion prod} is free and so $1$ is a regular value of $\mu_G$ 
(by \cite{AMM}  Proposition 4.1 (3)).
Thus
$\mu_G^{-1}(1)$ is a smooth subvariety of \eqref{eq: big fusion prod} with a free action of $G$ and so the quotient is smooth (by Luna's slice theorem or otherwise).
(The quasi-Hamiltonian reduction theorem applies to yield the structure of quasi-Hamiltonian $\bH$-space.)
Explicitly, if we write an element of 
$\cA(Q_i)=G\times H_i\times \ISto(Q_i)$ as $(C_i,h_i,\bS^i)$, then
we may identify $\mu_G^{-1}(1)/G$ with the 
affine subvariety of \eqref{eq: big fusion prod}
cut out by the equations $\mu_G=1$ and $C_1=1$.
In detail if we write $\ID^{\fus g}$ as 
$\{(a_i,b_i)\st a_i,b_i\in G, i=1,\ldots,g\}$,
then the relation
$\mu_G=1$ 
takes the form
\beq\label{eq: monod reln}
[a_1,b_1]\cdots[a_g,b_g] \mu_1\cdots \mu_m = 1
\eeq
where $\mu_i = C_i^{-1}h_i\cdots S^i_2S^i_1C_i$
and $[a,b]=aba^{-1}b^{-1}$. 
\epf

Recall (from Remark \ref{rmk: ab implies symp}) that 
if $\bH$ is a torus it follows that $\Hom_\IS(\Pi,G)$ 
is an algebraic symplectic manifold.
The special case of this when the leading coefficient of each irregular type
is regular semisimple (so each chain of reductive groups \eqref{eq: chain of levis} passes directly from $G$ to $T$) was considered in \cite{saqh}. 
In the simple example of $\Si=\IP^1$ with two marked points, and $Q_1,Q_2$ both with simple poles and regular semisimple coefficients, 
a covering of $\Hom_\IS(\Pi,G)$ 
was identified with the Lu--Weinstein double symplectic groupoid in \cite{saqh} Proposition 7.

A priori the space $\Hom_\IS(\Pi,G)$ depends on the irregular curve $\Si$, the basepoints $\{b_i\}$ and
the choice of the locations of the punctures used to define $\wt \Si$.
By shrinking the cilia, pulling the punctures close to the boundary circles, one may canonically identify the different possible groupoids $\Pi$ defined using different puncture locations. 
The next result shows that once we quotient by $\bH$ the result is independent of the basepoints as well.

\begin{cor}
The irregular curve $\Si$ canonically determines the Poisson variety
$\Hom_\IS(\Pi,G)/\bH$ of S-equivalence classes of Stokes representations.
\end{cor}
\pf
By Proposition \ref{prop: poisson quotients}, upon choosing 
basepoints $\{b_i\}$ the quotient is well defined and a Poisson variety, so we just need to check it is independent of this choice. 
Suppose we make two different choices $\{b_i\}, \{b_i'\}$ of basepoints.
Choosing any path $\ga_i$ in $\partial_i$ from $b_i$ to $b_i'$ (for each $i$) yields an isomorphism $\Hom_\IS(\Pi,G)\cong \Hom_\IS(\Pi',G)$ (where $\Pi'$ is defined using $\{b'_i\}$). Choosing a different set of paths corresponds to conjugating this isomorphism by the action of an element of $\bH$, by SR2), and so the $\bH$-invariant functions are canonically identified. 
\epf

Analogously to Simpson \cite{Sim94b} \S6 we will sometimes refer to the geometric invariant theory quotient $\Hom_\IS(\Pi,G)/\bH$
as the ``Betti moduli space'' of the irregular curve $\Si$, and denote it $\bM_B(\Si)$.
These varieties, or their symplectic leaves, will also sometimes be called ``wild character varieties'' (see \cite{sikora-charvars} and references therein for the usual case).

\begin{rmk}\label{rmk: stokes local system}
By considering Stokes representations we are in effect considering a special class of $G$-local systems on $\wt \Si$, the {\em Stokes $G$-local systems}, defined as follows.
Draw concentric circles (halos) $\IH_i$ on $\wt \Si$ near each boundary circle through the corresponding punctures, as in the proof of Theorem \ref{thm: qh str on streps}, thus determining a small annulus around each boundary circle, which we will refer to as the area {\em inside} of $\IH_i$.
Then a ``Stokes $G$-local system'' on $\wt \Si$ is a $G$-local system $\IL$ on $\wt \Si$ together with a flat reduction of structure group to $H_i$ inside of $\IH_i$ for each $i=1,\ldots,m$ (i.e. an $H_i$-local system $\IL_i$  defined inside $\IH_i$, such that $\IL = \IL_i \times_{H_i} G$ there), such that, for any basepoint inside $\IH_i$, the local monodromy around the puncture corresponding to $d\in \IA_i$ lies in $\ISto_d(Q_i)$.
The space $\Hom_\IS(\Pi,G)$ classifies Stokes $G$-local systems together with a framing of $\IL_i$ at the basepoint $b_i$ for each $i$.
Section \ref{ssn: stokes by level} may now be viewed as giving 
several equivalent 
descriptions of the category of Stokes $G$-local systems. 
\end{rmk}

\begin{rmk}\label{rmk: finite quotients}
Note that  the Betti space $\bM_B(\Si)$ is independent of the labelling of the boundary components $\partial_i$ by integers $1,\ldots,m$.
Further it is straightforward to verify that $\bM_B(\Si)$ only depends on each irregular type $Q_i$ up to the action of the Weyl group $W=N_G(T)/T$ of $G$.
More precisely if for some $i$ the irregular type $Q_i$ is replaced by $Q_i'=\Ad_w(Q_i)$ for some $w\in N_G(T)$, 
then there is a canonical isomorphism between the corresponding Betti spaces, independent of the choice of $w$. 
\end{rmk}

\begin{rmk}
More generally one may consider  ``twisted'' irregular types, defined by replacing the Cartan subalgebra $\lt\flp z \frp\subset \g\flp z \frp$ by a  nonconjugate Cartan subalgebra (which exist since $\IC\flp z \frp$ is not algebraically closed, cf. \cite{kazlus88} Lemma 2), and this yields a notion of ``twisted'' irregular curves, which will be studied elsewhere. 
Further one may replace the constant group $G$ by a local system of groups 
on $\Si\setminus\{a_i\}$ (e.g. as for twisted loop groups)---this yields a more general notion of twisted irregular curve. %
\end{rmk}

\section{Stability of Stokes representations}

Given an irregular curve $\Si$ and some basepoints $\{b_i\}$
we have defined a smooth affine variety $\Hom_\IS(\Pi,G)$
with an action of a reductive group $\bH$.
This is a situation much studied in geometric invariant theory  and one defines the ``stable'' points as follows.
Let $\bK\subset \bH$ be the subgroup of 
elements which act trivially on all points of $\Hom_\IS(\Pi,G)$ (the kernel of the action).

\begin{defn}(see e.g. \cite{Richardson}).
A point $\rho\in \Hom_\IS(\Pi,G)$ is stable if its orbit $\bH\cdot \rho$ is closed and of dimension equal to $\dim (\bH/\bK)$.
\end{defn}

In this section we will assume that 
$\bK$ 
equals the centre of $G$ embedded diagonally in $\bH$. 
(One may check this is the case unless the genus of $\Si$ is zero, 
there is just one marked point and $Q_1$ has at most a simple pole---in such exceptional cases $\Hom_\IS(\Pi,G)$ is a point.)

Our first aim is to characterise the stable Stokes representations $\rho$ in a more direct fashion, as follows.
Suppose we have a parabolic subgroup $P_i\subset G$ for each basepoint $b_i$, $i=1,\ldots,m$, and write $\bP=(P_1,\ldots,P_m)$
for this collection of parabolic subgroups.
We will say that $\bP$ is {\em compatible} with $\rho\in \Hom_\IS(\Pi,G)$ if
$$  \rho(\ga) P_i \rho(\ga)^{-1} = P_j$$
for any path $\ga$ in $\wt \Si$ from $b_i$ to $b_j$ (for any $i,j$).
Now let $Z_i\subset G$ be the identity component of the centre of $H_i$. Thus $Z_i$ is a torus and the group $H_i$ may be characterised as the centraliser of $Z_i$ in $G$. 
We will say that $\bP$ is {\em invariant} if $Z_i\subset P_i$ for each $i$. Finally $\bP$  is {\em proper} if some $P_i$ is a proper subgroup of $G$.

\begin{defn}
A representation $\rho\in \Hom_\IS(\Pi,G)$ is ``reducible'' if 
there is an 
invariant proper collection of parabolics compatible with $\rho$. Otherwise $\rho$ is ``irreducible''.
\end{defn}

In this definition it makes no difference if we only consider maximal proper parabolic subgroups: 
$\rho$ is reducible if and only if there is an 
invariant collection of parabolics $\bP$ compatible with $\rho$, such that some (and hence all) $P_i\subset G$ is a maximal proper parabolic.
We will establish the following:

\begin{thm}\label{thm: stable iff irred}
A point $\rho\in \Hom_\IS(\Pi,G)$ is stable if and only if it is irreducible.
\end{thm}

If there are no irregular singularities (each $Q_i=0$) this follows from Theorem 4.1 of Richardson \cite{Richardson}
who considered the diagonal conjugation action of $G$ on $G^N$.
On the other hand if there is just one singularity $m=1$ this follows from \cite{Richardson} Theorem 14.1,
considering the diagonal conjugation action of $H_1\subset G$ 
on $G^N$.
Further if $G$ is a general linear group one can translate this into a problem involving quivers and appeal to King \cite{king-quivers}.
The general case however does not seem to follow from the results of either of \cite{king-quivers, Richardson}, but as in those articles the result is again essentially ``an exercise in the Hilbert--Mumford theorem''.

\pf
If $H$ is any complex algebraic group let $Y(H)$ denote the set of one parameter subgroups $\la:\IC^*\to H$.
Recall that any $\la\in Y(G)$ determines a parabolic subgroup 
$$P_G(\la) = \{g\in G \st \lim_{t\to 0} \la(t)g\la(t)^{-1} \text{ exists}\}\subset G.$$ 
The Hilbert--Mumford theorem implies (see \cite{king-quivers} Proposition 2.5) that
$\rho \in \Hom_\IS(\Pi,G)$ is stable if and only if
whenever $\la\in Y(\bH)$ and $\lim_{t\to 0} \la(t)\cdot \rho$ exists, then $\la\in Y(\bK)$.
Thus given $\la\in Y(\bH)$ and $\rho$ it is important to determine if $\lim_{t\to 0} \la(t)\cdot \rho$ exists.
Write $\la=(\la_1,\ldots,\la_m)$ with $\la_i\in Y(H_i)$.
Since $H_i\subset G$ there are parabolic subgroups  
$P_i:=P_G(\la_i)$ of  $G$ for each $i$. 
Clearly if the limit exists then $\rho(\ga)$ should be in 
$P_i$ for any loop $\ga$ based at $b_i$ (since $\la$ acts on 
$\rho(\ga)$ by conjugation by $\la_i$), e.g. for each Stokes multiplier at $a_i$ and the formal monodromy $h_i$.
Now suppose $C=\rho(\ga)\in G$ for some path between two distinct basepoints, say from $b_2$ to $b_1$. 

\begin{lem}\label{lem: lim of C}
Suppose $C\in G$ and $\la_i\in Y(H_i)$ for $i=1,2$.
Then $\la_1(t) C \la_2(t)^{-1}$ has a limit as $t\to 0$ if and only if
$\la_1$ and $\la_2$ are conjugate in $G$ and 
$C P_2 C^{-1} = P_1$.
\end{lem}
\pf
If they are conjugate, say $\la_2 = g\la_1 g^{-1}$, then the desired limit exists if and only if $C P_2C^{-1} = P_1$.
Indeed $C P_2C^{-1} = CgP_1g^{-1}C^{-1}$ so that
 $C P_2C^{-1} =P_1$ if and only if $Cg\in P_1$, i.e. if and only if
$\la_1(t)Cg\la_1^{-1}(t)$ has a limit as $t\to 0$. Thus, multiplying on the right by the constant $g^{-1}$, this has a limit if and only if  $\la_1(t) C \la_2(t)^{-1}$ has a limit as $t\to 0$.
Conversely if the limit exists, and equals $C_0\in G$ say, then, since the limit is a fixed point, 
$\la_1(t)C_0\la_2(t)^{-1} = C_0$ for all $t$, and so $\la_1$ and $\la_2$ are conjugate. Then as above, because the limit exists, $C P_2C^{-1} = P_1$.
\epf

Note that if $P_i=P_G(\la_i)$ for some $\la_i\in Y(H_i)$ then 
$Z_i\subset P_i$. Indeed if $z\in Z_i$ and $g\in P_i$ then 
$zgz^{-1}\in P_i$ (since $z$ and $\la_i$ commute) so 
$z\in N_G(P_i) = P_i$.

\begin{prop}
Given $\rho\in \Hom_\IS(\Pi,G)$ then $\rho$ is reducible
if and only if there exists
$\la\in Y(\bH)$ such that $\la\not\in Y(\bK)$ and the limit 
$\lim_{t\to 0} \la\cdot\rho$ exists.
\end{prop}
\pf
If $\lim_{t\to 0} \la\cdot\rho$ exists then taking $P_i=P_G(\la_i)$ gives a collection of parabolics $\bP$. 
It is compatible with $\rho$ by Lemma \ref{lem: lim of C}, it is invariant by the remark after Lemma \ref{lem: lim of C}, and it is proper since $\la\not\in Y(\bK)$.

Conversely suppose we are given $\rho$ and a proper collection 
$\bP$ of invariant compatible parabolics.
Then for each $i$ there is a maximal torus $T_i$ of $G$ such that
$$Z_i\subset T_i\subset H_i\cap P_i$$
(for example since $Z_i$ is reductive it is contained in a Levi subgroup $L_i$ of $P_i$, so there is a maximal torus $T_i$ of 
$L_i$ containing $Z_i$---$T_i$ is also maximal in $G$ and 
clearly $T_i\subset H_i=C_G(Z_i)$).
Choose a Borel subgroup $B$ such that
$T_1\subset B\subset P_1$, i.e. so that $P_1$ is a standard parabolic (for this choice of $T_1$ and $B$).
Then we may choose $\la_1\in Y(T_1)$ so that $P_1=P_G(\la_1)$.
Now, due to the compatibility condition $P_i$ is conjugate in $G$ to $P_1$, and due to the conjugacy of maximal tori of $P_1$ (\cite{Bor91} 11.3) we may simultaneously conjugate 
the pair $T_i\subset  P_i$ to the pair
$T_1 \subset P_1$. 
Thus for each $i\ge 2$ we may  conjugate $\la_1$ to an element 
$\la_i\in Y(T_i)$ such that $P_G(\la_i)=P_i$.
Hence we have constructed $\la=(\la_1,\ldots,\la_m)\in Y(\bH)$
such that  $\la(t)\cdot\rho$ has a limit as $t\to 0$ (via Lemma \ref{lem: lim of C}). 
Moreover $\la\not\in Y(\bK)$ since each $P_i$ is proper. 
\epf

The result is now immediate from the Hilbert--Mumford theorem.
\epf

\subsection{Stability and differential Galois groups}\label{ssn: ssc I}

\begin{cor}
Suppose $\rho$ is a Stokes representation corresponding to a meromorphic connection $A$ on a $G$-bundle on $\Si^\circ$ (as in \S\ref{ssn: global stokes reps from cons}). 
Then $\rho$ is stable if and only if the differential Galois group $\Gal(A)\subset G$ of $A$ is not contained in any proper parabolic subgroup of $G$.  
\end{cor}

\pf
Define $G(\rho)\subset G$ to be the Zariski closure  of the subgroup of $G$ generated by the elements: 1) $\rho(\ga)$ for any loop $\ga$ in $\wt \Si$ based at $b_1$, and 2) $\rho(\ga_i)^{-1} t_i \rho(\ga_i)$
for any path $\ga_i$ in $\wt \Si$ from $b_1$ to $b_i$ and any element $t_i\in Z_i$
for any $i=1,\ldots,m$. 
By Theorem \ref{thm: stable iff irred} $\rho$ is not stable if and only if $G(\rho)$ is a subgroup of a proper parabolic  $P_1\subset G$.
The Ramis--Schlesinger density theorem (cf. \cite{ramis-factn}, 
\cite{L-R94} Theorem III.3.11, \cite{MR91} Theorem 21) says that $\Gal(A)$ is the Zariski closure of the subgroup of $G$
defined in the same way but with each torus $Z_i$ replaced by
the  Ramis exponential torus $\IT_i\subset T$ associated to $Q_i$.
One may verify that $\IT_i$ is characterised as the smallest subtorus of $T$
whose Lie algebra contains all the coefficients of $Q_i$.
Thus $\IT_i\subset Z_i$ (as $Z_i$ has this property) so that $\Gal(A)\subset G(\rho)$ and hence one direction of the corollary is clear.
Conversely it is sufficient to verify that if $\IT_i\subset P$ then $Z_i \subset P$, for any parabolic subgroup $P\subset G$.
To see this first note that $C_G(\IT_i)=C_G(Z_i)=H_i$ (since 
$\IT_i\subset Z_i$, and conversely if $g\in C_G(\IT_i)$ then $g$ centralises the Lie algebra of $\IT_i$ and thus all the coefficients of $Q_i$).
Thus if $T_\mu\subset G$ is any maximal torus containing $\IT_i$
then $T_\mu$ is a maximal torus of $H_i$ and so $Z_i\subset T_\mu$.
Thus if $\IT_i\subset P$ we can take $T_\mu$ to be in a Levi subgroup of 
$P$ containing $\IT_i$ and deduce $Z_i\subset T_\mu\subset P$.
\epf

\subsection{Sufficient stability conditions}\label{ssn: ssc 2}

Recall that for each marked point $H_i\subset G$ is a reductive group containing the maximal torus $T$, so that 
$\bH\subset G^m$ is a reductive group with maximal torus 
$\bT:=T^m$.
A conjugacy class $\cC\subset\bH$ is the same thing
as a product $\cC_1\times \cdots\times \cC_m$ of 
conjugacy classes $\cC_i\subset H_i$.
In this section we will show that if $\cC$ is sufficiently 
generic and
$\rho\in M:= \Hom_\IS(\Pi,G)$ has $\mu(\rho)\in \cC$,
then $\rho$ is stable, where $\mu:M\to \bH$ is the moment map.
(In the general linear case this is related to the sufficient stability conditions of \cite{wnabh} \S8.
Some aspects are similar to \cite{saqh} \S6.)

Recall the Jordan decomposition (cf. \cite{Bor91}), that any element 
$h\in \bH$ is conjugate to an element of the form 
$t \cdot u$ where $t\in \bT$ and $u\in \bH$ is unipotent and commutes with $t$ (so that any conjugacy class of $\bH$ 
may be specified by choosing such elements $t, u$).

\begin{cor}\label{cor: 1st genericity}
There is a Zariski open subset $\bT^\circ\subset \bT$ such that if
$t\in \bT^\circ$ and 
$\mu(\rho)\in \cC$ and $t\cdot u\in \cC$ for some unipotent 
$u\in \bH$ commuting with $t$,
then $\rho$ is stable.
\end{cor}
\pf
We will actually prove a slightly 
stronger and more precise statement.
Recall any connected complex reductive group $G$ has a finite cover
which is a product $Z(G)^\circ \times [G,G]$ of the identity component of the centre of $G$ and the (semisimple) derived subgroup $[G,G]$,
and there is a homomorphism
$$\pr_G : G \to \bar Z(G)$$
from $G$ onto the torus 
$\bar Z(G) := Z(G)^\circ/(Z(G)^\circ \cap [G,G]).$
(For example if $G=\GL_n(\IC)$ this map is the determinant, onto $\IC^*$, and the finite abelian group $Z(G)^\circ \cap [G,G]$ is the centre of $\SL_n(\IC)$.)
It follows that if $a,b,u\in G$ with $u$ unipotent then 
$\pr_G(u) = 1$ and  $\pr_G(aba^{-1}b^{-1})=1$.
Now recall (from Theorem 
\ref{thm: qh str on streps})
that $M$ is isomorphic to the set of $G$-orbits in the subvariety
$\mu_G=1$ of $\ID^{\fus g}\fus\cA(Q_1)\fus \cdots \cA(Q_m)$,
which is written explicitly in \eqref{eq: monod reln} (we have not set $C_1=1$ here).
Applying $\pr_G$ to both sides of \eqref{eq: monod reln}
implies $\Prod_1^m\pr_G(t_i)=1\in \bar Z(G)$, 
where $t_i\in T$ is the $i$th component of $t\in\bT$,
which is conjugate to the semisimple part of $h_i^{-1}$.
Thus $t$ is in the kernel $K$ of the surjective homomorphism of tori 
\beq\label{eq: main torus hom}
\bT\to \bar Z(G);\qquad t\mapsto \Prod\pr_G(t_i).
\eeq
Now we will define a Zariski  open subset of $K$ %
and show  $\rho$ is stable if $t$ is in this subset.
Suppose $P\subset G$ is a maximal standard proper parabolic subgroup (and so contains $T$)
and choose Weyl group elements $\bar w_i\in N_G(T)/T$ for $i=1,\ldots,m$.
Let $L=P/\Rad_u(P)$ be the Levi factor of $P$, 
so there is a homomorphism
$\pr_L: L \to \bar Z(L)$
onto the torus $\bar Z(L)$ associated to $L$.
Thus we may consider the 
surjective homomorphism 
\beq\label{eq: smaller torus homs}
\bT\to \bar Z(L);\qquad t\mapsto \Prod\pr_L(\bar w_i(t_i)).
\eeq
There are only a finite number of such maps (since the Weyl group and the number of standard parabolics is finite).
Let $\bT^\circ\subset \bT$ denote the complement of the 
kernels of all of the maps \eqref{eq: smaller torus homs}.
Since $P$ is a proper subgroup, $\dim(\bar Z(L))>\dim(\bar Z(G)),$
and so the kernel of each map \eqref{eq: smaller torus homs}
is of smaller dimension than the kernel $K$ of \eqref{eq: main torus hom},
and so $\bT^\circ\cap K$ is Zariski open in $K$. 
The precise statement we will prove is:
\begin{cor}\label{cor: generic ccls}
Suppose  $t\in\bT^\circ\cap K$, i.e. 
$t$ is in the kernel of \eqref{eq: main torus hom}, but not in 
the kernel of any of the maps
\eqref{eq: smaller torus homs}.
Then $\rho$ is stable.
\end{cor}
\noindent(The original statement is correct, but vacuous if the centre of $G$ has positive dimension).
To prove this note that if $\rho$ is not stable then there is a (maximal) proper parabolic $P\subset G$ such that all the elements
$C^{-1}_ih_iC_i, C^{-1}_iS^i_jC_i, a_k,b_k$ are in  $P$.
Using the $G$-action we may assume $P$ is standard (and thus contains $T$).
Let $\pr_P: P  \to \bar Z(L)$
be the map obtained by 
composing the canonical projection $P\to L$ with $\pr_L$. 
Applying $\pr_P$  to the  relation $\mu_G=1$ implies 
$\Prod\pr_P(C^{-1}_i h_i C_i) = 1 \in \bar Z(L)$,
noting that usually $C_i\not\in P$.
Now the semisimple part of $h_i$ is conjugate (in $G$) to $t^{-1}_i$ (recalling that the moment map for $H_i$ is $h^{-1}_i)$. 
Thus the semisimple part of $C^{-1}_i h_i C_i$ is conjugate in $P$
to an element of the form $\bar w_i(t^{-1}_i)\in T$, for some 
 Weyl group element $\bar w_i\in N_G(T)/T$. 
Thus $\Prod\pr_L(\bar w_i(t_i)) = 1 \in \bar Z(L)$, and so
 $t$ is  in the kernel of one of the maps 
\eqref{eq: smaller torus homs}.
\epf

\subsection{Examples of well-behaved quotients}

Let $\cC\subset \bH$ be a semisimple conjugacy class
which is generic in the sense of Corollary 
\ref{cor: generic ccls}.
Let $\mu:\Hom_\IS(\Pi,G)\to \bH$ be the moment map, 
and let $\bG:=\bH/\bK$.

\begin{cor}\label{cor: good quotients}

1) The subvariety $\mu^{-1}(\cC)\subset \Hom_\IS(\Pi,G)$ is a smooth affine variety,

2) There is a saturated open subset 
$\cU\subset \mu^{-1}(\cC)$ such that 
the quotient $\cU/\bH$ is a smooth algebraic symplectic manifold (over which $\cU$ is a principal $\bG$-bundle), and this quotient coincides with the set theoretic quotient,

3) If $G=\GL_n(\IC)$ then we may take $\cU=\mu^{-1}(\cC)$ in 2) so that 
$\mu^{-1}(\cC)/\bH$ is a smooth affine algebraic symplectic manifold.
\end{cor}
\pf
Let $\cC'\subset \bH$ be the inverse conjugacy class (if $h\in \cC'$ then $h^{-1}\in \cC$).
Consider the fusion
\beq \label{eq: fusion w ccl}
\Hom_\IS(\Pi,G)\fus \cC'.
\eeq
Let $\wh\mu$ be the corresponding moment map from 
\eqref{eq: fusion w ccl} to $\bH$.
As a variety \eqref{eq: fusion w ccl} is just the product, and so it is a smooth affine variety, since semisimple conjugacy classes are affine. 
Now consider the affine subvariety $\wh\mu^{-1}(1)$.
It is $\bH$-equivariantly isomorphic to 
$\mu^{-1}(\cC)$, and so every point is stable.
Firstly this implies every orbit is closed and so the geometric invariant theory quotient coincides with the set-theoretic quotient.
Stability also implies $\bG$ acts on $\wh\mu^{-1}(1)$ with finite stabilisers, and so implies the following.

\begin{lem}
$\wh\mu^{-1}(1)$ is a smooth affine variety.
\end{lem}
\pf
First suppose $G$ has finite centre, so $\bK$ is finite and $\bH$ acts with finite stabilisers.
Then, given any $p\in\wh\mu^{-1}(1)$, Proposition 4.1 (3) of \cite{AMM} 
implies $d\wh\mu_p$ is surjective, as in the usual Hamiltonian story.
(The proof in \cite{AMM} is for compact groups, but works provided $\wh\mu(p)\in\bH$ is semisimple, as is the case here.)
So $1$ is a regular value of $\wh \mu$ and the lemma follows.  
In general let $\wt \bG=\bH/\bK^\circ$ where $\bK^\circ$ is the identity component. Then \eqref{eq: fusion w ccl} is also a quasi-Hamiltonian $\wt\bG$ space with moment map $\bar\mu:=\pi\circ \wh \mu$ where $\pi:\bH\to \wt\bG$ is the projection. 
The above argument shows $\bar\mu^{-1}(1)$ is a smooth affine variety.
We claim $\wh\mu^{-1}(1)\cong\mu^{-1}(\cC)$ is the union of some connected components of 
$\bar\mu^{-1}(1)$ and so the result follows.
To establish the claim note 
$(\rho,c)\in \bar\mu^{-1}(1)$ iff $(\rho,z\cdot c)\in \wh\mu^{-1}(1)$
for some $z\in \bK^\circ$, i.e. if and only if $\mu(\rho)=zc$ for some 
$c\in \cC,z\in \bK^\circ$.
Then it follows  as in Corollary \ref{cor: 1st genericity} that $z$ is in the kernel of the map \eqref{eq: main torus hom}. 
But \eqref{eq: main torus hom} restricts to an {\em isogeny} 
$\varphi:\bK^\circ\to \bar Z(G)$ (noting that $Z(G)^\circ\cong \bK^\circ\subset \bT$) 
so $z$ is in the fixed finite abelian group $\ker(\varphi)$.
\epf

Since $\bG$ acts with at most finite stabilisers it follows that  
$\wh\mu^{-1}(1)/\bG=\wh\mu^{-1}(1)/\bH$ is an orbifold.
On the other hand it is any easy consequence (see \cite{drezet-slices} Proposition 5.7) of Luna's slice theorem
\cite{luna-slices} 
 that the subset 
$\cU\subset \wh\mu^{-1}(1)$
of points where $\bG$ acts with trivial stabilisers has the following properties:
1) it is a saturated open subset,
2) $\cU/\bG$ is smooth, and
3) $\cU$ is a principal $\bG$-bundle over $\cU/\bG$ (\'etale locally trivial), and so in particular the action of $\bG$ on $\cU$ is (scheme-theoretically) free (\cite{drezet-slices} p.17).

Finally if $G=\GL_n(\IC)$ we should check that the stabiliser in $\bH$ of any stable Stokes representation is $\bK\cong \IC^*$ (the centre of $G$ embedded diagonally in $\bH$). But this follows easily, as in the world of quiver 
representations (cf. e.g. \cite{yamakawa-mpa} Prop. 2.6), so we leave it as an exercise. \epf

Part 3) is reassuring since the irregular Riemann--Hilbert correspondence and  \cite{wnabh} show that, if nonempty, such spaces are complete hyperk\"ahler manifolds (one may check the complex symplectic forms match up as in \cite{smid, saqh}).

\begin{rmk}
Such results go back at least to Gunning \cite{gunning-lovb} \S9 
in the nonsingular case, where one explicitly differentiates the defining relation (see also Weil \cite{Weil64} which includes punctures). The quasi-Hamiltonian approach avoids this, and extends to the irregular case.
It is modelled on the case of symplectic (Marsden--Weinstein) quotients, which are treated algebraically for example in \cite{cas-slod}.
(The subset $\cU$ of stable representations with stabiliser $\bK$ is the analogue of the ``good'' representations in the usual set-up, \cite{johnson-millson} p.57.)
\end{rmk}

\begin{rmk}
Given this explicit description it is easy to write down a formula for the (complex) dimension of the symplectic manifolds in part 2) or 3) of Corollary \ref{cor: good quotients}, assuming they are nonempty:
Let $r=\dim(T)$ be the rank of $G$, and let $\cR\subset \lt^*$ 
denote the roots of $G$.
Given an irregular type $Q$ we have
$$\dim\cA(Q) =  \dim(G) + \dim(H) +
\sum_{\al\in \cR} \deg(\al\circ Q)$$
where $\dim(H)=r + \#\{\al\in \cR\st \deg (\al\circ Q)=0\}$.
Then given an irregular curve $\Si$ with irregular types $Q_1,\ldots,Q_m$
$$\dim\Hom_\IS(\Pi,G) = (2g-2)\dim(G) +\sum_1^m \dim\cA(Q_i)$$
and  in turn given $\cC\subset \bH=H_1\times\cdots\times H_m$
\beq\label{eq: dimns}
\dim\left(\Hom_\IS(\Pi,G)\spqa{\cC} \bH\right) =  
\dim\Hom_\IS(\Pi,G) + \dim\cC - 2(\dim\bH - \dim Z(G))
\eeq
where $Z(G)$ is the centre of $G$.
For example suppose
$Q_i$ has a pole of order $r_i$
and if $r_i=0$ the corresponding conjugacy class $\cC_i\subset G$ is regular semisimple and suppose the nonzero $Q_i$ have regular leading term
(as in \cite{saqh}). 
Then \eqref{eq: dimns} equals
\beq
(2g-2)\dim(G) + 2\dim Z(G)+(\dim(G)-r)(m+\sum r_i).
\eeq
Upon specialising further to $g=0$ and $G=\GL_2(\IC)$ this equals 
$2(m+\sum r_i)-6$, 
so for example one obtains moduli spaces of complex  dimension two when 
$$(m,r_1,r_2,\ldots) = (4,0,0,0,0), (3,1,0,0), (2,1,1), (2,2,0), (1,3)$$
as is well-known in the theory of Painlev\'e equations 
(these examples provided early motivation, cf. \cite{smid}). 
\end{rmk}

\ppb{
\begin{rmk}
The argument  in part 2) works more generally (again via \cite{drezet-slices} Proposition 5.7);
The closure $\bar \cC$ of $\cC\subset \bH$ is affine and thus so is 
$\mu^{-1}(\bar\cC)$.
Then the subset $\cU\subset \mu^{-1}(\bar\cC)$ of stable representations with stabiliser $\bK$ is open and saturated and is a principal $\bG$-bundle.
Further smooth points in $\cU$ map to smooth points of 
$\mu^{-1}(\bar\cC)/\bH$. 
\end{rmk}
}

\subsection{Irregular Deligne--Simpson problem}

Having defined and studied the notion of irreducible Stokes representations we can define the irregular analogue of the Deligne--Simpson problem in the present context.
Given an irregular curve $\Si$ 
with marked points $a_i$ and irregular types $Q_i$ 
(for $i=1,\ldots,m$) as above, 
 choose a conjugacy class
$$\cC_i\subset H_i$$
for each $i$ where $H_i=C_G(Q_i)$ as usual.
Let $\ga_i$
be the simple loop based at $b_i$ going once in a positive sense around the $i$th boundary component $\partial_i$ of $\wt \Si$ (as in SR2)). 

\noindent{\bf Question ($i$DS):} for which choices of conjugacy classes $\cC_i$
does there exist an {\em irreducible} Stokes representation
$$\rho\in \Hom_\IS(\Pi,G)$$
such that $\rho(\ga_i)\in \cC_i$ for each $i$? 

The original Deligne--Simpson problem is on the Riemann sphere with $G$ a general linear group and all the irregular types zero. 
We will make some conjectures in some irregular cases (again on the Riemann sphere with $G$ a general linear group) in \cite{cqv}.

The collection $(\cC_1,\ldots,\cC_m)$ of conjugacy classes is just a conjugacy class for the group $\bH$.
Let $\cC\subset \bH$ be the inverse conjugacy class, 
so a solution of the irregular Deligne--Simpson problem means that 
there are stable points in the subset
$$\mu^{-1}(\cC)\subset \Hom_\IS(\Pi,G)$$
since $\rho(\ga_i)=h_i$ is the inverse of the $H_i$ component of the moment map.
Thus the reduction $\mu^{-1}(\cC)^{stable}/\bH$ of the space of such stable points by $\bH$ is non-empty.
If $G$ is a general linear group then, via the irregular Riemann--Hilbert correspondence, such reductions are isomorphic to some of the hyperk\"ahler manifolds of \cite{wnabh} (specifically the cases here correspond to setting  $\Re(\la_i)=0$ in \cite{wnabh}, so that stability is equivalent to irreducibility).
Thus the irregular Deligne--Simpson problem translates into
the problem of characterising when certain hyperk\"ahler manifolds are nonempty.

\begin{rmk}
Kostov
has recently studied \cite{kostov-nfDSp}
an ``additive Deligne--Simpson problem for non-Fuchsian systems'', apparently suggested by Y. Haraoka.
This is not the additive analogue of our irregular Deligne--Simpson problem.\footnote{\cite{kostov-nfDSp} looks at the orbits of the residues rather than the orbits of the residues of the formal normal form.
The multiplicative analogue of the problem of \cite{kostov-nfDSp} would be to fix  conjugacy classes of the local monodromy around  singular points, rather than the conjugacy classes of the {\em formal} monodromy as we do here. These two notions coincide in the regular singular case.}
In fact some cases of the  additive analogue of our problem  have been studied (and solved) earlier in \cite{rsode} (see also \cite{iastalk07-nom, slims}).
As mentioned above the motivation for our version of the problem is from
the complex symplectic/hyperk\"ahler moduli spaces of \cite{smid, wnabh} appearing in the wild/irregular  extension of nonabelian Hodge theory. 
\end{rmk}

Part 2) of Corollary \ref{cor: good quotients} suggests that for
general groups a modified question ($i$DS$^+$) should also be considered:
For which choices of conjugacy classes $\cC_i$
is there a solution $\rho$ of $i$DS such that the stabiliser of 
$\rho$ in $\bH$ is minimal (i.e. equal to $\bK$)?

\section{Admissible deformations of irregular curves}\label{sn: adm def}

Given an irregular curve $\Si$ we have defined
a Poisson variety $\Hom_\IS(\Pi,G)/\bH$.
The aim of this section is to define the notion of an ``admissible deformation'' of $\Si$ over a base $\IB$, and show that the corresponding Poisson varieties fit together into the fibres of a fibre bundle with a canonical (complete) flat connection preserving the Poisson structures.
This leads to an algebraic Poisson action of the fundamental group of $\IB$ on the fibre $\Hom_\IS(\Pi,G)/\bH$.
We view this as the irregular analogue of the well-known mapping class/braid group actions on the character varieties.

Recall that, given a fixed connected complex reductive group $G$ with maximal torus $T$, 
an (algebraic) irregular curve is
a smooth compact algebraic curve $\Si$ with distinct marked points $a_1,\ldots, a_m\in \Si$ and an irregular type 
$Q_i$
at each marked point.
Now define a {\em family of irregular curves} over a (smooth) base $\IB$ to be a smooth family of curves
$$\pi: \Si\to \IB$$ 
so that each fibre $\Si_p=\pi^{-1}(p)$ is 
a curve (for $p\in \IB$),
with global sections $a_1,\ldots,a_m:\IB\to \Si$ (representing marked points of each fibre) and a smoothly varying family of irregular types $Q_i$.

\begin{defn}
An ``admissible deformation'' of an irregular curve $\Si_0$, consists of a family of irregular curves (with one fibre isomorphic to $\Si_0$),  such that 1) each fibre $\Si_p$  is smooth, 2) the marked points remain distinct, and 3) for any $i=1,\ldots,m$ and  any 
root $\al\in \cR$ the order of the pole of
$$\al \circ Q_i$$
does not change.
\end{defn}

Here $\cR\subset \lt^*$ denotes the roots of $G$ relative to $T$
and  if $\al$ is a root, $\al \circ Q_i$ is a germ of a  meromorphic function (well defined modulo holomorphic functions)  and its pole order is an integer $\ge 0$.
The simplest examples of admissible deformations  
of an irregular type $Q = A_r/z^r +\cdots $ were considered previously in \cite{JMU81} for $\GL_n$ and in \cite{bafi} for other $G$; they occur if all the terms are arbitrary except the leading
coefficient $A_r$ which is restricted to be regular, i.e. so that 
 $\al(A_r)\ne 0$ for all roots 
$\al$. 
In other words all deformations of the  coefficients are admissible, provided $A_r\in \lt$ stays off of all of the root hyperplanes.
Thus on a disk the space of such deformations is homotopy equivalent to 
$\lt_\reg=\{A\in \lt\st \al(A)\ne 0 \text{ for all roots $\al$}\}$, whose fundamental group is the (pure) $G$-braid group.
This brings the $G$-braid groups into play, much as deformations of curves with marked points involves mapping class groups and the usual Artin braid groups. 
For example one can show  (\cite{bafi} Theorem 3.6) that this gives the  geometric origins of the quantum Weyl group (which was constructed directly by Lusztig,  Soibelman and Kirillov--Reshetikhin by verifying explicit generators satisfied the desired relations).
Here we will consider general admissible deformations, so that the fundamental group of the space of admissible deformations will be more complicated.

Let $\pi: \Si\to \IB$ be an admissible family of irregular curves, and for any $p\in \IB$ let $M_p$ denote the Poisson variety $\Hom_\IS(\Pi,G)/\bH$ associated to the irregular curve $\Si_p$. 

\begin{thm} \label{thm: loc syst of vars}
The varieties $M_p$ assemble into a local system of Poisson varieties over $\IB$.
\end{thm}

This means that there is a fibre bundle $\pr:M\to \IB$ 
such that $\pr^{-1}(p) = M_p$ for any $p\in \IB$, and for any points $p,q\in \IB$ and path $\ga$ in $\IB$ from $p$ to $q$, there is a canonical algebraic Poisson isomorphism $M_{p}\cong M_q$, only dependent on the homotopy class of $\ga$.
Equivalently there is a covering $\cU=\{U_i\st i\in I\}$
of $\IB$ by contractible open sets (with contractible pairwise intersections), indexed by some set $I$, such that
if    $i\in I$ and $p,q\in U_i$ 
then there is a canonical Poisson isomorphism  
$\phi^{(i)}_{qp}:M_p\cong M_q$
and if $p\in U_i, q\in U_j, r\in U_i\cap U_j$
then  the isomorphism 
$$\phi^{(j)}_{qr}\circ\phi^{(i)}_{rp}:M_p\cong M_q$$
does not depend on the choice of $r$ in the two-fold intersection
(thus enabling us to define $M$ as a bundle with constant, Poisson, clutching maps).

\pf
The crucial point is that, modulo some direct spanning equivalences, one may locally on $\IB$ use the ``same'' generating paths for the fundamental groupoids $\Pi$, thereby identifying nearby fibres $M_b$ with the same explicit Poisson variety, built out of the fission spaces and internally fused doubles 
(and from Lemma \ref{lem: ds lem} we know the direct spanning equivalences induce Poisson isomorphisms).
This gives local trivialisations of $M$, over an open cover.
The clutching maps are constant Poisson isomorphisms since on each pairwise intersection they come from making a different choice of generating paths.

To make this precise we should first check carefully that using the ``same generating paths'' does indeed lead to direct spanning equivalences. This is local at a singularity, so we consider the unit disk $\Delta$ with an irregular type $Q$ at $0$, and choose generating paths as in Figure 
\ref{fig: Pi gens}. Then under a small admissible deformation the singular directions will move and may break up into more singular
directions. In effect the punctures along the singular directions may split into more punctures, but since everything is smooth, for sufficiently small deformations all the new punctures stay inside the chosen loops. We claim that the resulting Stokes groups (inside each loop) are direct spanning equivalent to the initial Stokes groups (and that no other truly new singular directions appear, e.g. outside the chosen loops). However this is clear from the definition of the Stokes groups, and of admissible deformations:
each root $\al$ supports $\deg(q_\al)$ punctures, and they vary continuously with $Q$. 

\begin{figure}[ht]
	\centering
	\input{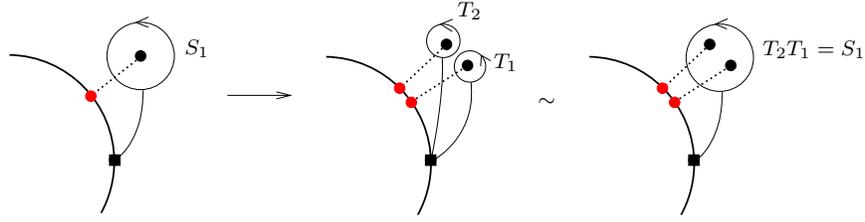}
\caption{Using the same loop, modulo direct spanning equivalence.}\label{fig: splitting}
\end{figure}

Finally to make sense of the notion of the ``same'' generating paths, note that in a neighbourhood of any point of $\IB$ we can choose local coordinates near each marked point on $\Si$ (so that in each such coordinate nothing is moving except the singular directions).
Then we may use these coordinates to make the punctures to define $\wt\Si$ (and as previously remarked, for different coordinate choices the resulting groupoids $\Pi$ may be canonically identified).
Then we just choose the neighbourhoods in $\IB$ small enough so the same loops may be used (up to direct spanning equivalence), as in Figure \ref{fig: splitting}. 
\epf

Thus roughly speaking the isomonodromy connection is defined by keeping the local products of Stokes multipliers constant (as in Jimbo--Miwa--Ueno \cite{JMU81}, and the extension to other groups in \cite{bafi}).
This situation is analogous to the action of the mapping class group of a Riemann surface on its spaces of fundamental group representations (via outer automorphisms of the fundamental group).
Indeed in the  nonsingular case ($m=0$) such results go back at least to Goldman \cite{Gol84}.
When the most singular coefficients of each irregular type are regular  semisimple
an analytic version of this was proved in \cite{smid} 
and an algebraic approach appears in \cite{Woodh-sat, Krich-imds, bafi} (in increasing generality).
The notion of ``local system of schemes'' was introduced by Simpson
\cite{Sim94b} \S6, where it was shown that 
moduli spaces of representations of fundamental groups of smooth projective varieties form local systems of schemes, when the varieties move in families.

\begin{rmk}
Note that set-theoretically the same argument shows that a path in $\IB$ from $p$ to $q$ yields a canonical bijection between the set of $\bH$ orbits in $\Hom_\IS(\Pi(p),G)$ and in 
$\Hom_\IS(\Pi(q),G)$
(not just at the level of $S$-equivalence classes)
where $\Pi(p)$ is the groupoid attached to $\Si_p$ 
(for any choice of basepoints) and similarly for $q$.
\end{rmk}

\subsection{Example}\label{ssn: example base}

If the underlying curve itself stays fixed then there is a stronger statement, since we may use the same basepoints.
For example suppose $G$ is semisimple and
let $Q_0$ be an irregular type at the origin of the unit disk $\Si=\Delta$, and let 
$\IB$ be the set of all admissible deformations of $Q_0$ (fixing the origin), so that 
\beq\label{eq: univ defmns}
\IB = \{ Q\in \lt\flp z \frp /\lt\flb z \frb \st 
\text{PoleOrder}(\al \circ Q) = \text{PoleOrder}(\al \circ Q_0)
\text{ for all roots $\al$}\}.
\eeq
Thus for example if $Q_0=A/z^r$ with $A$ regular semisimple, then $\IB$ is homotopy equivalent to $\lt_\reg$ and so the  fundamental group of $\IB$ is the pure $G$-braid group, for any $r\ge 1$.
Fix basepoints $b_1,b_2\in \partial \wh \Si$ as usual, and for any 
$p\in \IB$ let $\wt \Si_p\subset \wh \Si$ be the corresponding punctured surface, and let $\Pi_p= \Pi_1(\wt\Si_p,\{b_1,b_2\})$ be the corresponding fundamental groupoid.
The argument of the proof of Theorem \ref{thm: loc syst of vars} establishes the following (now at the quasi-Hamiltonian, rather than Poisson, level).

\begin{thm}\label{thm: local qh imd conn}  
The spaces $\Hom_\IS(\Pi_p,G)$ assemble into a local system of quasi-Hamiltonian spaces over $\IB$. 
\end{thm}

Explicitly, 
given two irregular types $Q,Q'\in \IB$ 
then each path in $\IB$ from $Q$ to $Q'$
determines an  algebraic isomorphism
$\Hom_\IS(\Pi_Q,G) \cong \Hom_\IS(\Pi_{Q'},G)$
relating the quasi-Hamiltonian structures and only depends on the homotopy class of the path.
This yields an (nonlinear) action of the fundamental group 
$\pi_1(\IB,Q)$ of the base  on the fibre 
$\Hom_\IS(\Pi_{Q},G)$.

For example if $H=C_G(Q)$ is abelian then the quotient $\Hom_\IS(\Pi_{Q},G)/G$ (just forgetting the framing on the outer boundary) is actually Poisson
and $\pi_1(\IB,Q)$ acts preserving the Poisson structure.
If further $Q=A/z$ then up to a covering $\Hom_\IS(\Pi_{Q},G)/G$ is isomorphic
\cite{smapg,bafi} to the Poisson Lie group dual to $G$; 
the computation of this nonlinear Poisson action appears in \cite{bafi}, and explicit formulae appear in \cite{k2p} (end of \S3) for $G=\GL_n(\IC)$---and, as explained there, in the case of $\GL_3(\IC)$ this action corresponds to the global monodromy of the Painlev\'e VI differential equation.

\begin{rmk}
Although the above example assumed $G$ was semisimple, in order to ensure the ``universal'' local deformation space 
\eqref{eq: univ defmns} was finite dimensional, this restriction is not necessary. In fact it is sometimes important not to make this restriction: for example in the case $G=\IC^*$ 
the infinite dimensional space of admissible deformations 
corresponds to the infinite number of ``times'' in many integrable hierarchies.
One can see this in several ways: 
1) the so-called Baker-Akhiezer functions are solutions of the corresponding irregular connections  (see e.g. \cite{schilling-bakerfns} and references therein), 
or 2) 
upon inverting the local coordinate $z$,  the element $e^Q$ 
is the same as the element $\exp(xz+t_2z^2+t_3z^3+\cdots)$ on p.9 of Segal--Wilson \cite{SW}.
\end{rmk}

\begin{rmk} (Full braid groups and bare curves.)
Since the Betti spaces do not dependent on the ordering of the marked points and only depend on the Weyl group orbit of each irregular type (cf. Remark \ref{rmk: finite quotients}) one can consider slightly  more general deformations, as follows. 
Define a ``bare irregular type'' to be an irregular type defined modulo the Weyl group $W$, i.e. an element $\bar Q\in (\lt(\wh\cK)/\lt(\wh\cO))/W$.
In turn a ``bare irregular curve'' is a curve $\Si$ with an unordered finite set of  marked points $S\subset \Si$ and a  bare irregular type $\bar Q_s$ at each point $s\in S$.
Thus any irregular curve $\Si$ has an underlying bare irregular curve $\bar \Si$ and the Betti moduli space $\bM_B(\Si)$ is determined by $\bar \Si$.
An admissible family of bare irregular curves then consists of 
a family $\pi :\Si\to \IB$ of smooth curves and a multisection, i.e. a subvariety $S\subset \Si$ finite \'etale over $\IB$, and a bare irregular type $\bar Q_s$ for each $s\in S\cap \pi^{-1}(b), b\in \IB$, such that locally over $\IB$ each $\bar Q_s$ is the $W$-orbit of an admissible family of irregular types $Q_s$.
Theorem \ref{thm: loc syst of vars} extends immediately to
show that the associated Betti spaces form a local system of Poisson varieties over 
$\IB$ for any admissible family of bare irregular curves (since any such family is locally isomorphic to an admissible family of irregular curves).
However the framed version Theorem \ref{thm: local qh imd conn} 
does not extend without some further choices (such as a pointwise lift of $W$ to $N(T)$) cf. \cite{bafi} Lemma 3.5 for an example.
  
\ppb{Suppose $\Si,\Si'$ are two irregular curves with marked points $S=\{a_i\},S'=\{a_i'\}$ and irregular types $Q_s,Q'_t$ (for each $s\in S,t\in S'$) respectively.
Define $\Si,\Si'$ to be {\em equivalent} if
there is an isomorphism $\varphi:\Si\to \Si'$ of underlying curves 
such that $\varphi(S)=S'$ and for each $s\in S$ there is some $w\in N_G(T)$ with  $\Ad_w(Q_s)=\varphi^*(Q'_{\varphi(s)})$. 
Similarly this extends to yield the definition of when an  admissible family $\Si\to \IB$ of irregular curves is {\em equivalent} to an  admissible family $\Si'\to\IB$.
Finally the notion of locally admissible family is obtained by gluing (over open subsets of the base) various admissible families to equivalent admissible families.

??? algebro-geometrically }
\end{rmk}

\appendix

\section{Stokes data from connections}\label{apx: stokes from connections}

This appendix summarises some results about the classification of meromorphic connections on curves. 
Most of these results are (well) known for $G$ a general linear group, and a path to extend them will be detailed elsewhere.
(The case when the leading coefficient is regular semisimple was established in \cite{bafi}, and one of the approaches there already involved multisummation.)
Note that, except to motivate the definitions, the results of this appendix are not used in the body of the article.
Section \ref{ssn: local classn} summarizes the ($G$-extension of the) local classification result, basically as it is presented in 
\cite{L-R94}.
Section \ref{ssn: local groupoid classn} explains how to relate this to 
Stokes representations and Stokes $G$-local systems.
Finally \S\ref{ssn: global stokes reps from cons} explains the global picture.

\subsection{Local classification}\label{ssn: local classn}
Fix a maximal torus $T\subset G$ and 
consider the closed unit disc  $\Delta$ in the complex plane with marked point $0$ and irregular type $Q$ at $0$.
If $z$ is a coordinate vanishing at $0$, consider the ring $\cO=\IC\{z\}$ of germs at $0$ of holomorphic functions on $\Delta$, 
its completion $\wh \cO = \IC\flb z\frb$,
and the corresponding fields of fractions $\cK\subset \wh \cK$.
Let $H=C_G(Q)\subset G$ be the stabiliser of $Q$,
so that $T\subset H$ and denote the Lie algebras 
$\lt\subset \lh\subset \g$.
Fix an element  $\La(z) \in \h(\cK)$ and thus a
connection
\beq\label{eq: formal type}
A^0 = dQ + \La(z) \frac{dz}{z} \in \h(\cK)dz
\eeq
on the trivial $H$-bundle over the disc.
Suppose $\La$ is such that the connection
$\La \frac{dz}{z}$ has a regular singularity at $0$.
In order to classify meromorphic connections formally equivalent to $A^0$ it  turns out to be simpler to first classify {\em marked pairs}, i.e. 
connections together with a choice of formal isomorphism with $A^0$.
To this end define
$$\cH(Q) = \left\{ (A,\wh F)\in \g(\cK)dz \times G(\wh \cK) \st
\wh F[A^0] = A \right\}/ G(\cK)$$
to be the set of isomorphism classes of marked pairs (here the Lie algebra valued meromorphic one-forms $A,A^0$ are viewed as connections on trivial bundles and the square brackets denote the gauge action).
Note that the formal transformations $\wh F$ appearing here cannot be arbitrary since they relate two convergent connections.
The main classification result may be stated as follows.

\begin{thm}\label{thm: local classn} $\cH(Q) \cong \ISto(Q)$, so that $\cH(Q)$ is isomorphic to a complex vector space.
\end{thm}

In more detail, such an isomorphism (taking a marked pair to its Stokes data) may be defined as follows. Let $\IA\subset S^1$ denote the set of singular direction of $Q$ at $0$. (A priori, without making the definitions more abstract, the isomorphism depends on the choice of a fundamental solution of $A^0$ near each singular direction $d\in \IA$.)

\begin{thm}\label{thm: multisums}
1) If $(A,\wh F)$ is a marked pair representing an element of $\cH(Q)$, then $\wh F$ is multisummable along each direction in $S^1\setminus \IA$.

2) Given $d\in\IA$ let $\Si^+_d(\wh F),\Si^-_d(\wh F)$ be the $G$-valued holomorphic maps obtained by multisumming $\wh F$ in a small sector on the positive (resp. negative) side of $d$, and let $\Psi_d$ be a fundamental solution of $A^0$ defined on a sectorial neighbourhood of $d$. Then, upon analytic continuation across $d$, both
$$
\Phi_d^+ := \Si^+_d(\wh F)\Psi_d,   \qquad\text{and}\qquad 
\Phi_d^- := \Si^-_d(\wh F)\Psi_d$$
are fundamental solutions of $A$, and moreover
$$\Phi_d^-=\Phi_d^+ S_d$$
for some ($z$-independent) element $S_d\in \ISto_d(Q)$.

3) Repeating for each $d\in \IA$ yields  a surjective map
 $$\left\{ (A,\wh F)\in \g(\cK)dz \times G(\wh \cK) \st
\wh F[A^0] = A \right\} \to \ISto(Q)$$
whose fibres are precisely the $G(\cK)$ orbits.
\end{thm}

\noindent(If such $\wh F$ depends holomorphically on some parameters then one can show that its multisums also vary holomorphically
and so $\ISto(Q)$ is an analytic moduli space for marked pairs.)
For $G$ a general linear group these results are known and are the result of work by many people, see especially \cite{BV89, L-R94, DMR-ci} and references therein such as \cite{sibuya77, jurkat78, BJLproper, Mal79,  BBRS91, MR91}. Note in particular that there is an alternative sheaf-theoretic description of the Stokes data due to Deligne \cite{DMR-ci}, although the above approach is more explicit.
For other groups one may adapt the above proofs. 
(Presumably one can also use a Tannakian approach, in essence considering homomorphisms from the wild fundamental group of \cite{MR91} in to $G$, but this seems to give less information---we really do want to work directly with $G$-valued Stokes data as in \cite{bafi} in order to better understand the isomonodromy deformations and resulting braid group actions.)

\subsection{Local groupoid representations} \label{ssn: local groupoid classn}
Now we will explain how Stokes representations arise
from connections,  
essentially rephrasing the above picture (to keep better track of the choices of fundamental solution of $A^0$).
Given the irregular curve $\Si=(\De,0,Q)$ as above, 
let $\wt \Si\subset \wh \Si\to \Si$ be the real blow up and the resulting punctured curve, as usual,
and let $\partial_1, \partial_2\subset \wt \Si$ denote the boundary circles, with $\partial_1$ lying over $0\in \Si$.
Draw a concentric circle (halo) $\IH$ through all the punctures in 
$\wt \Si$.
Choose a connection $A^0$ with irregular type $Q$ as in \eqref{eq: formal type}.
Now consider a meromorphic connection $A$ 
on the trivial $G$-bundle on $\Delta$ 
(singular only at $0$), 
together with a formal isomorphism $\wh F$
at $0$ between $A^0$ and $A$ (using the inclusion $H\subset G$ to view $A^0$ as a $G$-connection).

Given this data there is a canonically defined 
Stokes $G$-local system $\IL$ on $\wt \Si$ (as defined in Remark \ref{rmk: stokes local system}).
Namely inside $\IH$ (near $\partial_1$) 
we take the $H$-local system $\IL_0$ of solutions of $A^0$, 
and outside $\IH$ we take  $\IL$  to be 
the $G$-local system of solutions of $A$.
To glue them together on each component of $\IH$ (between two consecutive singular directions) we use the corresponding multisum of $\wh F$.
This gives the desired Stokes $G$-local system.

Now choose basepoints $b_i\in \partial_i$ for $i=1,2$,
and let $\Pi=\Pi_1(\wt \Si,\{b_1,b_2\})$ denote the corresponding fundamental groupoid.
If we 
choose a framing of $\IL_0$ at $b_1$ (equivalently this means choosing a 
fundamental solution of $A^0$ on a small cilium approaching $0$ in the direction $b_1$) and
a framing of $\IL$ at $b_2$, then taking the monodromy of $\IL$
with respect to these framings yields a 
Stokes representation $\rho\in \Hom_\IS(\Pi,G)$.

\subsection{Global picture}\label{ssn: global stokes reps from cons}

Now suppose $\Si$ is an arbitrary (algebraic) irregular curve, with irregular types $Q_i$ at marked points $a_1,\ldots,a_m$.
Let $\wt \Si\subset \wh \Si\to \Si$ be the real blow up and the resulting punctured curve  as usual, and draw halos $\IH_i$ on $\wt \Si$ through the punctures near $\partial_i$ for each $i$.
Let $A$ be a connection on an algebraic principal $G$-bundle $E^\circ$ on 
$\Si^\circ = \Si\setminus\{a_1,\ldots,a_m\}$.
We will say $A$ has irregular type $Q_i$
at $a_i$ if there is some extension $E$ of $E^\circ$ across $a_i$, and local trivialisation of $E$ in a neighbourhood of $a_i$  such that $A$ takes the form
\beq\label{eq: holom type}
dQ_i + \Ga(z) \frac{dz}{z}
\eeq
for some $\g$ valued map $\Ga$ (nonsingular at $z=0$),
where $z$ is a local coordinate vanishing at $a_i$. 
If this is the case it is possible 
pass to a new {\em formal} trivialisation of $E$ at $a_i$ in which the connection takes the form 
$A_i^0 = dQ_i+\La_i(z)dz/z$, as in \eqref{eq: formal type}, with $\La_i(z)$ a nonsingular $\h_i=\Lie(H_i)$ valued map (which may even be taken to be polynomial in $z$).
In contrast to \S\S\ref{ssn: local classn},\ref{ssn: local groupoid classn} above, here the element $\La_i$ in the normal form $A^0_i$ (and thus the formal monodromy) is not fixed a priori.

If we choose such formal trivialisations  $\wh F_i$ at each $a_i$
then there is a canonically determined Stokes $G$-local system $\IL$ 
on $\wt \Si$.
Namely we take the $H_i$-local system $\IL_i$ of solutions of $A^0_i$ inside $\IH_i$ and,  in the interior of $\wt \Si$ (outside all the halos) 
define $\IL$ to be the $G$-local system of solutions of $A$, and then we glue $\IL$ and $\IL_i$ as above using the multisums of $\wh F_i$ 
on each component  of 
$\IH_i$. 
This leads to the following equivalence of categories.

\begin{thm}
There is an equivalence between the category of connections on algebraic principal $G$-bundles on $\Si^\circ$ having irregular type $Q_i$ at $a_i$ ($i=1,\ldots,m$) and the category of 
Stokes $G$-local systems on $\wt \Si$.
\end{thm}

If we now choose basepoints $b_i\in \partial_i$ 
and define
the groupoid $\Pi=\Pi_1(\wt \Si,\{b_i\})$ as usual,
then
upon choosing a framing at $b_i$ of the local system $\IL_i$ (for each $i$),  a Stokes representation
$\rho \in \Hom_\IS(\Pi,G)$
is obtained, by
taking the monodromy of the Stokes $G$-local system $\IL$ via these framings.
Finally recalling that $\bH$ acts on $\Hom_\IS(\Pi,G)$ (and here we see this action corresponds to changing the choice of framings) the above result implies:
\begin{cor}
The isomorphism classes of such connections $(E^\circ,A)$
with irregular type $Q_i$ at each $a_i$, correspond bijectively to the $\bH$ orbits in $\Hom_\IS(\Pi,G)$. 
\end{cor}

\begin{rmk}
Various modifications of this  are useful for some applications 
(such as isomonodromy or wild non-abelian Hodge theory).
For example, rather than reducing at a conjugacy class $\cC\subset \bH$,
one may use one of the quasi-Hamiltonian $\bH$-spaces $\wh \cC$ of \cite{logahoric} Theorem B (weighted conjugacy classes) to obtain spaces of Betti data (filtered Stokes representations)  corresponding to meromorphic connections with unramified formal types on parahoric bundles, as in \cite{Sim-hboncc, yamakawa-mpa, logahoric} in the tame case.
The basic topological objects, {\em filtered Stokes $G$-local systems}, are defined by replacing $H_i$ by a weighted parabolic subgroup 
$P_i\subset H_i$ in the definition in Remark \ref{rmk: stokes local system}, generalising to the irregular case the filtered $G$-local systems of \cite{logahoric} Remark 2.
(Similarly in the general linear case $G=\GL(V)$
one may glue on to $\Hom_\IS(\Pi,G)$ some spaces 
$\cB(V_{0i},W_i)$ to obtain quasi-Hamiltonian spaces of Betti data for holonomic $\cD$-modules with unramified formal types, cf. \cite{DMR-ci} p.43 and \cite{malg-book} p.60---specifically
for each marked point $i=1,\ldots,m$ one takes $V_{0i}\subset V$ 
to be the kernel of $Q_i$, and $W_i$ to be arbitrary. The tame case of this, when each $V_{0i}=V$, gives the 
explicit description of the perverse sheaves on $\Si$ relative to $\{a_i\}$, cf. \cite{malg-book} p.34 and \cite{gmv-duke}.)
\ppb{
it is useful to have  a slightly more precise Riemann--Hilbert correspondence, involving filtered $G$ Stokes representations (analogously to \cite{Sim-hboncc, logahoric} in the tame case):
in the current set-up this amounts to adding the choice a parabolic subgroup of $\bH$ containing the formal monodromy $\bh=(h_1,\ldots,h_m)\in \bH$ (after reduction this corresponds to replacing $\cC\subset \bH$ by one of the spaces $\wh \cC$ of \cite{logahoric} (another quasi-Hamiltonian $\bH$ space).
On the other side this corresponds to adding a compatible parahoric level structure .
}
\end{rmk}

\ppb{
\sketch
Onto: construct by analytically gluing meromorphic connections on discs onto holomorphic connections on $\Si$ minus a disc around each $a_i$. On a disc around $a_i$ the result follows from Theorem \ref{thm: local classn} and we may glue to get a meromorphic connection on an analytic $G$ bundle $E$ over $\Si$.
This will necessarily be algebraic by GAGA (cf. \cite{serre-efa} Th\'eor\`eme 3), and we then restrict to $\Si^\circ$.

Injectivity.
Now suppose connections $A$ and $B$ (on bundles $E_A^\circ, E_B^\circ\to \Si^\circ$ say) lead to Stokes representations in the same $\bH$ orbit (via some choices of formal types $A_i^0, B_i^0$ and framings).
Then the formal monodromies are conjugate and this implies 
we may formally 
retrivialise so that $A_i^0=B_i^0$ for all $i$
(since regular singular connections on $H$-bundles are locally classified by the conjugacy class of their local monodromy).
In other words the algebraic bundles 
$E^\circ_A,E^\circ_B$ 
admit extensions $E_A,E_B\to \Si$ for which $A$ and $B$ are both isomorphic at $a_i$ to $A^0_i$, via a nonsingular formal isomorphism, for each $i$.
Then we may adjust the framings so the Stokes representations are actually equal. 
Then Theorem \ref{thm: local classn} implies $(E_A,A)$ and 
$(E_B,B)$ are holomorphically isomorphic in discs around each marked point, and in turn since the global monodromy is the same it follows they are globally holomorphically isomorphic over $\Si$.
Now GAGA implies  this is an algebraic isomorphism 
(\cite{serre-efa} \S6.3), which will restrict to $\Si^0$ to give the desired algebraic isomorphism
relating the connections.
\esketch

\subsection{Plan of proofs}\label{ssn: pop}
Here is a rough summary of steps to adapt existing (general linear) proofs
of Theorem \ref{thm: multisums}. 
In the case where the most singular term of $Q$ is regular such an adaptation
was already carried out carefully in \cite{bafi}---in fact multisummation was already used even in that case.

1) Establish a $G$-analogue of the asymptotic existence theorem for holomorphic fundamental solutions asymptotic to a given formal series solution on sufficiently small sectors. 
This can be established directly as in 
\cite{bafi} proof of Lemma 2.9 (cf. the usual case in \cite{Was76}), or by using multisummation of the matrix valued series obtained after choosing a faithful representation of $G$ (and using the fact that multisummation is a morphism of differential algebras to see that the multisum  will be a $G$-valued holomorphic map, as in \cite{bafi} Lemma A.5).

2) Use 1) to 
define the $G$-analogue of the Malgrange--Sibuya map 
$$\cH(Q)\to 
H^1(S^1, \Aut_1(A^0)),$$
 where $\Aut_1(A^0)$ is the sheaf of flat isotropies of $A^0$ (the sheaf of nonabelian groups on $S^1$ consisting of germs on small sectors with vertex  $0$ of holomorphic maps to $G$ which are automorphisms of $A^0$ and asymptotic to the identity).
Then check that Malgrange--Sibuya's  proof that this map is an isomorphism goes through verbatim (cf. \cite{bafi} Theorem A.7 for the extension of Malgrange's proof of surjectivity).

3) Rewrite the article \cite{L-R94} working with $G$-valued Stokes data, thus obtaining a proof that $H^1(S^1, \Aut_1(A^0))\cong \ISto(Q)$ and an algorithm to convert any cocycle representing a class in $H^1(S^1, \Aut_1(A^0))$ 
into a Stokes cocycle, and thus (after choosing fundamental solutions of $A^0$) an element of $\ISto(Q)$. 
(In turn one similarly obtains the existence of $G$-valued Ramis 
factorisations of the formal isomorphisms $\wh F$,
as the product of $G$-valued $k$-summable series for various $k$).
The key step is to see that Loday-Richaud's algorithm generalises: it consists of repeated use of various iterated semi-direct product decompositions, and these  decompositions generalise immediately, using the notion of direct spanning subgroups---indeed the iterated semi-direct decompositions are perhaps better rephrased as direct spanning decompositions.
More precisely the main steps of the algorithm are 

1) to factor 
the stalk  $\Aut_1(A^0)_\th$  (at $\th\in S^1$) in terms of levels, and 

2) to factor each fixed level subgroup $\Aut_1(A^0)^k_\th\subset \Aut_1(A^0)_\th$ as a product of level $k$ groups of Stokes cocycles (isomorphic to our groups $\ISto_d(k)\subset G$).

3) and then to use various commutation relations amongst the Stokes cocycles to order them by singular directions

The key point is that all of these steps go through, using the present $G$ analogue of the  definition of the Stokes groups and singular directions, and the proofs are similar to that of Lemma \ref{lem: stokes gps} using the notion of direct spanning subgroups.

}

\section{Stokes groups are groups}\label{sn: stokes gps}

We will prove Lemma \ref{lem: stokes gps}, which claimed that 
each  of the sets $\cR(d), \cR(d,k)$ is a closed subset of some system  of positive roots.

\pf
First we must show that if $\al,\be\in \cR(d)$ then any root of the form $\ga = n\al + m\be$ (for integers $n,m>0$) is in $\cR(d)$.
Thus $e^{q_\al(z)}$  and $e^{q_\be(z)}$ have maximal decay as $z\to 0$ in the direction $d$. Now $q_\ga = nq_\al + mq_\be$, so if the degrees of $q_\al$ and $q_\be$ are different then the result is clear. Otherwise suppose the degrees are both $k$ and the leading terms of $q_\al,q_\be$ are $c_\al/z^k, c_\be/z^k$ (resp.).
Then the maximal decay condition for $e^{q_\al}$ means that  $c_\al/z^k$ is real and negative when $\arg(z)=d$ (and similarly for $e^{q_\be}$). Thus the leading term $nc_\al/z^k+mc_\be/z^k$ of $q_\ga$ is again real and negative when $\arg(z)=d$, so $\ga\in \cR(d)$.
Clearly the same argument also works for  $\cR(d,k)$.
To show $\cR(d)$ (and thus also $\cR(d,k)$) is in some subset of positive roots we will find $\la\in \lt_\IR$ such that $\al(\la)>0$ for all $\al\in \cR(d)$. (By taking the derivative, the set of 1-parameter subgroups $\Hom(\IC^*,T)$ embeds as a lattice in $\lt$ and $\lt_\IR$ is defined to be its real span, so that $\lt = \lt_\IR\otimes \IC$, and $\cR$ is a subset of the real dual of $\lt_\IR$.)
Fix $z\in \IC^*$ with $\arg(z)=d$ and for each $i$ let 
$R_i = -\Re(A_i/z^{k_i})\in \lt_\IR$,
where $Q=\sum A_i/z^{k_i}$.
Thus $\al(R_i)>0$ if $\al \in \cR(d,k_i)$ and 
$\al(R_i)=0$ if $\al \in \cR(d,k_j)$ with $j<i$.
 Now set $\la_r = R_r$, and $\la_{r-1}= N\la_r + R_{r-1}$ for a large real number $N$. Since $\cR$ is finite we can choose $N$ large enough such that 
$\al(\la_{r-1})>0$  for any $\al \in \cR(d,k_r)\cup \cR(d,k_{r-1})$.
Iterating (with $\la_{i-1}= N\la_i + R_{i-1}$ for various $N$) yields
$\la_1$ with $\al(\la_{1})>0$  for any $\al \in \cR(d)$, as required. 
(Note that we can move $\la_1$ off all the root hyperplanes by going one step further: setting $\la_0 = N\la_1 + R_0$ for any regular $R_0\in \lt_\IR\setminus \bigcup_{\al\in \cR} \ker \al$.)
\epf

\renewcommand{\baselinestretch}{1}              %
\normalsize
\bibliographystyle{amsplain}    \label{biby}
\bibliography{../thesis/syr} 

\vspace{0.5cm}   
DMA \'Ecole Normale Sup\'erieure and CNRS, 
45 rue d'Ulm, 
75005 Paris, 
France

www.math.ens.fr/$\sim$boalch

boalch@dma.ens.fr

\end{document}